\documentclass{svjour3}
\usepackage{geometry}
\usepackage{epsfig}
\usepackage{amsmath,bm}
\usepackage{amsfonts}
\usepackage{amssymb}
\usepackage{cuted}
\usepackage{caption}
\usepackage{subfiles}
\usepackage[numbers,sort&compress]{natbib}
\usepackage{pgfplots}
\usepackage{filecontents}
\usepackage{tikz}
\usepackage{tikzscale}
\usepackage{pgf}
\usepackage[normalem]{ulem}
\usetikzlibrary{patterns}
\usetikzlibrary{calc}
\usetikzlibrary{decorations.pathmorphing,decorations.markings}
\RequirePackage{fix-cm}
\definecolor{Color_FS}{rgb}{.7,0,.7}
\definecolor{Color_NF}{rgb}{0,0,.7}
\definecolor{Color_MD}{rgb}{.75,.35,0}
\definecolor{Color_SMD}{rgb}{1,.7,0}

\newcommand{\red}[1]{{\color{black!87!black} {#1}}}

\newcommand{\ot}[1]{{\color{black!55!black} {#1}}}

\newcommand{\N}{\mathbb N}

\newcommand{\avec}{\mathbf{a}}
\newcommand{\bvec}{\mathbf{b}}
\newcommand{\cvec}{\mathbf{c}}
\newcommand{\uvec}{\mathbf{u}}
\newcommand{\fvec}{\mathbf{f}}
\newcommand{\mvec}{\mathbf{m}}
\newcommand{\Avec}{\mathbf{A}}
\newcommand{\Bvec}{\mathbf{B}}
\newcommand{\Evec}{\mathbf{E}}
\newcommand{\Fvec}{\mathbf{F}}
\newcommand{\Gvec}{\mathbf{G}}
\newcommand{\Hvec}{\mathbf{H}}
\newcommand{\Mvec}{\mathbf{M}}
\newcommand{\Nvec}{\mathbf{N}}
\newcommand{\Pvec}{\mathbf{P}}
\newcommand{\Rvec}{\mathbf{R}}
\newcommand{\Kvec}{\mathbf{K}}
\newcommand{\Xvec}{\mathbf{X}}
\newcommand{\Yvec}{\mathbf{Y}}
\newcommand{\Wvec}{\mathbf{W}}

\newcommand{\xvec}{\mathbf{x}}
\newcommand{\yvec}{\mathbf{y}}
\newcommand{\zvec}{\mathbf{z}}
\newcommand{\x}{\mathbf{x}}

\newcommand{\kvec}{\mathbf{k}}

\newcommand{\svec}{\mathbf{s}}
\newcommand{\X}{\mathbf{X}}

\newcommand{\R}{\mathbf{R}}
\renewcommand{\S}{\mathbf{S}}
\newcommand{\M}{\mathbf{M}}
\newcommand{\K}{\mathbf{K}}
\newcommand{\Lvec}{\mathbf L}

\newcommand{\G}{\mathbf{G}}
\renewcommand{\H}{\mathbf{H}}
\newcommand{\g}{\mathbf{g}}
\newcommand{\h}{\mathbf{h}}
\newcommand{\phit}{\boldsymbol \phi}
\newcommand{\V}{\mathbf{V}}
\newcommand{\I}{\mathbf{I}}
\newcommand{\Ot}{{\boldsymbol \Omega}^2}
\newcommand{\at}{{\mathbf{a}}}
\newcommand{\bt}{{\mathbf{b}}}
\newcommand{\ct}{{\boldsymbol {\gamma}}}

\newcommand{\At}{{\mathbf{A}}}

\newcommand{\Bt}{{\mathbf{B}}}

\newcommand{\alphavec}{\boldsymbol \alpha}
\newcommand{\betavec}{\boldsymbol \beta}
\newcommand{\gammavec}{\boldsymbol \gamma}
\newcommand{\phivec}{\boldsymbol \phi}

\newcommand{\Psivec}{\boldsymbol \Psi}
\newcommand{\Upsvec}{\boldsymbol \Upsilon}
\newcommand{\etavec}{\boldsymbol \eta}
\newcommand{\Omegavec}{\boldsymbol \Omega}

\newcommand{\xivec}{\boldsymbol \xi}
\newcommand{\Lambdavec}{\boldsymbol \Lambda}
\newcommand{\lambdavec}{\boldsymbol \lambda}

\newcommand{\tp}[1]{{#1^{\operatorname{T}}}}
\newcommand{\dd}{{\operatorname{d}}}
\newcommand{\DD}{{\operatorname{D}}}

\newcommand{\diver}{{\mathbf{\operatorname{div}}}}
\newcommand{\grad}{{\mathbf\nabla}}

\DeclareMathOperator{\expo}{e}

\newcommand{\zervec}{\mathbf{0}}

\begin{document}
\title{Model order reduction methods for geometrically nonlinear structures:
a review of nonlinear techniques}
\titlerunning{Model order reduction methods for geometrically nonlinear structures: a review} 
\author{	Cyril Touz\'e\and
Alessandra Vizzaccaro\and
Olivier Thomas
}

\institute{C. Touz\'e	\at
	IMSIA, ENSTA Paris, CNRS, EDF, CEA, Institut Polytechnique de Paris, 828 Boulevard des Mar{\'e}chaux, 91762 Palaiseau Cedex, France\\
	\email{cyril.touze@ensta-paris.fr}
	\and
	A. Vizzaccaro \at
	Department of Mechanical Engineering, University of Bristol, Bristol, BS8 1TR, UK\\
	\email{alessandra.vizzaccaro@bristol.ac.uk}
\and
O. Thomas  \at
Arts et M{\'e}tiers Institute of Technology, LISPEN, HESAM Universit{\'e}, F-59000 Lille, France\\
\email{olivier.thomas@ensam.eu}
}

\date{Received: date / Accepted: date}

\maketitle

\begin{abstract}
This paper aims at reviewing nonlinear methods for model order reduction of structures with geometric nonlinearity, with a special emphasis \red{on} the techniques based on invariant manifold theory. Nonlinear methods differ from linear based techniques by their use of a nonlinear mapping instead of adding new vectors to enlarge the projection basis. Invariant manifolds have been first introduced in vibration theory within the context of nonlinear normal modes (NNMs) and have been initially computed from the modal basis, using either a graph representation or a normal form approach to compute mappings and reduced dynamics.  These developments are first recalled following a historical perspective, where the main applications were first oriented toward structural models that can be expressed thanks to partial differential equations (PDE). They are then replaced in the more general context of the parametrisation of invariant manifold that allows unifying the approaches. Then the specific case of structures discretized with the finite element method is addressed. Implicit condensation, giving rise to a projection onto a stress manifold, and modal derivatives, used in the framework of the quadratic manifold, are first reviewed.  Finally, recent developments allowing direct computation of reduced-order models (ROMs) relying on invariant manifolds theory  are detailed. Applicative examples are shown and the extension of the methods to deal with further complications are reviewed. Finally, open problems and future directions are highlighted.

\keywords{Reduced order modeling \and geometric nonlinearity \and thin structures \and invariant manifold \and nonlinear mapping}
\end{abstract}

\section{Introduction}

The scope of the present contribution is the nonlinear dynamics exhibited by elastic structures subjected to large amplitude vibrations, such that geometric nonlinearities are excited. The focus is set on the derivation of efficient, predictive and simulation-free reduced-order models (ROM). 

Geometric nonlinearities are associated \red{with} large amplitude vibrations of thin structures such as beams, plates and shells, because of their relatively low bending stiffness. By \red{its} nature, it is a distributed nonlinearity, which means that all the degrees of freedom of the model are nonlinearly coupled. On the contrary, other types of nonlinearities, such as those related to contact, are associated \red{with} localized nonlinearities. In this latter case, reduction methods are often associated \red{with} substructuring (see {\em e.g.} \cite{Substruct,krack2017,lai2019,corradi2021}), which cannot be transposed to the present case of geometric nonlinearities. Applications to real-world engineering problems are numerous and tend to increase \red{since} lightweight thinner structures are \red{being increasingly} used. The range of applications thus spans from aeronautics to transportation industry~\cite{PESHECK2001Blade,PesheckASME2002,PATIL2004,KerschenAircraft13,MIAN2014,thomas16-ND,quaegebeur2020,martin2019}, wind energy systems~\cite{Manolas2015,Dou2020}, musical instruments~\cite{monteil15-AA,jossic18,ChaosGong} and micro/nanoelectromechanical systems (M/NEMS), in which those nonlinearities must be mastered to design efficient structures~\cite{FRANGI2019,LazarusThomas2012,sajadi2019,massimino2020,zega2020}. The vibratory phenomena arising from the nonlinear dynamics of geometrically nonlinear structures can also be intentionally used for the purpose of new designs, especially in the M/NEMS domain or for energy harvesting, where for example internal resonances or parametric driving are conceived with specific goals~\cite{Younisbook,shoshani2021,qalandar2014,thomas13-APL,li2017,mam2017,nitzan2015,iurasov2020,fallahpasand2019}.

With regard to the aim of deriving effective ROM, geometric nonlinearity presents two important difficulties. The first one is connected to the nonlinear dynamics itself and the number of different solutions arising when nonlinearity comes into play. Instabilities, bifurcations, important changes in the nature of the solutions, \red{the} emergence of more complex dynamics including quasiperiodic solutions, chaotic solutions and even wave turbulence in structural vibrations, have been observed experimentally and numerically studied with models (see {\it e.g} \cite{Nayfeh79,Nayfeh00,mangussi2016,guillot2020MAN} for examples of dynamics exhibited by \red{oscillating} systems and \cite{Amabilibook,Lacarbo,thomas07-ND,WTbook,nayfehcarboChin99} for nonlinear phenomena in beams, plates and shells). The second issue is connected to the distributed nature of the nonlinearity and the resulting nonlinear coupling \red{that} appears in the equations of motion. Of course, these two characteristics are the two faces of the same coin since the couplings are the most important drivers of the complexity observed in dynamical solutions. Nonetheless, while the first problem concerns analysis and treatment of complex dynamical phenomena often observed when nonlinearities are present,  the second is more directly related to model order reduction, which needs specific methods for geometrical nonlinearities,  often alleviated to an  efficient choice of a reduction basis that could take into account these couplings.

Most of the model order reduction methods can be seen as linear methods, where the aim is to find the best orthogonal basis to represent the dynamics, and add new basis vectors until convergence. In this setting, the main problem is to have a computational method allowing one to automatically compute the basis vectors. Linear modes have been used for a long time for such problems for their ease of computation and their clear physical meaning~\cite{NICKELL1976,Nayfeh79,Meirovitch80,MeirovitchF}. However, their main drawback is the number of nonlinear couplings. In a finite element context, it imposes reduction bases with a prohibitive number of modes to reach convergence, most of them having natural frequencies out of the frequency band of interest~\cite{givois2019,Vizza3d}. Those drawbacks can be compensated \red{for} with the addition of extra vectors in the basis such as modal derivatives~\cite{IDELSOHN1985,IDELSOHN1985b} or dual modes~\cite{KIM2013}.  Proper orthogonal decomposition~\cite{HolmesLumley,BerkoozPOD93}, arising from statistical methods (Karhunen-Loeve decomposition)~\cite{Karhunen,Kosambi} and having a direct link with the singular value decomposition~\cite{LIANG2002}, have also been used with success in numerous applications related to nonlinear vibrations~\cite{Krysl01,LinLeeLu02,AmabiliPOD1,Kerschen2005POD,SampaioSoize,Goncalves08}. The major drawback of this strategy is the need of preliminary data to compute the basis vectors, often obtained by time integration of a full scale model. More recently, proper generalized decomposition (PGD), aiming at generalizing the POD approach~\cite{Chinesta2011,metoui2018}, \red{has} also been used in a context of nonlinear vibration problems~\cite{GroletPGD,MEYRAND2019}. Since the topic of this review article is focused on {\it nonlinear reduction methods}, all these linear methods will thus not be covered (or only cited for illustrative purposes), and the reader is referred to already existing reviews on these methods for further information~\cite{Chinesta2011,mignolet13,LU2019,Substruct}.

The focus of this paper is to review the reduction methods that are essentially nonlinear, in the sense that they are based on  defining new variables that are nonlinearly related to the initial ones, instead of producing a linear change of basis as in the above-cited techniques. Being nonlinear, they also associate the ROM with a curved structure in phase space: a {\em manifold}. In this realm,  particular subsets are of main importance. The {\em invariant} manifolds of a dynamical system are indeed  particularly suitable domains on which reducing the dynamics. By definition, it is a region of the phase space which is invariant under the action of the dynamical system. In other words, any trajectory  of the system  that is initiated in the invariant manifold is entirely contained in it for all time.
Hence the invariance property ensures that the trajectories of the ROM are also trajectories of the full system, which is a very important prerequisite to define accurate reduction techniques. If this is not fulfilled, then the meaning of the simulation produced by the ROM with regard to the full system will remain unclear. 
Moreover, the curvature of the invariant manifolds in the phase space directly  embeds the non-resonant couplings and thus represents in a single object the slave coupled modes, without the need of adding new vector basis to catch these couplings, nor finding a correct computational method for sorting them out. In short, using invariant manifolds transform\red{s} the question of reduction to a problem of geometry in the phase space. 

The perspective of this review paper is thus strongly related to the application of invariant manifold based techniques for model order reduction of nonlinear vibratory systems. A special emphasis is also put on methods applicable to finite element problems. Since the vast majority of engineering calculations are nowadays performed using a finite element (FE) procedure to discretize \red{the spacial geometry of the structure}, numerous reduction techniques have tackled the problem of geometric nonlinearity with special adaptation to comply with FE formulation. An important feature arising from this viewpoint led to \red{a distinction between} intrusive and non-intrusive methods. By non-intrusive, it is meant that the reduction method can take advantage of the basic capabilities of any FE code, without the need to enter at the elementary level to perform specific calculations. In practice, a standard FE code is used with its already existing features, \red{which} are then specifically post-processed to build a ROM. On the other hand, intrusive method\red{s} compute the needed quantity at the elementary level, such that an open (or in-house) code is needed. 

The paper is organized as follows. Section~\ref{sec:framework} details the starting point and the typical equations of motion one has to deal with when geometric nonlinearity is taken into account. A short review of models is \red{given} in \red{Section}~\ref{sec:EOM} and the general form of the equations of motion that will be used in the rest of the paper is \red{also} given. Section~\ref{sec:which} gives a rapid survey of the most classical nonlinear phenomena and their consequence in the correct choice of a ROM. Section~\ref{sec:IMfull}  reviews the derivation of ROMs for nonlinear vibratory systems expressed in the modal basis. It starts with a short survey of the underlying mathematical developments in \red{Section}~\ref{subsec:dynasys}. Graph style and normal form styles are then reviewed in Sections~\ref{sec:IMSP} and~\ref{sec:NFmodal}, and Spectral Submanifolds in \red{Section}~\ref{subsec:SSM}. Section~\ref{sec:ROMFE} then discusses the application to FE models. The Stiffness Evaluation Procedure (STEP) is first \red{reviewed} in \red{Section}~\ref{subsec:STEP}, then implicit condensation is detailed in \red{Section}~\ref{subsec:IC}. The construction of a quadratic manifold with modal derivatives is reported in \red{Section}~\ref{subsec:MDQM}. These last two methods produce different manifolds that are compared to the invariant ones. Finally, direct computations of invariant manifolds from the physical space and thus directly applicable to FE discretization, \red{are} explained in Sections~\ref{sec:DNF} and~\ref{sec:directSSM}. The paper closes with a discussion on open problems and future directions in \red{Section}~\ref{sec:future}.

\section{Framework}\label{sec:framework}

This Section is devoted to \red{delineate} the framework of the problems addressed in this review. First, the different kinds of models used to tackle large amplitude vibrations of thin structures with geometric nonlinearity are surveyed. Some simplified analytical models, obtained thanks to a selected number of assumptions and  allowing the derivation of partial differential equations (PDE), are first recalled, to introduce the main physical consequences of geometric nonlinearity. Then, basic features of finite element  modeling and its specificities are exposed. Finally, the most common types of dynamical solutions exhibited by such systems are described, the analysis of which being necessary to better ascertain the ROM needed.

\subsection{Equations of motion}\label{sec:EOM}

\begin{figure}
\begin{center}
\includegraphics[width=.7\textwidth]{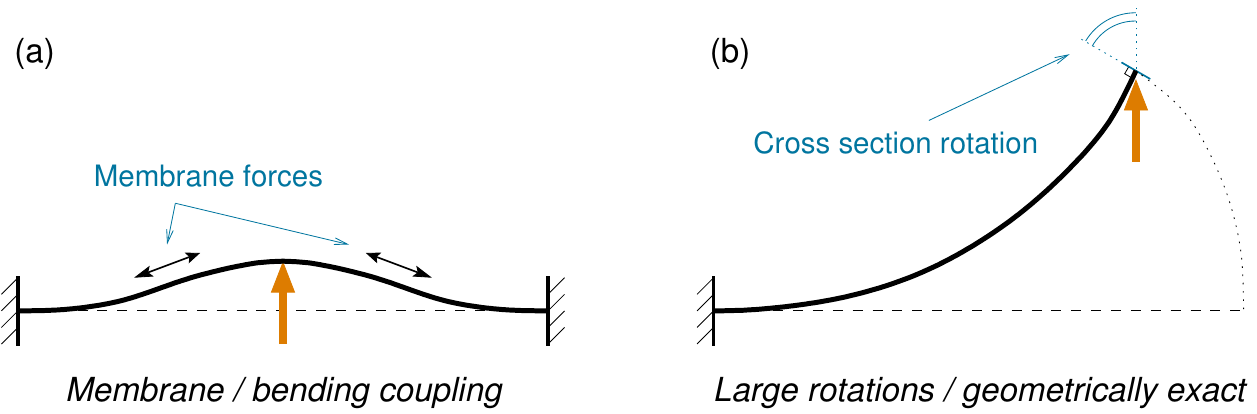}
\caption{Illustration of the two main families of models for geometrically nonlinear thin structures.}
\label{fig:vonkarGR}
\end{center}
\end{figure}

Geometric nonlinearity is \red{a} consequence of \red{a} large amplitude \red{change of geometry}, beyond the small motion assumptions ensuring the validity of a linearized model. For thin structures, they \red{become evident} once the transverse vibration amplitude is of the order of the structure's thickness~\cite{Nayfeh79,Amabilibook,neukirch2021}. Other nonlinearities can also be observed at large amplitude, such as material nonlinearities (plasticity, large deformations of soft materials, nonlinear piezoelectricity\ldots). In this paper, only a linear elastic constitutive law for the material is  considered, valid for small strains, to focus on geometrical nonlinearities. In practice, this situation is often encountered for thin structures, for which the small thickness allows large transverse displacements with small strains. 
Importantly, the nonlinearity is polynomial as a function of the unknowns (generally displacements and velocities), again in contrast to contact and friction that involve \red{steeper functional forms and may be modeled with non-smooth terms}. 
\red{The aim of this section is twofold. First, it will be shown that for common structural models (beams, plates, shells, three-dimensional continuous media of arbitrary shape), either analytical or discretized by a finite-element method, quadratic and cubic nonlinear terms are sufficient to describe the geometric nonlinearity. Second, the focus is set on the main physical and structural mechanisms that give birth to geometric nonlinearities, and which involve either membrane/bending coupling or large rotations in thin structures. The objective is not to provide a complete discussion on approximate beam/plate models, which can be seen themselves as a reduction technique, but rather on delivering simple keys to the vibration analyst to understand the source of the nonlinearities and also help him to choose a correct model.}

\subsubsection{Analytical PDE models}
\label{sec:PDE}

A number of models have been derived for beams, arches, plates and shells, based on simplifying assumptions, and only a few, \red{representative}, are recalled here, mainly to survey the related physical phenomena. Most thin structure models are based on the assumption that the cross sections are subjected to a rigid body kinematics. Then, depending on the range of amplitude of the rotation of the cross sections, two main families of models are considered, as illustrated in Fig.~\ref{fig:vonkarGR}. 

The first one can be denoted as the ``von~K\'arm\'an'' family of models. It is based on a  clever truncation of the membrane strains (the only nonlinearities kept in the strains expressions are  quadratic terms in the rotation angles of the cross section) due to Theodore von~K\'arm\'an \cite{vonkarman10}. This assumption, directly linked to a truncation of the rotations of the cross section, enables one to write very simple analytical (and numerical) models, that have been used in a large number of contributions. For a straight beam of length $l$ with homogeneous cross section, \red{the governing PDE} reads \cite{woinowsky50,Bazant,Nayfeh79}:
\begin{equation}
\label{eq:VKbeam}
\rho S \ddot{w} +EI w'''' - N w'' = p,\quad N=\frac{EI}{2l}\int_0^l w'^2\,\dd y,
\end{equation}
and for a plate, one has \cite{Landau1986,Bazant,ThomasBilbao08,Audoly,BilbaoVK15}:
\begin{equation}
\label{eq:VKplate}
\rho h \ddot{w} + D\Delta\Delta w - L(w,F) = p,\quad \Delta\Delta F=-\frac{Eh}{2}L(w,w).
\end{equation}
In those models, $w(y,t)$ (resp. $w(\mathbf{y},t)$) is the transverse displacement at time $t$ and location $y$ in the middle line of the beam (resp. location $\mathbf{y}$ in the middle plane of the plate), $\dot{w}=\partial w/\partial t$, $w'=\partial w/\partial y$, $\Delta$ is the bidimensional Laplacian, $L(\circ,\circ)$ is a bilinear operator, $(E,\rho)$ are the density and Young's modulus of the material, $(S,I)$ the area and second moment of area of the cross section of the beam, $h$ the thickness of the plate, $D=Eh^3/12(1-\nu^2)$ its bending stiffness and $p(y,t)$ ($p(\mathbf{y},t)$) an external force per unit length (resp. area). The axial (resp. membrane) inertia is neglected, enabling \red{one} to obtain a uniform axial force $N(t)$ in the beam and a scalar Airy stress function $F(\mathbf{y},t)$ to represent the membrane strains in the plate. Their main characteristics \red{are} that they accurately model the {\it axial/longitudinal} (resp. {\it membrane/bending}) coupling, that is the first physical source of geometric nonlinearities, illustrated in Fig.~\ref{fig:vonkarGR}(a). Indeed, when the structure is subjected to a transverse displacement $w$, its length (for the beam) or the metric of its \red{neutral} surface (for the plate) is modified, thus creating axial/membrane stresses proportional to the square of $w$, which thus nonlinearly increase the bending stiffness. This can be seen in Eqs.~(\ref{eq:VKbeam}) and (\ref{eq:VKplate}), in which $N$ and $F$ are quadratically coupled to $w$, which in turn creates cubic terms in the equations of motion, classically related to a {\it hardening behaviour}. 

If the initial geometry of the structure shows a curvature, Eq.~(\ref{eq:VKplate}) can be modified to the following shallow shell equation~\cite{Ostiguy92,Amabiliplatetheoexpe,CamierCTOT09}:
\begin{equation}
\label{eq:VKshell}
\rho h \ddot{w} + D\Delta\Delta w - L(w_0,F) - L(w,F) = p,\quad \Delta\Delta F + Eh L(w_0,w) = -\frac{Eh}{2}L(w,w),
\end{equation}
in which the geometry of the middle surface is represented by the static deflection $w_0(\mathbf{x})$ (that has to be small to ensure the validity of the model). Compared to the plate equation~(\ref{eq:VKplate}), two additional terms appear. They are responsible \red{for} a linear membrane / bending coupling, but also \red{for} a quadratic nonlinear coupling, both directly linked to the non-flat geometry of the shell (a non zero $w_0$)~\cite{thomas04}. This quadratic coupling can be responsible \red{for} a {\it softening behaviour} of the vibration modes \cite{Amabiliplatetheoexpe,TOUZE2004CS,touze-shelltypeNL,touze08-MPE,jossic18}.

The second family of models, illustrated in Fig.~\ref{fig:vonkarGR}(b) and usually known as ``geometrically exact models'', is more refined since no simplifying assumption on the modelling of the spatial rotation of the cross sections is formulated. The writing of those models as partial differential equations (PDE) is explicit only for simple geometries and their solving often relies on numerical discretization techniques like FE (see {\it e.g.} \cite{simo1989-shell} for shells and \cite{geradin2001} for beams) or finite-difference \cite{lang2012}. Because of the untruncated rotation operator, the nonlinearities appear in the PDEs in terms of sine and cosine functions\footnote{In fact, the parametrisation of the rotations in geometrically exact models can take several forms (full rotation matrix, quaternions, Lie groups\ldots \cite{cottanceau17-FEAD,geradin2001}), mainly to avoid singularities in the case of very large rotations (of several turns). The minimal sine/cosine parametrisation discussed here allows \red{the description of} the planar motion of a beam, but not a full 3D motion. It is shown here only to formally understand the nature of the nonlinear couplings.} of the cross section rotations (see the case of a straight cantilever beam in {\it e.g.} \cite{thomas16-ND}). In the case of a cantilever beam, one can obtain a very accurate and widely used model due to Crespo da~Silva \& Glynn \cite{crespo78-1}, which reads \cite{thomas16-ND}:
\begin{equation}
\label{eq:GRbeam}
\rho S \ddot{w} +EI w'''' + \underbrace{EI(w'w''^2+w'^2w''')}_\text{curvature} +\underbrace{\frac{\rho S}{2}\left[w'\int_l^x\left(\int_0^x w'^2\dd y\right)^{..}\dd y\right]'}_\text{axial inertia} = p.
\end{equation}
It is obtained from the geometrically exact model of the cantilever beam by (i) using an inextensibility condition to condense the axial motion into the transverse equation of motion and (ii) truncating Taylor expansion of the trigonometric functions of the cross-section rotation to the third order. To this end, it is interesting to compare Eqs.~(\ref{eq:GRbeam}) and (\ref{eq:VKbeam}). In the case of a cantilever beam, the axial force is null ($N=0$) because of the free end boundary condition, and the von~K\'arm\'an model~(\ref{eq:VKbeam}) becomes linear. On the contrary, the large rotation model~(\ref{eq:GRbeam}) makes appear two higher order nonlinear terms, related (i) to the large rotation of the cross section (curvature term) and (ii) to its axial inertia, condensed in the transverse motion. 

The conclusion is that geometric nonlinearity can be created by different physical effects. In the first family (Fig.~\ref{fig:vonkarGR}(a)), it comes from a membrane/bending coupling, which is effective only if the structure is constrained in the axial direction. For a 1D structure, this effect is observed only if the ends are axially restrained~\cite{lacarbonara2006}, thus explaining  why the von~K\'arm\'an model of a cantilever beam is linear. For plates and shells, the validity of the model depends on the deformed shape  during vibrations. Most of the time, the deformation changes the metric of the middle surface (most of mode shapes of a plate/shell are not developable surfaces) and the von~K\'arm\'an model can be thus used safely since offering accurate predictions.  On the contrary, the second family (Fig.~\ref{fig:vonkarGR}(b)) of models is mandatorily needed if the rotations of the cross-sections are large (from several tens of degrees to several turns). They make appear higher order stiffness couplings as well as nonlinear terms due to the inertia (with time derivatives). 

Another important point is the hardening/softening effects, the latter being created either by a loss of transverse symmetry of the geometry of the structure in its transverse direction (due to curvature and/or a nonsymmetric laminated structure \cite{thomas13-APL}) in the case of a von~K\'arm\'an model, or because of inertia effects like in Eq.~(\ref{eq:GRbeam}) (the first mode of a cantilever beam is hardening, whereas the others are softening \cite{NayfehPai}).

A key feature of the models described above is that, thanks to given assumptions, they are able to provide the equations of motion under the form of a PDE. In essence, they are limited to simple geometries, due to the fact that analytical models for producing PDE need to rely on a specific coordinate system. This limitation is generally relaxed for the shape of a shell model, since $w_0(\mathbf{y})$ can be chosen arbitrary in Eq.~(\ref{eq:VKshell}). But even in this case, the shape of the imperfect plate needs to be rectangular or circular to coincide with a simple coordinate system. All these models also clearly underline the nature of the geometric nonlinearity, which is distributed and involves only quadratic and cubic nonlinear terms as a function of the displacement. Finally, separating the models within the two families underlined above helps in  understanding   numerical simulations, in particular in term\red{s} of hardening / softening behaviour, in relation \red{to} membrane/bending coupling, curvature or inertia nonlinearities.

\subsubsection{Finite elements and space discretization}

\red{Most of the engineering applications now use space discretization based on the finite element (FE) procedure~\cite{Bathe}, mostly because of the geometry of the structural elements that can be more easily accounted for}. \red{From a} modeling point of view, this choice has for main consequence that one cannot rely anymore on a PDE to derive mathematical tools for reduced-order modeling. \red{Nowadays, a number of codes are available so that one can easily perform standard operations such as computing the eigenvalues and eigenvectors of a vibratory problem}. Since all these codes are routinely used for engineering applications, the notion of non-intrusiveness has emerged as an important feature for deriving reduced-order models.

FE discretisation techniques can be applied to all classical PDEs of mechanical models (after a proper variational formulation) and in particular to the nonlinear beam, plate and shell models discussed in section~\ref{sec:PDE}, for which 1D/beam or 2D/plate FE are used. It is also possible to avoid the thin structure cross-section kinematical constraint and to use 3D FE. In this case, the framework considered in this article is a full Lagrangian formulation with a Green-Lagrange strain measure  $\mathbf{E}$, conjugated with the second Piola-Kirchhoff stress measure $\mathbf{S}$, for which the strong form of the problem reads \cite{holzapfel00,Touze:compmech:2014,givois21-CS}:
\begin{equation}
\label{eq:3DMMC}
\diver(\mathbf{F}\mathbf{S})+\mathbf{f}_\text{b}=\rho\ddot{\mathbf{u}},\quad \mathbf{S}=\mathbf{\mathbb{C}}\mathbf{E},\quad \mathbf{E}=\frac{1}{2}\left(\grad\mathbf{u}+\tp{\grad}\mathbf{u}+\tp{\grad}\mathbf{u}\grad\mathbf{u}\right),
\end{equation}
where $\mathbf{u}(\mathbf{y},t)$ is the displacement field at point $\mathbf{y}$ of the 3D domain occupied by the structure. The first of the above equation is the \red{equation of motion}, in which $\mathbf{f}_\text{b}(\mathbf{y},t)$ is an external body force field and $\mathbf{F}$ the deformation gradient; the second equation is the linear elastic constitutive law, with $\mathbf{\mathbb{C}}$ the four-dimensional elasticity tensor; the third equation is the definition of the strain $\mathbf{E}$. The operators $\diver$ and $\grad$ are the vector divergence and the tensor gradient\footnote{For any tensor field $\mathbf{A}$, the $i$-th cartesian component of its divergence is $\sum_j  \partial A_{ij}/\partial y_j$ where $y_i$ is the $i$-th component of the position vector; for any vector field $\mathbf{v}$, the $(i,j)$ cartesian component of its tensor gradient is $\partial v_i/\partial y_j$.}. Since $\mathbf{F}=\mathbf{I}+\grad\mathbf{u}$ (with $\I$ the identity tensor), formally eliminating $\mathbf{F}$, $\mathbf{S}$ and $\mathbf{E}$ as a function of $\mathbf{u}$ in the equilibrium equation leads to obtain an equation of motion with a polynomial stiffness operator with quadratic and cubic nonlinear terms. The scope of this formulation is generic, exact (no assumption on the kinematics of the continuous media have been formulated) and theoretically embeds the thin structure models of section~\ref{sec:PDE}.

The starting equations for this contribution is the semi-discretised equations of motion: discretised in space and continuous in time. If one starts from a PDE (like those of section~\ref{sec:PDE}), then the space discretisation can be obtained using any Rayleigh-Ritz or Galerkin procedure or any other method that fits to the problem at hand, including a FE procedure. Another choice could be a 3D FE discretization of Eqs.~(\ref{eq:3DMMC}). In all cases, all unknowns resulting from the space discretization procedure are   gathereed in the displacement vector $\X(t)$. In case of a PDE and {\em e.g.} a Rayleigh-Ritz method, $\X$ contains all the unknown generalised coordinates related to the shape functions used to discretise the problem. In case of a FE procedure, $\X$ gathers all the degrees of freedom (dofs) of the model (displacements/rotations at each nodes). Denoting by $N$ the size of $\X$, the semi-discretised equations of motion for our geometrically nonlinear problem reads:
\begin{equation}
\M\ddot{\X}+\K \X+\fvec_\text{nl}(\X)=\fvec_\text{e},
\label{eq:eom_phys}
\end{equation}
where $\M$ is the mass matrix, $\K$ the tangent stiffness matrix at the equilibrium of the structure and $\fvec_\text{e}(t)$ is a vector of external forcing. In a general framework, $\fvec_\text{e}$ may also depend on the displacement vector $\X$ (an example of which being follower forces), however this case is not considered here for the sake of simplicity. In the present framework, it is first assumed that the geometric nonlinearities give rise only to quadratic and cubic polynomial terms involving only the displacement vector $\X$ (see the models of section~\ref{sec:PDE} and Eqs.~(\ref{eq:3DMMC})). They are expressed in the following internal force vector:
\begin{equation}
\fvec_\text{nl}(\X)=\G(\X,\X)+\H(\X,\X,\X),
\end{equation}
thanks to the terms $\G(\X,\X)$ and $\H(\X,\X,\X)$, using a functional notation for the quadratic and cubic terms with coefficients gathered in third-order tensor\footnote{The term ``tensor'' used here simply refers to a multidimensional array of dimension larger than two, and not to a multilinear map, as used in  mechanical models of stress and strain for instance.} $\G$ and fourth-order tensor $\H$. Their explicit indicial expressions read:
\begin{subequations}
\label{eq:tensor_product_phys}
\begin{align}
\G(\X,\X)&=\sum^N_{r=1}\sum^N_{s=1}\G_{rs}X_r X_s,\\
\H(\X,\X,\X)&=\sum^N_{r=1}\sum^N_{s=1}\sum^N_{t=1}\H_{rst} X_r X_s X_t,
\end{align}
\end{subequations}
where $\G_{rs},\H_{rst}$ are the $N$-dimensional vectors of coefficients $G^l_{rs}$, $H^l_{rst}$, for $l=1,\, ...,\, N$. In practice, the components of $\G$ and $\H$ are rarely computed, since it would lead to a huge computational burden and memory requirement for large values of $N$ ($\G$ has $N^3$ components while $\H$ has $N^4$). In standard FE codes, $\fvec_\text{nl}(\X)$ is computed by standard assembly procedure.

Simple extensions of this framework could include systems with polynomial terms involving the velocities, arising in different fields. For ease of reading, these cases are not considered in detail, but will be emphasised when needed. Note that all the methods explained \red{henceforth} are extendable to handle such cases.

\subsection{Modal expansion}

Eq.~\eqref{eq:eom_phys} expresses the semi-discretized equations of motion in physical space. For all the next developments, the equations in the modal space needs to be defined.  Let $(\omega_p,\phit_p)$ be the $p$-th eigenfrequency and eigenvector of the linearized problem, \red{which satisfy}:
\begin{equation}
(\K - \omega_p^2 \M)\phit_p = \zervec.
\label{linear_eigsys}
\end{equation}
Using normalization with respect to mass, one has
\begin{equation}
\V^T\M\V=\I,\qquad\mbox{and}\qquad\V^T\K\V=\Ot,
\label{eq:diagonalisation}
\end{equation}
with $\V$ the matrix of all eigenvectors, $\V = [\phit_1, ..., \phit_N]$, $\I$ the identity matrix, and $\Ot$ a diagonal matrix composed of the square of the natural frequencies, $\Ot = \mbox{diag}(\omega_p^2)$.
The linear change of coordinate $\X = \V \x$ is used to go from the physical to the modal space, where $\x$ is the $N$-dimensional full vector of modal displacements. The dynamics reads:
\begin{equation}
\ddot{\x}+\Ot\x+\g(\x,\x)+\h(\x,\x,\x) = \zervec,
\label{eq:eom_modal}
\end{equation}
where the third- and fourth-order tensors $\g$ and $\h$ expresses the nonlinear modal coupling coefficients. They are linked to their equivalent $\G$ and $\H$ in the physical basis via:
\begin{subequations}
\begin{align}
&\g_{ij}
=\V^T
\G(\phit_i, \phit_j),
\label{eq:g_ij}\\
&\h_{ijk}
=\V^T
\H(\phit_i, \phit_j, \phit_k).
\label{eq:h_ijk}
\end{align}
\end{subequations}
The modal equations of motion can be detailed line by line, $\forall \; p =1,\ldots N$:
\begin{equation}
\ddot{x}_p + \omega_p^2 x_p +  \sum_{i=1}^{N} \sum_{j=1}^{N} g_{ij}^p x_i x_j + \sum_{i=1}^{N} \sum_{j=1}^{N} \sum_{k=1}^{N}  h_{ijk}^p x_ix_jx_k = 0.
\label{eq:EDOmodalproj}
\end{equation}
It can be noticed that the above equation is not written with the upper-triangular form\footnote{In Eq.~(\ref{eq:EDOmodalproj}), some coefficients can be grouped together since being related to the same monomial: $g_{ij}^px_ix_j+g_{ji}^px_jx_i=\hat{g_{ij}^p}x_ix_j$ with $\hat{g}_{ij}^p=g_{ij}^p+g_{ji}^p$, $j>i$, one of the upper triangular coefficient, see appendix~\ref{sec:ghsym}.}  of the tensors $\g$ and $\h$, which is often used due to the commutative property of the usual product (see {\it e.g.} \cite{touze03-NNM,givois2019}). As explained in appendix~\ref{sec:ghsym}, we assume in this contribution that the internal force vector $\kvec(\X)=\K\X+\fvec_\text{nl}(\X)$ derives from a potential energy, which leads to particular  symmetry relationships in the nonlinear quadratic and cubic coefficients. These symmetry relationships are different, depending on the fact that the upper-triangular form is adopted or not. Appendix~\ref{sec:ghsym} recalls all these relationships in a unified manner.

In the above modal expansion, the maximum number of modes $N$ has been formally retained since in the present paper, Eq.~(\ref{eq:EDOmodalproj}) is not used for computational purpose. This point will be addressed in section~\ref{subsec:STEP}. The number of nonlinear coupling terms (scaling as $N^4$) being a very large number, it is important to understand the different roles played by the monomials. Appendix~\ref{app:class} recalls the terminology used to classify these terms, that will be used throughout the paper. Among them, some play a very important role for understanding the idea of {\em invariance} that is key for the computation of invariant manifolds.  Let us assume that $m$ is the main mode having most of the vibrational energy (for example in the case of a harmonic forcing in the vicinity of $\omega_m$). Then all terms $g^p_{mm} x_m^2$ and $h^p_{mmm} x_m^3$ for all other equations labelled $p$ are {\em invariant-breaking} terms. Indeed, the sole presence of these terms creates a coupling that breaks the invariance of the linear eigensubspace~\cite{touze03-NNM,TouzeCISM}, and thus feeds energy to the other linear modes that cannot be easily neglected. Tracking those specific terms will thus be of importance in all the next derivations. 

The second classification criterion is linked to the fact that the nonlinear terms can be interpreted as a forcing on \red{the} $p$-th oscillator. This interpretation leads to the definition of {\em resonant}  and {\em non-resonant} monomials. For a given monomial, its resonant (or non resonant) nature depends on the oscillator to which it belongs, its order and also to eventual internal resonances between the oscillation frequencies of the oscillators. Resonant terms have a major importance in the mode couplings and the related exchanges of energy. This is more detailed in Appendix~\ref{app:class} and in the remainder of the text.

\subsection{Which ROM for which dynamics ?}\label{sec:which}

The choice of a ROM capable of producing accurate predictive results for a structure with geometric nonlinearity must rely on a correct analysis and understanding of the dynamics one wants to reproduce or predict with the model. Since nonlinearities are present, the behaviour of the system is amplitude-dependent. A correct two-dimensional parameter space to classify possible dynamics and advise \red{on} the choice of a ROM should include the {\em frequency content} of the forcing and the {\em vibration amplitude} of the structure. Depending on this vibration amplitude, very complex phenomena can be observed and the analysts should have a clear idea \red{of} the type of dynamical solutions \red{they} want to simulate with \red{the} ROM. Thus, the nature of the ROM will also condition the type of dynamical solutions one wants to represent.

In the rest of the paper, we will denote as ``master coordinates" the variables kept in order to describe the dynamics of the reduced model, and ``slave coordinates" all the others. Of course, one looks for ROM strategies in which the number of master coordinates is as small as possible. In the case where the vibration amplitude is moderate so that the system is close to linear vibrations and weakly nonlinear, the number of master coordinates needed to describe the dynamics should follow the same rules as in the linear case. This means that the number of master modes must be nearly equal to the number of eigenfrequencies contained in the frequency band of the forcing. In particular, a good ROM should account for the non-resonant couplings existing between the linear modes, even if they are out of the frequency band of interest, and embed them in the reduction process. An example of this is the membrane/bending coupling in thin structures, discussed in section~\ref{sec:PDE}, for which the low frequency bending modes are nonlinearly coupled to high frequency axial modes, the latter being sometimes very far from the frequency band of interest \cite{givois2019}. In the case of 3D FE models, some similar couplings occur with very high frequency thickness modes, as investigated in \cite{Vizza3d}. 

Consequently, an accurate ROM should contain only the driven transverse modes and enslave the axial/thickness motions directly in a transparent and automatic manner, such that the analyst does not need to derive \red{a} cumbersome convergence study to verify the accuracy of the reduction. This is one of the propert\red{ies} of the invariant manifold-based approach, thus making them particularly appealing. As long as the vibration amplitude is moderate, then the number of master modes can be selected according to the frequency band of the forcing. If the forcing is harmonic with moderate amplitude, then reduction to a single master mode should be targeted in order to describe the backbone curve. If a band-limited noise excitation drives the structure, then the number of selected master coordinates should be equal to the number of modes in the excitation frequency band.

This simple picture is however complicated by the presence of resonant interactions between the modes, which are linked to the existence of internal resonance relationships between the eigenfrequencies of the structure. A second-order internal resonance is a relationship of the form $\omega_p \pm \omega_k = \omega_j$ between three eigenfrequencies of the structure, which can degenerate in the simple 1:2 internal resonance when one has $\omega_j= 2\omega_l$. These second-order internal resonances are directly connected to the quadratic terms of the nonlinear restoring force, and can be linked to three-waves interactions in the field of nonlinear waves~\cite{zakbook,nazabook,during06,DURING201742}. Third-order internal resonance involves commensurability between four eigenfrequencies $\omega_p \pm \omega_k \pm \omega_l= \omega_j$ and are also related to four-waves interactions. This includes the simplest case of 1:1 internal resonance when $\omega_j= \omega_l$ as well as 1:3 resonance when $\omega_j= 3\omega_l$, and is connected to cubic nonlinearity. When such internal resonances exist, resonant couplings occur and strong energy exchange may take place (see Appendix~\ref{app:class} for more details). In such a case, a ROM should then retain as master, all the modes whose eigenfrequencies \red{possess} internal resonance relationships with the directly driven ones, since peculiar couplings leading to bifurcations can appear. This complicates a little the analysis but fortunately, the first analysis of internal resonance can be done on the eigenfrequencies which are generally known.

Moving to larger amplitudes, the picture complexifies again with the appearance of internal resonance between the {\em nonlinear} frequencies of the system. Indeed, since the oscillation frequencies depend on amplitude, an internal resonance relationship can be fulfilled at moderate to large amplitudes, \red{with the response frequencies of the system}. This is more difficult to predict beforehand since it can be analysed only by computing the backbone curves of each mode and verify that no strong internal resonance can be fulfilled at larger amplitudes.

\begin{figure}[ht]
\begin{center}
\includegraphics[height=.26\linewidth]{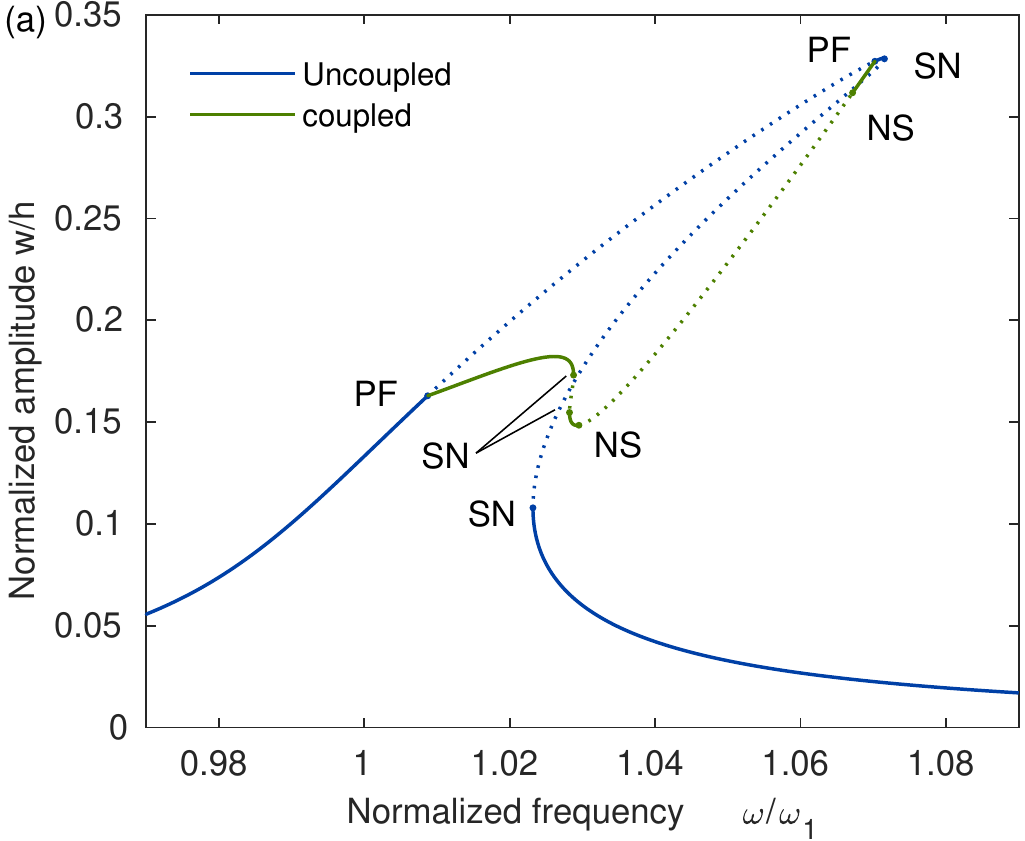}\hspace{.6em}
\includegraphics[height=.26\linewidth]{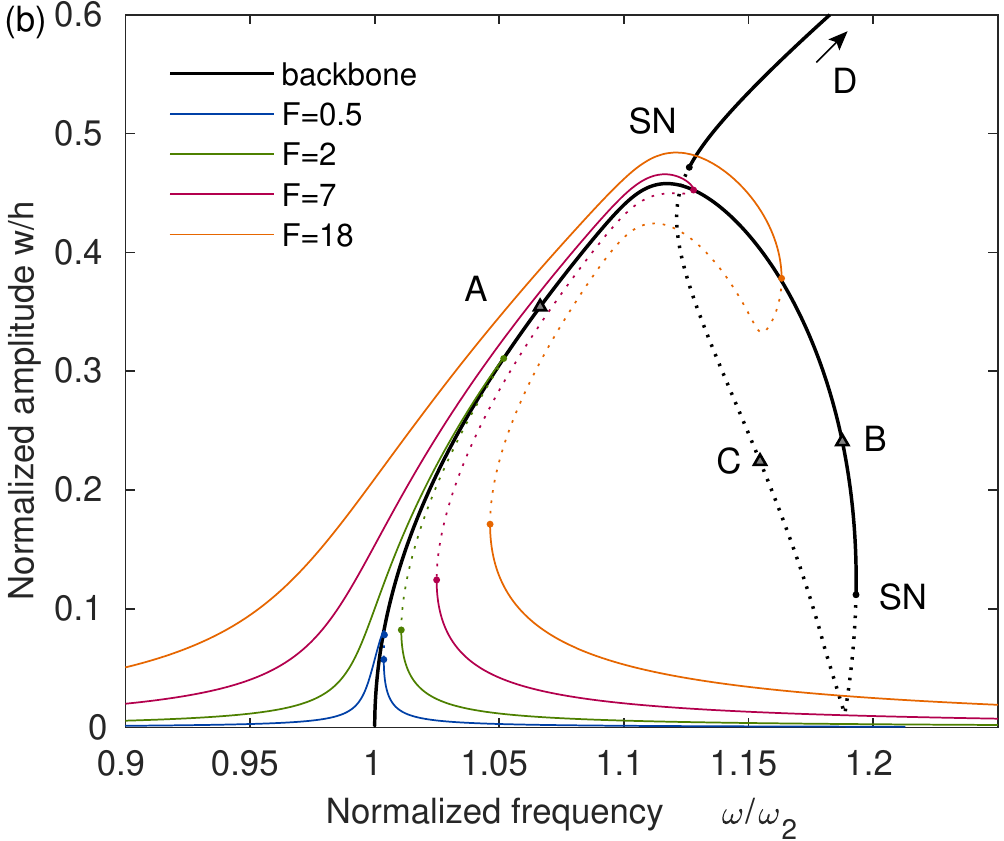}\hspace{.6em}
\includegraphics[height=.26\linewidth]{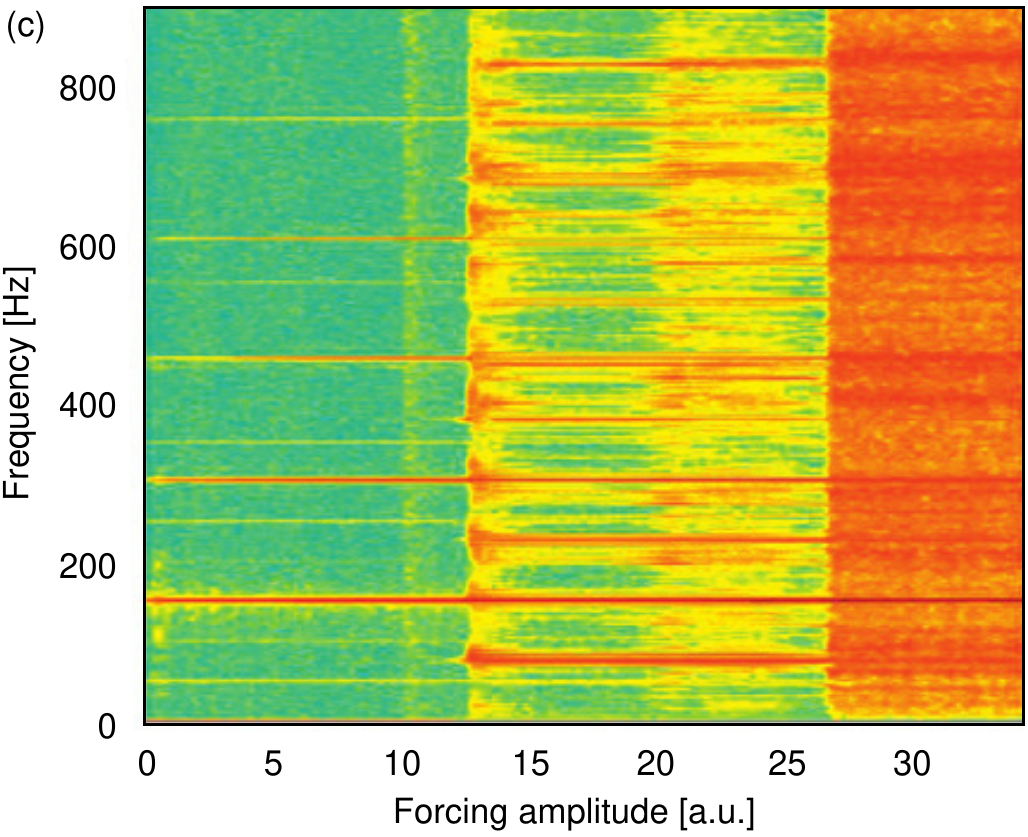}
\end{center}
\caption{(a) Frequency response of a clamped-clamped beam excited in the vicinity of its first bending mode, in 1:1 resonance with the companion mode in the other bending direction~\cite{YichangVib}. Maximum amplitude over one period of the directly excited mode at the driving point (at $0.275$ times the length of the beam from one end), scaled by the thickness of the beam. Horizontal axis scaled by the first eigenfrequency. The forced response \red{shows the existence of} bifurcation points, typical of 1:1 resonance: pitchfork (PF), saddle-node (SN) and Neimark-Sacker (NS). '---': stable branches; '$\cdot\cdot\cdot$': unstable branches. (b) Backbone curve and forced responses, at various amplitudes, of the second bending mode of a clamped-clamped beam, showing the characteristic loop due to activation of 1:3 internal resonance between the nonlinear frequencies of modes 2 and 4~\cite{givois2019,artDNF2020}. Same vertical axis as (a),  the horizontal axis being scaled by the second natural frequency. (c) Experimental spectrogram of the vibration response of a rectangular plate harmonically forced with frequency 151 Hz and increasing amplitude, showing transitions from periodic solutions to wave turbulence~\cite{Touze:JSV:2012}. \red{Points A, B, C and D in (b) refers to Fig.~\ref{fig:manifolds} and are used subsequently.}}
\label{fig:resonanceexample}
\end{figure}

These two cases are illustrated in Fig.~\ref{fig:resonanceexample}(a-b) with a clamped-clamped beam. Fig.~\ref{fig:resonanceexample}(a) shows the frequency response curve of such a beam that is allowed to vibrate out of plane, in both transverse directions and having a square cross-section. Consequently, the two fundamental bending modes in each polarization have the same eigenfrequency and the structure naturally \red{possesses} a 1:1 resonance. The beam is excited with a force in only one direction. Out of the resonance, only the driven, directly excited mode, participates to the vibration, its companion staying quiescent (blue curve).  A pitchfork bifurcation (PF) gives rise to a coupled solution where both modes are vibrating (green curve).  Along this coupled branch, two Neimark-Sacker bifurcations are observed(NS), from which a quasiperiodic regime emerges. Two saddle-node (SN) bifurcations also exist, as it is the case for \red{an equivalent} single dof nonlinear oscillator. This example shows that a simple system composed of only two master modes in 1:1 resonance can already display very different dynamical solutions. It also underlines that the minimal ROM should contain two master modes.

\begin{figure}[ht]
\begin{center}
\includegraphics[width=.8\textwidth]{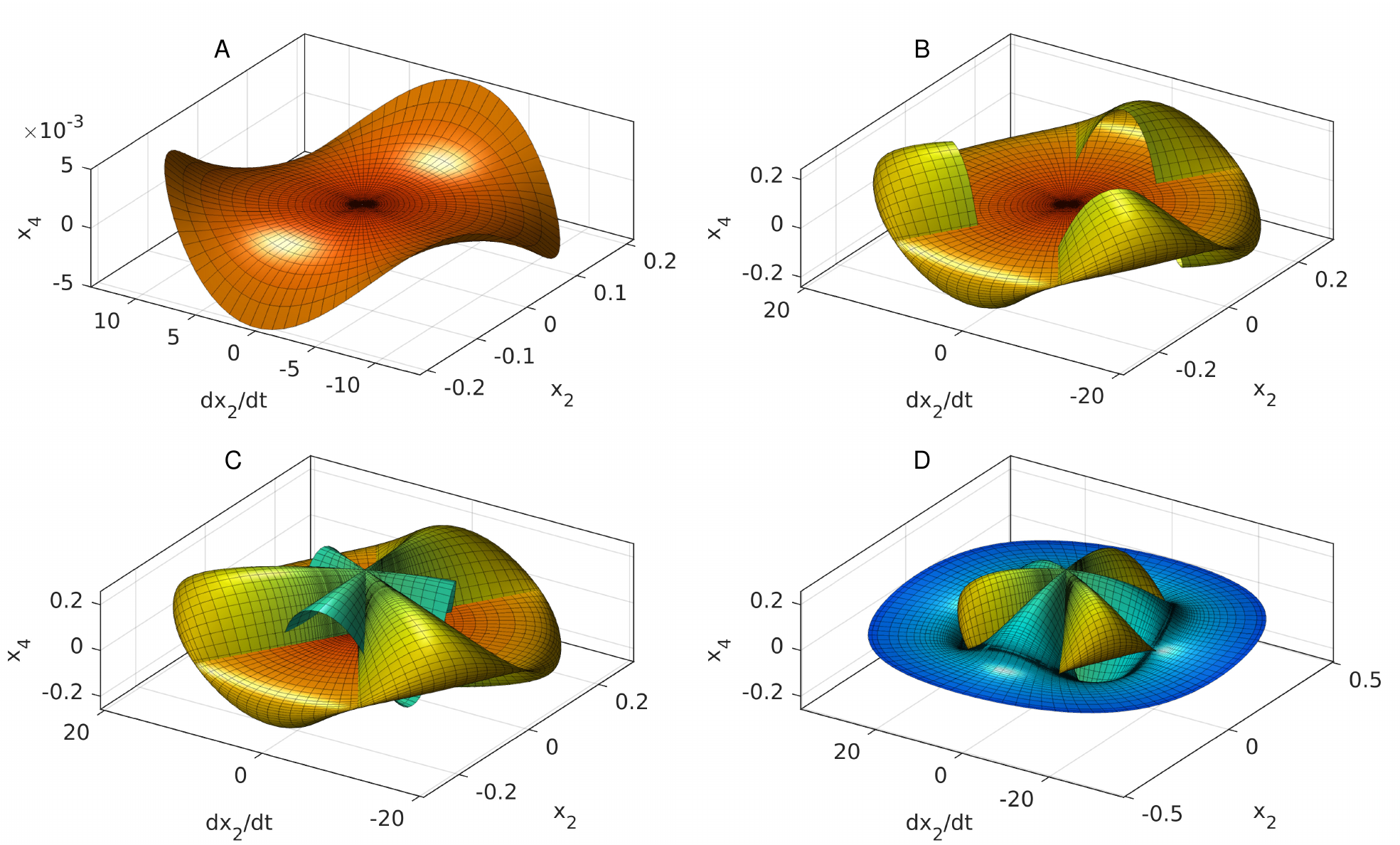}
\end{center}
\caption{3-d representations of the invariant manifold  (LSM) associated to the backbone curve of Fig.~\ref{fig:resonanceexample}(b), in the subspace spanned by the modal coordinates $(x_2,\dot{x}_2,x_4)$, showing a 1:3 internal resonance between the 2nd and 4th modes of a clamped/clamped beam. The four views of the developing manifold correspond to points $A$, $B$, $C$ and $D$ indicated in Fig.~\ref{fig:resonanceexample}(b), plotted by assembling the periodic orbits for increasing arclength of the numerical continuation with the software Manlab \cite{Guillot2019,guillot2020MAN}, and showing an apparent folding in this 3-d representation.}
\label{fig:manifolds}
\end{figure}

A second example is shown in Fig.~\ref{fig:resonanceexample}(b), where the backbone curve of the second mode of a straight clamped-clamped beam is plotted. Its cross section is chosen without symmetries to avoid a 1:1 internal resonance. For small amplitudes, a hardening behaviour is observed, and this could be reported by a ROM having a single master mode. However, for larger amplitudes, a loop appears in the solution branch, denoting a strong interaction and the emergence of an internal resonance. It is also responsible of a folding of the corresponding invariant manifold, as shown in Fig.~\ref{fig:manifolds} and discussed in Section~\ref{sec:IMfull}. What is interesting in this case is that the resonance relationship occurs between the nonlinear frequencies, whereas the natural frequencies were not close enough to fulfill the resonance relationship. In this particular case, a 1:3 resonance occurs with mode 4, creating a strong interaction. A correct ROM should then include mode~4 as additional master coordinates to fully recover the coupling. This means in particular that the choice of the master modes is made difficult and is strongly amplitude-dependent, since possible internal resonance between nonlinear frequencies could appear. Consequently, the simple analysis of the relationships between the natural frequencies might not be enough.

As mentioned \red{at} the beginning of this section, the parameter space allowing one to get a rough idea of the possible dynamics should include the frequency content but also the vibration amplitude. This amplitude dependence, already \red{addressed} above concerning Fig.~\ref{fig:resonanceexample}(b), can also be illustrated by inspecting how a thin structure bifurcates to complex regimes when it is forced harmonically with increasing amplitudes. Following numerous experiments and numerical simulations reported in~\cite{ChaosGong,TOUZE:IJNLM:2011,Touze:JSV:2012}, a general scenario for the transition can be observed. It is illustrated in Fig.~\ref{fig:resonanceexample}(c) reporting an experimental measurement on a plate, harmonically excited at 151~Hz. For small vibration amplitudes, the regime is weakly nonlinear, and only harmonics of the solution appear in the response. The ROM targeted for reproducing such a dynamics should contain one master mode. A first bifurcation occurs where the spectrum of the vibration response is enriched by a number of extra peaks. The appearing peaks correspond to internally resonant modes, such that the energy is now spread between all the modes that are strongly coupled to the driven one. In the case reported in the figure, only one mode appears through a 1:2 internal resonance, with eigenfrequency at 75~Hz. For this range of vibration amplitude, a ROM containing only the internally resonant modes must be enough to reproduce this dynamics. At larger amplitudes, a second bifurcation occurs, leading to a more complex regime characterized by a broadband Fourier spectrum. This regime is typical of {\em wave turbulence}. Wave turbulence has been studied in a number of physical contexts and the interested reader is referred to~\cite{zakbook,Newellreview,nazabook} for a complete view \red{of} the theory and its applications. Application to plate vibrations have been investigated since the pioneering work by D{\"u}ring, Josserand and Rica~\cite{during06,DURING201742}, including numerous experimental and numerical studies, see {\em e.g.}~\cite{boud2008,mordant2008,R21,Humbert,DucceschiPhysD,Yoko:PRL2013,miquel14} as well as the review chapter~\cite{WTbook}. In this dynamical regime, an energy cascade occurs with a flux from the low to the high-frequency range, typical of a turbulent behaviour following a Richardson's-like cascade. Consequently, all the modes are excited through this mechanism. One then understands that building \red{a} ROM to reproduce such complex nonlinear dynamics including a complete transfer of energy is difficult and not achievable with small order subsets.

We now turn to the presentation of nonlinear methods for model order reduction, with a special emphasis on methods based on invariant manifold theory.

\section{Nonlinear methods and invariant manifolds}\label{sec:IMfull}

The aim of this section is to introduce the nonlinear methods for model order reduction based on the concept of invariant manifold. A special emphasis is put on understanding the problem from a geometric perspective in the phase space. 
In this Section, we will focus on explaining the methods  from the equations of motion in modal space, taking Eq.~\eqref{eq:eom_modal} as a starting point. Section~\ref{sec:ROMFE} will consider the case of equations in physical space as starting point, Eq.~\eqref{eq:eom_phys}, with a special attention to methods in a FE context. In the course of this Section, we will also see that a key point is on the extension of the definition of linear modes to \red{the} nonlinear regime. We begin with a short introduction on the mathematical foundation and the developments in the theory of invariant manifolds from the dynamical system point of view.

\subsection{Invariant manifolds for dynamical systems}\label{subsec:dynasys}

Dynamical systems theory offers a geometrical point of view on the global organization of trajectories inside the phase space, thus giving a complete understanding of the long-term behaviour of  solutions. The phase space is structured by the fixed points and the invariant manifolds emerging from their linearized eigendirections with their stability dictated by the eigenspectrum~\cite{gucken83,Kuznetsov,Wiggins}. The center manifold theorem~\cite{Carr,Kelley,CarrMuncaster,ST2015} has long been used as a major tool in the spirit of model-order reduction. Using the terminology introduced in~\cite{Haken81,MannevilleENG}, the long-term dynamics is driven by the central modes and reduction to the center manifold allows an {\em adiabatic} elimination of the slave coordinates. 

Reduction to center manifolds and invariant manifolds has then been used in a number of context in the community of applied mathematics, see {\em e.g.}~\cite{Coullet83,Muncaster1,Muncaster2,Roberts89,RobertsRod93,ROBERTS1997,Llave97}, see also the concept of inertial manifold as exponentially attracting invariant and finite subspaces~\cite{FOIAS1988,FOIAS89,Temam1989,DEBUSSCHE1991}. The method has been used in particular in fluid dynamics for model-order reduction of different problems, see {\em e.g.}~\cite{Mielke,TITI1990,MannevilleENG,LIU2004297,HaragusIooss,Kawahara,Carini2015}, but also in unsteady magnetic dynamos~\cite{LIU2004297} and plasma physics~\cite{BURBY2020}. For conservative or near-conservative systems, a straightforward application of center manifold is however more difficult due to the small (or vanishing) decay rates. 

An important step with regard to the general understanding of the invariant manifold theory and its link to other important theorems from dynamical systems (center manifold,  normal form approach), has been realized with the introduction of the parametrisation method for invariant manifolds by Cabr{\'e}, Fontich and de la Llave~\cite{Cabre1,Cabre2,Cabre3}. The book by Haro {\em et al.}~\cite{Haro} gives a complete presentation and the reader is referred to the first chapters for an accurate understanding. Here a very short presentation of the main ideas is given following their notations and for the case of the computation of invariant manifolds of vector fields in the vicinity of a fixed point. An autonomous dynamical system is considered as:
\begin{equation}
\dot{\zvec} = \Fvec (\zvec),
\label{eq:eomphys}
\end{equation}
with $\zvec = [z_1, ..., z_n]^t$ a $n$-dimensional vector and $\Fvec$ the nonlinear vector field. Let  $\zvec_{\star}$ be a fixed point,  such that $\Fvec(\zvec_{\star}) = \zervec$ and  let  ${\mathcal W}$ be a $d$-dimensional invariant manifold (with $d \ll n$), tangent to the  linear $d$-dimensional subspace $V^L$ at $\zvec_{\star}$. A parametrisation is introduced as a nonlinear mapping between the original coordinates $\zvec$ (of dim. $n$) and newly introduced coordinates $\svec = [s_1, ..., s_d]^t$ ($d$-dimensional vector of master coordinates). This nonlinear change of coordinates is written in general form as 
\begin{equation}
\zvec = \Wvec (\svec),
\label{eq:NLmap}
\end{equation}
where $\Wvec$ is  unknown at this stage. The reduced-order dynamics, {\em i.e.} the dynamics on the invariant manifold, also unknown at this stage,  writes
\begin{equation}
\dot{\svec} = \fvec (\svec),
\label{eq:ROMdyna}
\end{equation}
with $\fvec (\zervec) = \zervec$. 
To compute both $\Wvec$ and $\fvec$, we replace 
the nonlinear mapping \eqref{eq:NLmap} in Eq.~\eqref{eq:eomphys}. Using the chain rule, one obtains by differentiating \eqref{eq:NLmap} with respect to time $\dot{\zvec} = \DD\Wvec (\svec) \dot{\svec} = \Fvec (\Wvec (\svec))$, with $\DD\Wvec$ the derivative of $\Wvec$. Using Eq.~\eqref{eq:ROMdyna}, one finally obtains: 
\begin{equation}
\Fvec (\Wvec (\svec)) = \DD\Wvec (\svec) \fvec (\svec),
\label{eq:invariance}
\end{equation}
sometimes written $\Fvec \circ \Wvec = \DD\Wvec \, \fvec$ \cite{Haro}. Since this equation is independent of time, it enforces the invariance property of ${\mathcal W}$. It is known as the {\em invariance} equation and enables to compute high-order expansions of both $\Wvec$ and $\fvec$.

The remaining of the calculation as presented in~\cite{Haro} introduces polynomial expansions for the two unknowns $\Wvec$ and $\fvec$ into the invariance equation, from which order-by-order identification leads to the so-called co-homological equations, related to the tangent (master coordinates) and normal (slave coordinates) parts. Full details are given in~\cite{Haro} and a short summary is proposed in Appendix~\ref{app:IMparam}. Those co-homological equations enables \red{one} to compute, order by order, the two unknowns $\Wvec$ and $\fvec$. However, their solution is not unique and a choice on the parametrisation has to be done.  Haro {\em et al.}  introduces the two main parametrisation methods that one can use to solve the problem. 

The first one is called {\em the graph style} and leads to a functional relationship between slave and master coordinates, in which the master coordinates are only linearly related to the original ones. 
The second one is the {\em normal form style} and leads to the introduction of new coordinates, nonlinearly related to the original ones. The idea in this case is to simplify as much as possible the reduced-order dynamics, by keeping only the resonant monomials, and discarding all other non-essential terms for the dynamical analysis. This leads to a more complex calculation, and a full nonlinear mapping between original coordinates and reduced ones. The drawback is that calculations are a bit more involved (which is particularly true when there are numerous internal resonances to handle). The advantage is that the parametrisation is able to go over the foldings of the manifold.
Finally, since other parametrisations exist (an infinite number), {\em mixed styles} can also be used, but the first two are the extreme cases and mixed styles are only variations using both graph and normal form styles.

Now restricting ourselves to the case of vibratory systems, it is important to distinguish the conservative and dissipative case. In the conservative case, the eigenspectrum is purely imaginary with pairs of complex conjugates $\{ \pm i \omega_p\}_{p=1,...,N}$. A center theorem from Lyapunov then states the existence of two-dimensional manifolds densely filled with periodic orbits, for each couple of imaginary eigenvalues~\cite{Lyapunov1907,Kelley2,Gordon71,Weinstein73}, under the assumption of non-resonance condition. These invariant manifolds are named Lyapunov subcenter manifold (LSM). The existence of these LSM leads to the definition of {\em nonlinear normal modes} (NNM), which are the extension of the (linear) eigenmodes (LM) to the nonlinear range. Two properties of the linear modes can be extended to the nonlinear case, giving two complementary definitions of an NNM. The first one, historically proposed by Rosenberg in the sixties, and modernised by many contributions since, is to define an NNM as a {\it family of periodic orbits} \cite{Rosenberg62,Rand74,Pecelli79,KingVakakis94,KingVakakis96,VakakisMSSP97,VakakisNNM,KerschenNNM09}.
\red{From this definition, numerous investigations tackled the problem of constructing NNMs thanks to perturbative approaches, that could be inserted directly into the PDE of motion, also including internal resonances~\cite{nayfehnayfeh94,nayfeh95,nayfeh96,nayfehcarbo97,NayfehROM98,nayfehcarboChin99,LACARBONARA2004nnm}
}. Then, Shaw \& Pierre proposed in 1991 to define an NNM as an {\em invariant manifold} of the phase space. This second definition naturally allows the derivation of accurate reduced-order models: this will be the subject of the next sections. In the conservative case, both definitions are equivalent. For dissipative vibratory systems, existence theorems for the manifolds have been proven only recently \red{by Haller and Ponsioen}~\cite{Haller2016}, leading to the notion of spectral submanifolds (SSM). This case will be more deeply analysed in Section~\ref{subsec:SSM}.

The presentation will now follow the chronological order, which is also coherent with the separation into graph style and normal form style proposed by Haro {\em et al.} in~\cite{Haro}.

\subsection{The graph style: nonlinear normal modes as invariant manifolds}\label{sec:IMSP}

The first step for defining ROMs based on invariant manifold theory has been proposed by Steve Shaw and Christophe Pierre in the early 1990s. 
The key idea is to use the center manifold theorem, as given in most classical textbooks on dynamical systems (see {\em e.g.}~\cite{Carr,gucken83,MannevilleENG,Wiggins}) as a technical method in order to derive the equations describing the geometry of the invariant manifold in phase space. 
Replacing this calculation in the light of the parametrisation method, one understands that the technique as proposed by Shaw and Pierre~\cite{ShawPierre91,ShawPierre93,ShawPierre94} for conservative nonlinear vibratory systems is equivalent to the parametrisation method of LSM following the graph style.

In the next sections, the method and main results from the graph style approach, following the developments led by Shaw, Pierre and coworkers, will be reviewed. In order to introduce progressively the details, section~\ref{subsec:SP1} considers the case of a single master mode. Then section~\ref{subsec:SP2} extends the results to multiple master coordinates, opening the doors to more complex geometry of invariant manifolds. Finally, section~\ref{subsec:SP3} summarizes all the results obtained with the method, including the addition of damping and forcing, piecewise linear restoring force, and numerical computations.

\subsubsection{Two-dimensional invariant manifold}\label{subsec:SP1}

In this section, we restrict ourselves to the case of a single master coordinate, labelled $m$.
Rewriting the $p^{\mathrm{th}}$ equation of~\eqref{eq:eom_modal} at first-order, one obtains, $\forall \, p =1,...N$:
\begin{subequations}\label{eq:eomIMo1}
\begin{align}
\dot{x}_p &= y_p \\
\dot{y}_p &= - \omega_p^2 x_p - f_p(x_1, ..., x_N),
\end{align}
\end{subequations}
with $y$ the velocity and $f_p$ the function grouping quadratic and cubic nonlinear terms:
\begin{equation}
f_p(x_1, ..., x_N) = \sum_{i=1}^N \sum_{j=1}^N g^p_{ij} x_i x_j + \sum_{i=1}^N \sum_{j=1}^N \sum_{k=1}^N h^p_{ijk} x_i x_j x_k.
\end{equation}
The idea is to assume the existence of a functional relationship between all the slave coordinate $s$ and the master one $m$, {\em i.e.} $\forall \, s \neq m$, there exist two functions $U_s$ and $V_s$, solely depending on the displacement and velocities of the master coordinates $(x_m,y_m)$, such that
\begin{subequations}\label{eq:IMFR}
\begin{align}
x_s &= U_s (x_m,y_m), \\
y_s &= V_s (x_m,y_m).
\end{align}
\end{subequations}
At this stage $U_s$ and $V_s$ are the unknowns, and it is important to remark that: 
\begin{itemize}
\item the dependence is written for both displacements and velocities. Since oscillations occur on two-dimensional surface involving two independent coordinates, the velocities shall not be neglected. This also reflects the fact that the eigenspace of a mode is two-dimensional, with eigenvalues $\pm i\omega$.
\item a functional dependence between the modal variables is searched for, which is different from a change of coordinate or nonlinear mapping introducing new coordinates. This is typical of the {\em graph style} for the parametrisation of the invariant manifold.
\end{itemize}

The methodology to find the unknown functions  $U_s$ and $V_s$ consists in deriving Eqs.~\eqref{eq:IMFR} with respect to time and substitute in the dynamical equations~\eqref{eq:eomIMo1} whenever possible in order to eliminate all explicit dependence on time, thus following a similar development as the one shown in Section~\ref{subsec:dynasys} to arrive at the invariance equation. The development leads to, $\forall \, s \neq m$:
\begin{subequations}\label{eq:IMgeometry}
\begin{align}
\frac{\partial U_s}{\partial x_m} y_m + \frac{\partial U_s}{\partial y_m} \left( -\omega_m^2 U_m - f_m \right) & = V_s (x_m,y_m), \\
\frac{\partial V_s}{\partial x_m} y_m + \frac{\partial V_s}{\partial y_m} \left( -\omega_m^2 U_m - f_m \right) & = -\omega_s^2 U_s   - f_s.
\end{align}
\end{subequations}
Eqs.~\eqref{eq:IMgeometry} are a set of $2N-2$ partial differential equations depending on the master coordinates $(x_m,y_m)$. They describe the geometry of the two-dimensional invariant manifold in the $2N$-dimensional phase space. 
The solutions of Eqs.~\eqref{eq:IMgeometry}  will give the $N-1$ unknown functions $(U_s,V_s)$. Unfortunately, these equations contains all the nonlinearities of the initial problem through the $f_p$ functions. Consequently obtaining simple solutions to \eqref{eq:IMgeometry} is generally out of reach. In their first papers, Shaw and Pierre proposed to solve them using asymptotic expansions. This will be detailed next since it gives the first significant terms in the developments, that can be used for direct comparisons to other methods. In subsequent developments, They also propose to solve \eqref{eq:IMgeometry} numerically. This will be reviewed in Section~\ref{subsec:SP3}.

Since the invariant manifold is tangent to its linear counterpart close to the origin, the functions $(U_s,V_s)$ shall contain neither constant terms, nor linear ones. Consequently the asymptotic expansion begins with second-order terms. The analytical developments to arrive at the coefficients are given in~\cite{PesheckBoivin}, we here simply recall the obtained results. Up to third order, the solution reads:
\begin{subequations}\label{eq:IMgeometryo3}
\begin{align}
x_s &= a_{sm} x_m^2 + b_{sm} y_m^2 + c_{sm} x_m^3 + d_{sm} x_m y_m^2,\\
y_s &= \alpha_{sm} x_m y_m + \beta_{sm} x_m^2 y_m^2 + \gamma_{sm} y_m^3.
\end{align}
\end{subequations}
One can note in particular that all the coefficients of the multivariate polynomials $(x_m,y_m)$ are not present. Indeed, some of them are vanishing due to the conservative nature of the nonlinear restoring force assumed from the beginning. However, adding more terms to the initial problem ({\em e.g.} damping, gyroscopic force, ...) will complete the polynomial expansions with other coefficients.  
The expressions of the quadratic coefficients, which will be used after for explicit comparisons with other reduction methods, reads:
\begin{subequations}\label{eq:coefquadIM}
\begin{align}
a_{sm} &= \frac{2\omega_m^2 - \omega_s^2}{\omega_s^2 (\omega_s^2 - 4 \omega_m^2)}g^s_{mm}, \\
b_{sm} &= \frac{2}{\omega_s^2 (\omega_s^2 - 4 \omega_m^2)}g^s_{mm}, \\
{\alpha}_{sm} &= \frac{-2}{\omega_s^2 - 4 \omega_m^2}g^s_{mm}.
\end{align}
\end{subequations}
One can note the two following important features: (i) the coefficients are proportional to $g^s_{mm}$ which is the coefficient of the invariant-breaking term $X_m^2$ on slave mode $s$. 
(ii) The formulas are valid as long as no second-order internal resonance $\omega_s = 2 \omega_m$ exist between slave and master coordinates. This is fully logical since in that case a strong coupling exists between the two modes and reduction to a single master mode $m$ is not meaningful.

The reduced dynamics on the invariant manifold is found by substituting the functional relationships \eqref{eq:IMFR} into the equation of motion for the master mode $m$~:
\begin{equation}\label{eq:romIM1m}
\ddot{x}_m + \omega_m^2 x_m + f_m (U_1(x_m,y_m), ..., x_m, y_m, ..., U_N(x_m,y_m)) = 0.
\end{equation}
Given the expressions of the coefficients in Eq.~\eqref{eq:coefquadIM}, Eq.~\eqref{eq:romIM1m} can be explicitly written as~\cite{PesheckBoivin,Pesheck00,YichangICE}:
\begin{equation}\label{eq:romIM1mdev}
\ddot{x}_m + \omega_m^2 x_m + g^m_{mm} x_m^2 + x_m \left(\underset{s\neq m}{\sum_{s=1}^N} 2\,g^m_{ms}g^s_{mm}\left[\frac{2\omega_m^2 - \omega_s^2}{\omega_s^2 (\omega_s^2 - 4 \omega_m^2)}x_m^2 +  \frac{2}{\omega_s^2 (\omega_s^2 - 4 \omega_m^2)}y_m^2 \right] \right) + h^m_{mmm} x_m^3 = 0.
\end{equation}
One can note in particular that the ``self-quadratic'' term $g^m_{mm}x_m^2$ stays in the reduced dynamics. The cubic term $h^m_{mmm} x_m^3$ is balanced by two other cubic terms, one involving  the $x_m^3$ monomial, while the other involves $x_m y_m^2$ and the coefficient is a summation on all the slave modes, showing how their effect is gathered in the nonlinear dynamics on the invariant manifold. The expression assumes a third-order truncation in both the relationship between slave and master coordinates as well as for the reduced dynamics. Asymptotic developments can be pushed further at the expense of more involved derivations. We now turn to the generalization with a multi-mode manifold.

\subsubsection{Multi-dimensional invariant manifold}\label{subsec:SP2}

The multi-dimensional extension of the previous development has been first given in~\cite{PesheckBoivin}, in order to propose ROMs with a larger number of master modes that can handle internal resonance and more complex nonlinear dynamical phenomena. The methodology is unchanged as compared to the previous case but is complexified by the fact that numerous master modes are taken into account. The starting point is to distinguish master and slave coordinates. For the sake of simplicity, let us note as $1,..., m$ the index of the $m$ master modes and $m+1, ..., N$ the index of the remaining slave modes. The functional relationship now reads, $\forall s \in [m+1,N]$ (slave coordinates):
\begin{subequations}\label{eq:IMmulti}
\begin{align}
x_s &= U_s (x_1, y_1, ..., x_m,y_m), \\
y_s &= V_s (x_1, y_1, ..., x_m,y_m).
\end{align}
\end{subequations}
In order to derive the unknown functions $(U_s,V_s)$, $s=m+1, ... N$, one has to solve:
\begin{subequations}\label{eq:imgeomeq}
\begin{align}
&\sum_{r=1}^{m}\left( \frac{\partial{U_s}}{\partial{x_r}}y_r+\frac{\partial{U_s}}{\partial{y_r}} \left[ -\omega_r^2 U_r-f_r\right] \right)\, = \, V_s,\\
&\sum_{r=1}^{m}\left( \frac{\partial{V_s}}{\partial{x_r}}y_r+\frac{\partial{V_s}}{\partial{y_r}}\left[-\omega_r^2 U_r-f_r\right]\right)\, = \, -\omega_s^2 U_s-f_s.
\end{align}
\end{subequations}
These $2(N-m)$ equations describe the geometry of the $2m$-dimensional invariant manifold in the phase space. Again, the solution of these PDE is generally out of reach, and asymptotic solutions up to order three are a convenient way to work it out. The method can also be written in a systematic manner, highlighting the repeating structures appearing at each order and thus opening the doors to automated high-order solutions. The individual coefficients up to order three are given in~\cite{PesheckBoivin,Pesheck00} under matrix form instead of explicit expressions. The reduced dynamics on the manifold is simply found by replacing \eqref{eq:IMmulti} in the master coordinates in \eqref{eq:eomIMo1}.

\subsubsection{Applications}\label{subsec:SP3}

The first applications of the invariant manifold approach have been mainly proposed on beam examples: a simply supported beam resting on a nonlinear elastic foundation is considered in~\cite{ShawPierre94,Shaw94}, a linear beam with local nonlinear springs attached either at the ends (torsional springs) in~\cite{ShawPierre94}, or at center (transverse spring) in~\cite{PesheckBoivin}, and  a nonlinear rotating beam in~\cite{PesheckASME2002}. \ot{Applications to planar frames and simply-supported beam have also been reported in~\cite{MAZZILLI2001,MAZZILLI2004}.}

An important advantage of the method, based on the center manifold theorem, is to express the geometry of the invariant manifold (the reduction subspace) in terms of a partial differential equation describing its geometry in phase space, Eq.~\eqref{eq:IMgeometry} for single master coordinate and \eqref{eq:imgeomeq} for the multi-dimensional manifold with $m$ master coordinates. Consequently, all the numerical tools for solving PDEs can be implemented in order to propose a fully numerical yet accurate computation of the manifold and the reduced dynamics, thus bypassing the intrinsic limitation of any asymptotic development. However, the starting point assuming a graph relationship inherently precludes the method to overcome the possible folding of the manifold~\cite{PesheckJSV,Haro,BLANC2013}, so that in any case the method will have a limit in terms of amplitude at the first folding point.

\begin{figure}[h!]
\includegraphics[width=.57\linewidth]{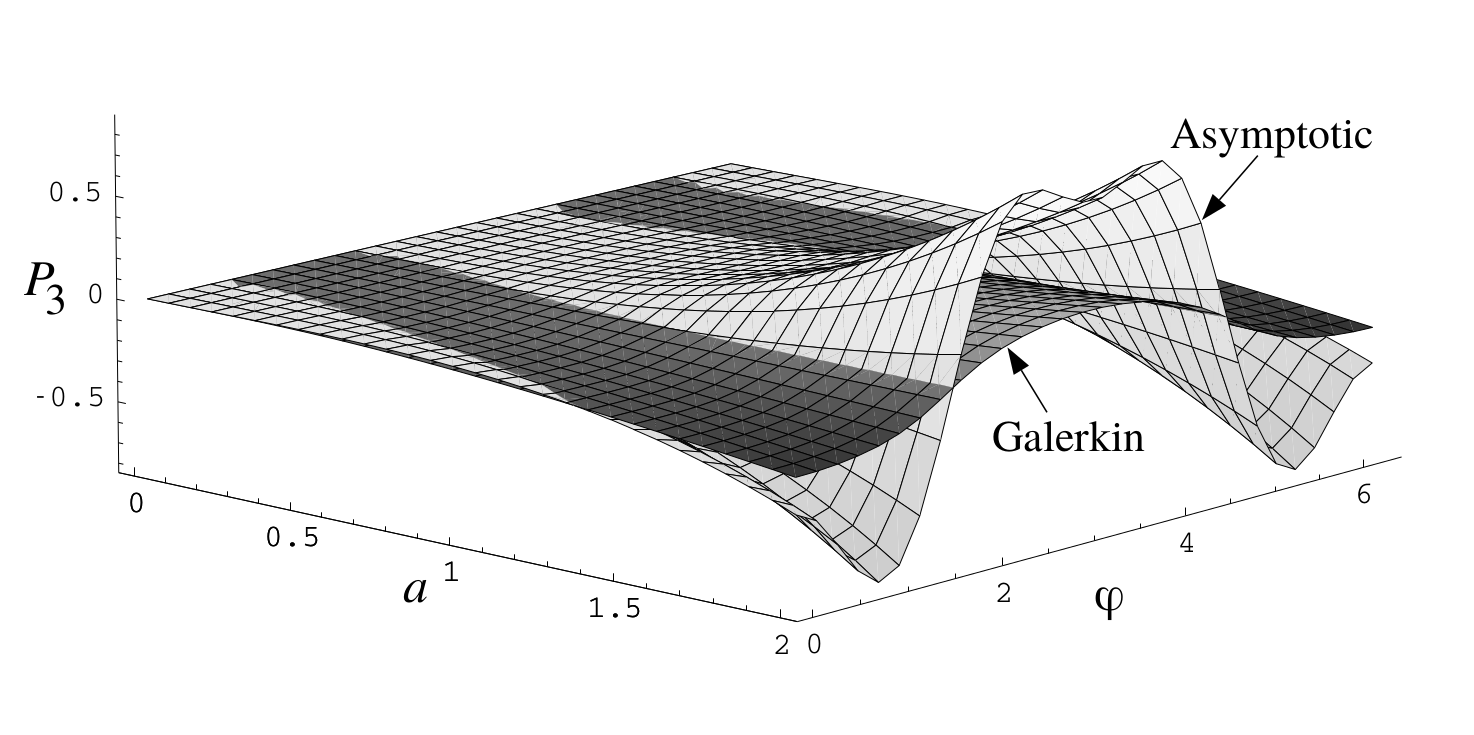}$\quad$\includegraphics[width=.38\linewidth]{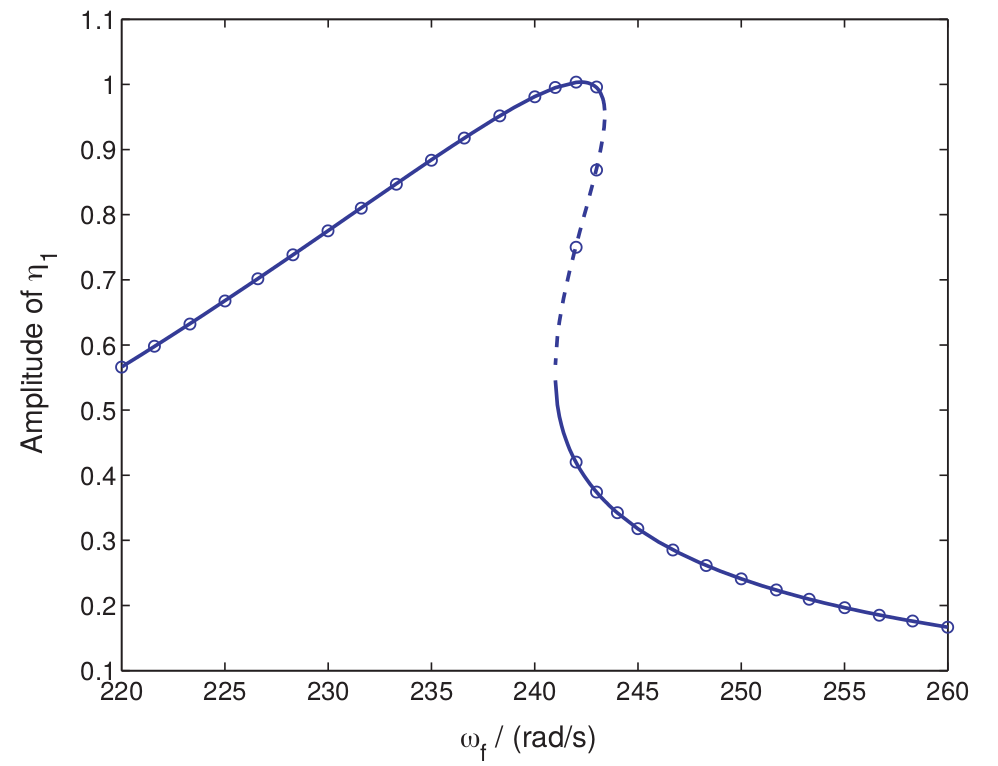}
\caption{(a) Comparison of the invariant manifold as computed from a third-order asymptotic development and numerically obtained by direct numerical solution of \eqref{eq:IMgeometry} (solution depicted as ''Galerkin'' in the figure) for a linear beam with a nonlinear torsional spring at boundary, \red{in an amplitude-phase $(a,\phi)$ representation}. 
Reprinted \red{with permission} from~\cite{PesheckJSV}. (b) Frequency response for the same beam with damping and forcing, comparison between reference full-order solution (continuous line: stable solution, dotted line: unstable solution) and ROM with one master mode, graph style parametrisation (\red{circles}). Reprinted from~\cite{JIANG2005H}.}
\label{fig:ShawPierreFig}
\end{figure}

Based on this idea, a numerical procedure has been developed in~\cite{PesheckJSV} for numerical computation of two-dimensional manifold, and has then been extended to the case of multiple mode invariant manifolds in~\cite{JIANG2005}. Using this numerical procedure, extension of the method in order to properly take into account forcing and damping in order to compute frequency responses has been proposed in~\cite{JIANG2005H}, whereas the forced case is also considered in~\cite{GABALE2011} using series expansions. Also, the case of piecewise linear systems has been tackled in~\cite{ShawChen96,JIANG2004PW}. With regard to applications, the case of a rotating beam is considered in~\cite{JIANG2005}, and a rotating shaft in~\cite{LegrandIM}. \red{Along} the same lines, different numerical procedures have been proposed in~\cite{NORELAND2009,BLANC2013,RENSON2014} to solve the nonlinear PDEs of the invariant manifold, and a more general review of numerical methods (including other approaches) is reported in~\cite{RENSON2016}.  Fig.~\ref{fig:ShawPierreFig} shows two illustrations from these works. Finally, one can also note that the invariant manifold parametrisation with graph style has also been used in combination with Lyapunov-Floquet transform for systems with periodic coefficients~\cite{SINHA2005985}, and the technique for augmenting the state space for forced systems has been investigated in~\cite{Redkar2008,GABALE2011}.

With regard to finite element applications, one can note that several examples using a FE procedure in order to semi-discretize the problem, have been implemented, for example a linear FE beam with a nonlinear rotational spring at one end is considered in~\cite{PesheckJSV,JIANG2005H} and a one-dimensional finite-element model representing the axial and transverse motions of a cantilever rotating beam is selected in~\cite{Polarit}. \ot{Applications to planar frames discretized by the FE method are also shown in~\cite{SOARES2000,Mazzilli11}.} But in all these cases, a relative simple geometry is considered and the first step is the full projection of the system equation on the modal equations. As it will be discussed in Section~\ref{sec:ROMFE}, the problem of very large FE models having millions of DOFs -- thus preventing  such a first step -- has not been addressed in these studies.

As a short conclusion, the method strictly follows the graph style  for the parametrisation of \red{an} invariant manifold. Expressing the geometry of the invariant subset as a PDE is an advantage since opening the doors to numerical solution. However, the assumption of a graph relationship between slave and master coordinates puts a clear limitation to the method that will never be able to pass through folding points of the manifold. The method has however an important generality and versatility and shall be used in a number of contexts.

\subsection{Normal form approach}\label{sec:NFmodal}

The normal form approach, used with the purpose of analysis and model-order reduction  of vibratory systems, has been  proposed and developed from the following papers~\cite{Jezequel91,touzeLMA,touze03-NNM}. It relies on a complete normal form calculation, following the general guidelines of dynamical systems theory~\cite{Poincare,Dulac1912,ELPHICK1987,MannevilleENG,IoossAdel,HaragusIooss,Iooss88,Murdock,IoossLombardi05,IoossLombardi09},
 adapted to the framework of vibratory systems, and then followed by a truncation to achieve reduction by selecting only a few of the resulting coordinates as master. By doing so, one retrieves an equivalent procedure \red{to} the one proposed \red{for} the parametrisation method of invariant manifolds, \red{but now} with a normal form style~\cite{Haro}.

In its first derivation reported in~\cite{touzeLMA,touze03-NNM}, the complete normal form is computed by 
keeping oscillator-like equations (with second-order derivatives in time), to better fit the usual mechanical framework, thus arriving at a real-valued normal transform. On the other hand, all mathematical derivations use a complex formulation with diagonalized linear part~\cite{IoossAdel,ELPHICK1987,Jezequel91}.  A complete nonlinear mapping is thus derived, allowing one to express the dynamics with new coordinates related to the individual invariant manifolds ascertained in the previous section. Consequently, the method generalizes the asymptotic approach described in \ref{subsec:SP2}, since the complete change of coordinates is derived. The master coordinates are selected after the transform thus offering versatility to the method and easy implementation of ROMs with arbitrary number of master modes. On the other hand, the calculation as shown in~\cite{touzeLMA,touze03-NNM} has been limited to the third-order.


\subsubsection{Method and main results}\label{sec:NFmethod}

The derivation of the complete nonlinear mapping for conservative nonlinear vibratory systems expressed in the modal basis, \red{{\em i.e.}} taking Eq.~\eqref{eq:eom_modal} as starting point, is established in~\cite{touzeLMA,touze03-NNM,TouzeCISM}, following the general guidelines of normal form theory~\cite{IoossAdel,Murdock}. In essence, the calculations are led order by order, and the procedure at each order is to inject an unknown nonlinear mapping, derive the associate\red{d} homological equation~\cite{touzeLMA,LamarqueUP}, which is solved by assuming that the goal is to eliminate as \red{many} monomials as possible, to arrive at a reduced dynamics (the normal form) having the simplest expression. In case of no internal resonance, the normal form is linear (Poincar{\'e}'s theorem), whereas existence of nonlinear resonance leads to a more complex normal form where only the resonant monomials finally stay (Poincar{\'e}-Dulac's theorem).

An important feature related to conservative vibratory systems is the presence of {\em trivial resonance} relationships (see Appendix~\ref{app:class} and~\cite{touze03-NNM,TouzeCISM} for the definition), meaning that a vibratory system can never be linearized: the normal form will always contain resonant monomials. Importantly, the monomials connected to trivial resonances have an odd order, meaning for example that cubic terms are especially important as compared to quadratic ones. In particular, the processing of the calculation is to eliminate at order $n$ the non-resonant terms thanks to an order $n$ nonlinear mapping, creating in turn new terms at order $n+1$. Consequently, quadratic terms can be eliminated (under the assumption of no second-order internal resonance), and the effect of this elimination will result in modified cubic terms that can be derived. While the presence of trivial resonance is not a good news \red{from} the mathematical point of view (leading to more involved calculations), it is meaningful in a nonlinear vibration context since the resulting cubic terms will drive the hardening/softening behaviour.

Up to the third order, the nonlinear change of coordinates, following the real formalism proposed in~\cite{touzeLMA,touze03-NNM,TouzeCISM}, \red{can be written}, for each pair of displacements and velocities $(x_k,y_k)$, $\forall \; k=1...N$, \red{as}
\begin{subequations}
\begin{align}
&x_k=
 R_k + 
\sum^N_{i=1}\sum^N_{j=1}
a^k_{ij}
R_i R_j+
\sum^N_{i=1}\sum^N_{j=1}
b^k_{ij}
S_i S_j+
\sum^N_{i=1}\sum^N_{j=1}\sum^N_{l=1}
r^k_{ijl}
R_i R_j R_l+
\sum^N_{i=1}\sum^N_{j=1}\sum^N_{l=1}
u^k_{ijl}
 R_i S_j S_l,
\label{eq:nonlinear_change_x}
\\
&y_k=
S_k + 
\sum^N_{i=1}\sum^N_{j=1}
\gamma^k_{ij}
R_i S_j+
\sum^N_{i=1}\sum^N_{j=1}\sum^N_{l=1}
\mu^k_{ijl}
S_i S_j S_l+
\sum^N_{i=1}\sum^N_{j=1}\sum^N_{l=1}
\nu^k_{ijl}
S_i R_j R_l,
\label{eq:nonlinear_change_y}
\end{align}
\label{eq:nonlinear_change_modal_compact}
\end{subequations}
where the newly introduced {\em normal} coordinates $R_i$ and $S_i=\dot{R}_i$ are respectively homogeneous to a displacement and a velocity. The calculation has been done once and for all with $N$ variables, and the full expressions of all the reconstruction coefficients $a^k_{ij}$, $b^k_{ij}$, $\gamma^k_{ij}$, $r^k_{ijl}$, $u^k_{ijl}$, $\mu^k_{ijl}$, and $\nu^k_{ijl}$ are  given in~\cite{touze03-NNM,artDNF2020}. The nonlinear mapping takes velocities into account, based on the fact that in vibration theory, velocities are mandator\red{ily} needed as second independent variables in order to construct oscillations as closed orbit\red{s} in a two-dimensional subspace. It is identity-tangent meaning that at the lowest order, the usual eigenspaces are retrieved. Higher-order (quadratic and cubic) terms \red{lead to expressions for}  the curvature of the invariant manifold in phase space, and thus the dependence of modal quantities with respect to amplitude. 

As shown in~\cite{touzeLMA,touze03-NNM}, the method expresses the  reduced dynamics in an invariant-based span of the phase space. \red{These can be} written for the general case where no internal resonance exists between the eigenfrequencies of the system. When an internal resonance is present, some terms are vanishing in Eqs.~\eqref{eq:nonlinear_change_modal_compact}, leading to extra terms staying in the normal form of the system. 

The reduction step consists \red{of} selecting a few master {\em normal} coordinates, say $m \ll N$, and \red{eliminating} all the others. Assuming for simplicity that the master coordinates are for $p=1...m$, this means that $\forall j=m+1, ..., N$, $R_j=S_j=0$, hence transforming the one-to-one diffeomorphism \eqref{eq:nonlinear_change_modal_compact} to a nonlinear mapping parametrising the invariant manifold associated \red{with} the master coordinates.

In case of no internal resonance, the reduced dynamics on this $m$-dimensional manifold can be written explicitly as, $\forall r=1,\, ...,\, m$:
\begin{equation}
\begin{split}
&
\ddot{R}_r+\omega^2_rR_r+
({A}^r_{rrr}+{h}^r_{rrr})R_r^3+
({B}^r_{rrr})R_r\dot{R}_r^2\\
&
+R_r
\sum^m_{j\neq r}
({A}^r_{jjr}+{A}^r_{jrj}+{A}^r_{rjj}+3{h}^r_{rjj})R_j^2
+R_r
\sum^m_{j\neq r}
({B}^r_{rjj})\dot{R}_j^2
+
\dot{R}_r
\sum^m_{j\neq r}
({B}^r_{jjr}+{B}^r_{jrj})R_j\dot{R}_j
=0.
\end{split}
\label{eq:ROMNF}
\end{equation}
This dynamical equation is the real normal form of the problem, where only the resonant monomials corresponding to trivial resonances are present, all other terms being cancelled. As \red{stated}, quadratic terms have disappeared and only cubic terms are present. The result of this operation appears through the new fourth-order tensors $\At$ and $\Bt$, that gathers the elimination of the quadratic terms and whose expression only contains quadratic coupling coefficients $g^p_{ij}$. Their expressions from the modal basis can be found in~\cite{touze03-NNM} and are here recalled:
\begin{subequations}\label{eq:ABmodal}
\begin{align}
A^r_{ijk}= \sum^N_{s=1}2\,g^r_{is}{a}^s_{jk},\\
B^r_{ijk}= \sum^N_{s=1}2\,g^r_{is}{b}^s_{jk}.
\end{align}
\end{subequations}

One can note in particular that the same invariant subspaces are computed as in Section~\ref{sec:IMSP},   only the parametrisation and thus the meaning of the reduced coordinates, is different. In the graph style, the master coordinates are a subset of the original one\red{s} $(\x,\dot{\x})$. In the normal form style, new coordinates $(\R,\dot{\R})$, nonlinearly related to the original ones, are introduced. 

To be more specific, let us compare the geometry of the manifold given by the two methods when restrict\red{ed} to the case of a single master coordinate. From the normal form approach, the geometry of the manifold is expressed by Eqs.~\eqref{eq:nonlinear_change_modal_compact}. Assuming only mode $m$ as master, limiting to the second order for the sake of simplicity, and replacing the coefficients $\at$, $\bt$ and $\ct$ by their explicit expressions given in~\cite{touze03-NNM}, the geometry is given by, $\forall s\neq m$
\begin{subequations}\label{eq:NFgeommaster1}
\begin{align}
x_s &= \frac{(2\omega_m^2 - \omega_s^2)g^s_{mm}}{\omega_s^2 (\omega_s^2 - 4 \omega_m^2)} R_m^2 + \frac{2\,g^s_{mm}}{\omega_s^2 - 4 \omega_m^2} \dot{R}_m^2,\\
y_s &= \frac{2\,g^s_{mm}}{4 \omega_m^2-\omega_s^2} R_m\dot{R}_m.
\end{align}
\end{subequations}
These equations are exactly those given in \eqref{eq:IMgeometryo3}-\eqref{eq:coefquadIM}, meaning that at second-order of the development, the two different styles of parametrisation gives the same quadratic terms for the geometry of the manifold on the slave modes. The developments then start to depart one from another at the next orders, due \red{to} the use of different coordinates. For the reduced dynamics, the difference starts to appear from the second-order as shown next.

The reduced dynamics obtained with \red{the} normal form approach restricted to a single master coordinate $R_m$ reads
\begin{equation}\label{eq:NFsingleNNM}
\ddot{R}_m + \omega_m^2 R_m + (A^m_{mmm} + h^m_{mmm}) R_m^3 + B^m_{mmm} R_m \dot{R}_m^2 = 0.
\end{equation}
Comparing to Eq.~\eqref{eq:romIM1mdev}, one can observe in particular that Eq.~\eqref{eq:romIM1mdev} contains a quadratic term which is not present in~\eqref{eq:NFsingleNNM}. This difference is only related to  the  meaning of the variables used in each method and their nonlinear relationship. Introducing the {\em normal} variables defined by Eqs.~\eqref{eq:nonlinear_change_modal_compact} in the reduced dynamics given by Eq.~\eqref{eq:romIM1mdev}, the same equation is obtained. This is demonstrated in Appendix~\ref{app:IMtoNF}. In particular the two methods predict exactly the same and correct hardening/softening behaviour. Using a perturbative expansion, the nonlinear frequency/amplitude relationship can be written as $\omega_{NL} = \omega_m (1 + \Gamma_m a^2)$, with $\omega_{NL}$ the nonlinear radian frequency, $a$ the amplitude, and $\Gamma_m$ the nonlinear coefficient dictating the type of nonlinearity. In each case,  the same coefficient is found as:
\begin{equation}
\Gamma_m = - \dfrac{5}{12\;\omega^2_{m}} \left(\dfrac{g^m_{mm}}{\omega_{m}}\right)^2  + \dfrac{3}{8\;\omega^2_{m}}\left( h_{mmm}^m - \sum_{\underset{s\neq m}{s=1}}^N 2\left(\dfrac{g^s_{mm}}{\omega_{s}}\right)^2 \left( 1+ \dfrac{4\,\omega^2_m}{3(\omega^2_s-4\,\omega^2_m)}
\right)\right),
\label{eq:gammaNF}
\end{equation}
where the symmetry relationships on the $g^p_{ij}$ coefficients have been used, see Eq.~\eqref{eq:zesymg} in Appendix~\ref{sec:ghsym}.

\subsubsection{Applications}

The normal form approach and its use in model-order reduction has been first extended to handle the case of linear modal damping ratio in the change of coordinates~\cite{TOUZE:JSV:2006}, thus opening the doors to the computation of forced-damped dynamics and frequency responses, by also adding an external forcing with a first-order assumption under its modal formulation. In this case, special care has to be taken in order to follow the trivial resonances, that are destroyed with added damping. As shown in~\cite{TOUZE:JSV:2006,TouzeCISM}, this can be done using parameter-dependent normal forms as derived in~\cite{IoossAdel,HaragusIooss}, enforcing the dissipative case to tend to the conservative case when damping is vanishing. Thanks to this derivation, the reduced dynamics driven by the master modes displays a damping factor that takes into account the damping coefficients of all the slave modes, ensuring  a more proper estimate of the decay rates on the invariant manifold. As an interesting particular result, it has been shown in~\cite{TOUZE:JSV:2006} that the damping can \red{affect} the type of nonlinearity.

\begin{figure}[h!]
\includegraphics[width=\linewidth]{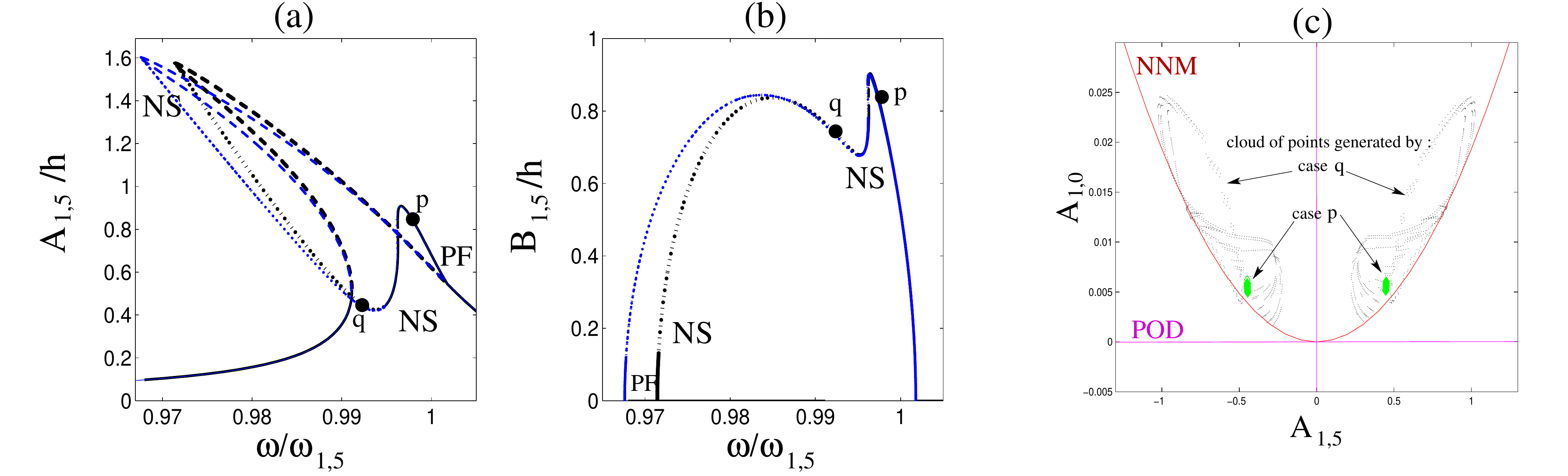}
\caption{Reduced-order models using normal form approach for circular cylindrical shells featuring 1:1 resonance. (a)-(b) Frequency response to harmonic excitation $\omega$ in the vicinity of mode (1,5), with eigenfrequency $\omega_{1,5}$, from~\cite{TOUZE:JSV:2006}. $A_{1,5}$ is the coordinate of the driven mode and $B_{1,5}$ the companion mode. Black: reference, full-order solution. Blue: ROM with two master coordinates. NS: Neimark-Sacker bifurcation, PF: pitchfork bifurcation. \red{Solid line: stable solutions, dashed and dotted lines: unstable solutions.} (c) Partial representation of the phase space with two retained coordinates, the driven mode $A_{1,5}$ and the axisymmetric slave mode $A_{1,0}$. Poincar{\'e} section of the temporal solutions obtained from points $p$ and $q$. POD axes in magenta, invariant manifold (NNM) in red. \red{Figure reworked from~\cite{TOUZE:JFS:2007}.}}
\label{fig:normalformFIG}
\end{figure}

Applications to beams have been first reported in~\cite{touze03-NNM,TOUZE2004CS}. Then the case of circular cylindrical shells \red{has} been tackled in~\cite{TOUZE:JSV:2006}, including a comparison with the POD method in~\cite{TOUZE:JFS:2007}. Interestingly, these shells have degenerate eigenmodes leading to 1:1 internal resonance and a complex dynamics including Neimark-Sacker (NS) bifurcation points. The bifurcation diagram in the frequency response function (FRF) \red{was} very well predicted by the ROM including two master coordinates, as shown in Fig.~\ref{fig:normalformFIG}(a-b). This underlines that a minimal model with only two master coordinates, computed directly from the model equations, is able to retrieve all the dynamical features of the full-order solution, including the quasiperiodic solutions developing in between the two NS bifurcations. Fig.~\ref{fig:normalformFIG}(c) shows a geometrical interpretation in phase space. Two clouds of points obtained by Poincar{\'e} section and generated from the full-order system are shown. They have been respectively obtained for a periodic solution (point $p$) and a quasiperiodic solution (point $q$). The magenta axes are the reduction directions given by the POD method, undoubtedly showing that two directions are necessary in this plane to correctly represent the data~\cite{TOUZE:JFS:2007,AmabiliPOD1}. On the other hand, the \red{section through}  the 4-dimensional invariant manifold in this plane shows that the reduced subspace goes exactly in the vicinity of the data, underlining the geometrical accuracy of the reduction process, and thus the need of \red{fewer} master coordinates.

Shallow spherical shells have been investigated in~\cite{touze-shelltypeNL} and the method has been used to predict the correct type of nonlinearity for each mode of such structure\red{s} as a function of the curvature. FRFs for different type of shells (hyperbolic paraboloid panel, circular cylindrical panel and closed circular cylindrical shell), have been exhibited in~\cite{TOUZE:CMAME:2008}. Also, the transition to chaotic vibrations has been investigated with a  ROM composed of only the two modes in 1:1 resonance in~\cite{TOUZE:JFS:2007}, showing the limitation of the method (based on an asymptotic expansion) for very large amplitudes. Finally, applications of the method to FE structure\red{s} have been considered in~\cite{Touze:compmech:2014,Vizza3d}, but still taking the modal equations as a starting point. Direct computation of the normal form from the FE model will be discussed in Section~\ref{sec:DNF}.

Another interesting aspect of the normal form approach is to provide the simplest formulation of the reduced-order dynamics with only resonant monomials, thus opening the doors to \red{the} derivation of efficient {\em ex-nihilo} models~\cite{touze03-NNM,TouzeCISM}. In short, the normal form is the skeleton of the  dynamics and contains the correct qualitative picture and the same bifurcations \red{as} the full system. It is thus a powerful tool to understand the minimal models driving dynamical solutions and to build {\em ad-hoc} models containing the observed bifurcations. Important consequences are in the field of identification methods, where minimal nonlinear models can be used reliably, see {\em e.g.} circular plates with 1:1 internal resonances \cite{thomas03-JSV,Givois11,GivoisPiezo}, shallow shells with 1:1:2 and 1:2:2:4 internal resonances \cite{thomas07-ND,monteil15-AA,jossic18}, MEMS structure with 1:2 and 1:3 resonance~\cite{Gobatres12,ShawshaniAPL}, and  the identification of the hardening/softening behaviour of particular modes of a structure \cite{denis18-MSSP}.

The normal form approach has also been used by numerous other authors in the context of vibration, and the first introduction can be traced back to J{\'e}z{\'e}quel and Lamarque~\cite{Jezequel91}. The method has then be investigated by Nayfeh  who reduces it to a simple perturbation method~\cite{NayfehNF}, and by Leung and Zhang who developed close approaches~\cite{LEUNG1998,LEUNG1998b}. Higher-order approximations of normal transforms have also been developed using symbolic processors, see {\em e.g.}~\cite{ZHANG2000,Huseyin2000,LEUNG2003}, and application to plate vibration featuring 1:1 resonance is investigated in~\cite{YU2001}. More recently, it has been introduced for second-order vibratory systems, in a manner very similar to the presentation given in this section~\cite{NeildNF00,NeildNF01}, with in view the derivation of solutions for nonlinear vibration problem by using a single-harmonic assumption for the normal dynamics to derive analytical predictions. Also, only the first term in the normal form expansion was taken into account, leading to an incorrect prediction of the type of nonlinearity for systems with quadratic and cubic nonlinearity, as underlined in~\cite{BreunungHaller18}. The problem has then been corrected and the link to reduced-order models underlined in~\cite{Wagg2019}. Other contributions also tackled the problem of systems with periodic coefficients and/or periodic forcing, combining the Lyapunov-Floquet with a normal transform, see {\em e.g.}~\cite{SINHA1994687,SINHA07,Waswa}, or the computation of time-dependent normal form for handling the harmonic forcing~\cite{IoossForcing,Gabale009}. 

\subsection{Spectral submanifold}\label{subsec:SSM}

Spectral submanifolds (SSMs) have been first introduced by Haller and Ponsioen in~\cite{Haller2016} with the aim of emphasising the problem of existence and uniqueness  in the case of dissipative systems that had not been clearly elucidated in the previous works, as well as eliminating ambiguities in the terminology \red{being used} in the field of nonlinear normal modes. As underlined in Section~\ref{subsec:dynasys}, the problem of existence and uniqueness is tackled for conservative systems thanks to strong results by Lyapunov and Kelley~\cite{Lyapunov1907,Kelley3}, and LSM (Lyapunov subcenter manifolds) are known to be densely filled with periodic orbits under non-resonance conditions, thus naturally extending the linear modal subspaces. The picture is completely different for dissipative systems, with an immediate loss of uniqueness which has been underlined in different investigations~\cite{NeildNF01,CIRILLO2016,Haller2016}. In that case the structure of the phase space is dominated by strongly decaying modes leading to fast contraction of the flow. As a consequence, there are infinitely many invariant manifolds tangent to any subspace spanned by modes having small damping ratios (see~\cite{NeildNF01,Haller2016,CIRILLO2016} for simple illustrations on linear and nonlinear systems).
SSMs have been introduced in~\cite{Haller2016} with the aim of proving existence and uniqueness of the searched subspaces, using involved mathematical tools from the most recent development in dynamical system theory. They are defined as  the smoothest nonlinear continuation of a spectral subspace of the linearized system.  \red{From} the theoretical point of view, the understanding of the transition between conservative and dissipative structures is always delicate and challenging, and only very recent studies draw out the connection between LSM and SSM~\cite{Llave2019}.

The existence theorem for SSM given in~\cite{Haller2016} is linked to conditions on the regularity of the nonlinear vector field (that are fulfilled in our case of geometric nonlinearity since polynomial restoring forces are infinitely smooth) and non-resonance conditions on eigenfrequencies. Furthermore, the {\em spectral quotient} $\sigma_{\mathrm{out}}$ is defined as the integer part of the ratio between the largest damping rate of the slave modes to the smallest damping rate of the master modes. To be more specific, let us suppose that linear viscous damping of the form $\sigma_k \dot{x}_k$ is appended to each modal oscillator equation in Eq.~\eqref{eq:eom_modal}. Then the eigenvalues reads $\lambda_k =  - \sigma_k \pm i \omega_k \left( 1-\frac{\sigma_k^2}{\omega_k^2} \right)^{1/2}$. Assume one wants to construct the SSM associated to the first $d$ linear modes $(x_1, ..., x_d)$, with $d\ll N$. Then the spectral quotient reads:
\begin{equation}
\sigma_{\mathrm{out}} = \mathrm{Int} \left[ \frac{\mathrm{max}_{j=d+1, ..., N} \; |\sigma_j| }{\mathrm{min}_{p=1, ...,d} \; |\sigma_p| }  \right],
\end{equation}
where $\mathrm{Int}$ refers to the integer part. Existence, uniqueness and persistence of $d$-dimensional SSM are stated in~\cite{Haller2016}, under the general conditions given above. Furthermore, the SSM is unique among all other invariant manifolds of smoothness $\sigma_{\mathrm{out}}+1$ that share the same properties. In other words, uniqueness is reached only when the SSM can be computed by an asymptotic expansion which has an order at least equal to $\sigma_{\mathrm{out}}+1$. All prior developments, of lower order, are not unique and are only an approximation of the exact SSM. \red{In most real structures, it is commonly observed that the damping ratios of the modes are increasing with frequency. Consequently, the spectral quotient  is expected to be very large in a number of applications, and shall go to infinity in some  cases.} Importantly, lower-order truncations of the series will also approximate an infinity of other invariant manifolds which share locally similar properties. This very important result a posteriori justifies older developments shown in the previous sections, that give approximate low-order development of the searched LSM and/or SSM once the damping is taken into account.

In view of model order reduction, the computational procedure proposed in~\cite{Haller2016} is technically detailed in~\cite{PONSIOEN2018}, by restricting to the case of a single master mode (two-dimensional SSM). 
Note that the developments shown in the two previous sections were not relying on the formalism of the parametrisation method of invariant manifolds developed in~\cite{Cabre1,Cabre2,Cabre3,Haro}, which was not available at that time, but rather used more classical techniques proposed in center manifold theorem and homological equations for \red{the} normal form derivation. On the contrary, the computational scheme proposed in~\cite{PONSIOEN2018} closely follows the general guidelines of the parametrisation method as given in~\cite{Haro}. The mechanical equations of motion are set into the first-order, and Eq.~\eqref{eq:invariance} is rewritten for mechanical oscillatory systems. Asymptotic polynomial expansions are then introduced, and the tangent and normal cohomological equations are derived, while the {\em normal form style} is used to solve out the coefficients. The procedure has been automated and coded in the software SSMtool~\cite{PONSIOEN2018}. A special care is taken for the {\em near-inner resonance}, occurring for small damping values. Indeed, as stated in~\cite{TOUZE:JSV:2006,TouzeCISM,Szalai2017}, {\em trivial} resonances occur easily at third-order due to the particular shape of the eigenspectrum of a conservative vibratory system $\{ \pm i \omega_k\}_{k=1,...,N}$. When damping is added then these resonances are destroyed, however with the assumption of small damping it is important to retain these terms in the resulting dynamics (normal form) to avoid small divisors and linearisation. In such a case of near-inner resonance, the nearly resonant monomials are kept in the normal form (reduced dynamics) and the associated term in the nonlinear mapping are set to zero. Even though this procedure is named as using a {\em mixed style} in~\cite{PONSIOEN2018}, the calculation is equivalent to the one used in~\cite{TOUZE:JSV:2006,TouzeCISM} \red{which is} justified thanks to parameter-dependent normal form, and thus can be sorted as a normal form style.  The main advantage\red{s} of the computational procedure derived in~\cite{Haller2016,PONSIOEN2018} are: (i) to take easily into account dissipative forces and external harmonic forcing \red{by} resorting to a time-dependent SSM,  (ii) to offer an integrated and high-order solution since series expansions up to any order are technically possible. The only limitation to higher orders then rely in the memory requirements; in particular, computations up to order 15 are shown, also underlining that orders higher than 15 are generally too expensive in terms of memory consumption for standard computers.

In \cite{PONSIOEN2018}, the derivation is limited to a two-dimensional invariant manifold, meaning that extending the method to more than one master mode needs an extra calculation. The reduced-order dynamics is given in polar coordinates, thus providing a direct amplitude-frequency relationship that \red{does} not require the reduced model to be integrated in time. This is again an advantage since directly providing the coefficients of the development of the nonlinear amplitude-frequency relationship (backbone curve),  paid at the price of loosing oscillator-like equations. Whereas the developments shown in the two previous sections tried to fit the mathematical theories to the framework of mechanical systems, the point of view developed in~\cite{Haller2016,PONSIOEN2018} is  to fit the mechanical context into dynamical system formalism, with the important gain of more versatility, more generality and possibility of high-order expansions, up to converged results to large amplitude vibrations.

SSMs have already been applied to a number of different problems and contexts. A nonlinear Timoshenko beam is used as illustrative example in~\cite{PONSIOEN2018} while a linear Rayleigh beam resting on a cubic nonlinear foundation is tackled in~\cite{Kogel18}. Backbone curves and their relationship to frequency responses are investigated in~\cite{BreunungHaller18}, and the link with the identification problem using experimental data is illustrated in~\cite{Szalai2017}. Forced response calculations, involving \red{a} non-autonomous manifold, are shown on a linear Euler-Bernoulli beam with a nonlinear spring at its end in~\cite{PONSIOEN2020}. Time domain simulation\red{s} are reported on a von K{\`a}rm{\`a}n beam in~\cite{JAIN2018VK}. Finally, isolated solutions have also been predicted thanks to SSM~\cite{HallerIsola}.

\section{ROMs for finite element problems}\label{sec:ROMFE}

This section is specifically devoted to reduction methods for the case of geometrically nonlinear structures discretized with a finite element (FE) procedure. Four main reasons explain the need of having a dedicated section for this case. First, FE procedures are nowadays the most commonly used methods in engineering. Application of reduction methods to this class of problem is thus \red{of} specific interest,  since the potential applications are numerous. Second, the starting point is not given under the form of a PDE or under the modal form as in Eq.~\eqref{eq:eom_modal}. Instead, the starting point is semi-discrete equations in the physical space as in Eqs.~\eqref{eq:eom_phys}. Consequently reduction methods need to comply with this formalism. Third, numerous FE codes (commercial or  open source) exist and offer large capabilities in terms of computational power. Hence the idea of using the existing codes as such, without entering deeply inside their core, and use \red{of} classical features already developed, to produce a ROM, has led to the emergence of {\em non-intrusive} or {\em indirect} methods~\cite{mignolet13}. A fully non-intrusive method has a lot of inherent advantages in this respect \red{as it is} easily applicable to any existing FE code. Fourth, FE models generally use fine meshes involving millions of dofs, thus raising the curse of dimensionality. Consequently the methods need to be adapted to overcome this specific issue. In particular, using the modal basis as starting point is generally out of reach, and thus the methods presented in the previous sections need to be adapted.

\subsection{FE procedure and Stiffness Evaluation}\label{subsec:STEP}

The derivation of ROMs for FE structures featuring geometric nonlinearity is made difficult by the fact that existing codes don't give access to the quadratic and cubic terms, either in physical space, $\Gvec$ and $\Hvec$, or in modal space, $\g$ and $\h$. Non-intrusive methods started to develop in the early 2000's with the idea of using static computations of FE software to derive some of these coefficients, needed to build a ROM, since \red{the} access to $\kvec (\X) = \K \X+\fvec_\text{nl}(\X)$ is easily provided.  Two different methods have then been proposed: \red{the} first one where a prescribed displacement is imposed to the structure, and a second one where an imposed static force is applied. While the first method has then  been named as stiffness evaluation procedure (STEP), the second method give\red{s} rise to  implicit condensation which is fully detailed in Section~\ref{subsec:IC}.

The STEP has been first introduced by Muravyov and Rizzi in~\cite{muravyov}, with the aim of computing in a non-intrusive way, the nonlinear modal coupling coefficients $g^p_{ij}$ and $h^p_{ijk}$ appearing in Eq.~\eqref{eq:EDOmodalproj}.  The idea is to use a set of well-chosen prescribed displacements  as inputs for a static computation, a standard operation that is easily performed by any FE code. Then from the resulting deformed structure, a simple algebra  allows one to retrieve all the coefficients from the internal force vector given by the FE code, the key idea being to impose plus/minus the displacement with selected combinations of modes.

The method is fully explained in~\cite{muravyov}, here we illustrate the procedure by deriving the sole computation of coefficients $g_{pp}^k$ and $h_{ppp}^k$. The following static displacements are prescribed to the structure:
\begin{equation}
\label{eq:x1}
\Xvec_p =\pm \lambda \phivec_{p}\quad\Rightarrow\quad 
\left\{\begin{array}{l}
x_p=\lambda, \\
x_j=0\quad\forall j\neq p,
\end{array}\right.
\end{equation}
where $\lambda$ refers to an amplitude, the value of which has to be carefully selected (see {\em e.g.}~\cite{muravyov,givois2019} for discussions on this choice). Introducing Eq.~\eqref{eq:x1} into  Eqs.~\eqref{eq:eom_phys} and \eqref{eq:EDOmodalproj} leads to, for all $k=1,\ldots, N$:
\begin{subequations}
\label{eq:SysLinMur1Mode}
\begin{align}
\lambda^2 g_{pp}^k  + \lambda^3  h_{ppp}^k & = \tp{\phivec_k} \fvec_\text{nl} (\lambda\phivec_p)/m_k, \\ 
\lambda^2 g_{pp}^k  - \lambda^3  h_{ppp}^k & = \tp{\phivec_k} \fvec_\text{nl}(-\lambda \phivec_p)/m_k,
\end{align}
\end{subequations}
where $m_k$ is the modal mass (which could be unity with mass normalisation). The unknown quadratic and cubic coefficients are thus obtained by solving this linear system in $(g_{pp}^k,h_{ppp}^k)$, which depends on the computation of $\fvec_\text{nl}(\pm\lambda \phivec_p)$, which requires only the computation of reaction forces due to the prescribed displacements in the FE code, and no nonlinear Newton-Raphson procedure. 
Similar algebraic manipulations with more modes involved in the prescribed displacements then allows one to get the full family of quadratic and cubic coefficients, \red{which are solutions of other linear systems, not reported here for the sake of brevity.}  In addition to be non-intrusive, this procedure is then very time efficient since only linear operations are required, with no nonlinear system solving. 

At this stage it is important to underline that the STEP is not a reduction method, but only a non-intrusive  algebraic manipulation that can be used as a tool to derive some desired nonlinear characteristics from simple FE calculation. In its first development as given by~\cite{muravyov}, it was used to get access to the nonlinear modal coupling coefficients, meaning that the anticipated ROM that can be simply derived from that is the projection onto the linear modes basis, with all the known problems due to this projection (loss of invariance and nonlinear cross-coupling terms). One can note also that further developments used the procedure with different inputs, that are not necessarily combinations of modes, see {\em e.g.}~\cite{artDNF2020}. A key feature is related to the amplitude $\lambda$ the user has to select to obtain the modal coupling coefficients. As shown in a number of studies, see {\em e.g.}~\cite{muravyov,givois2019}, there exists a large range of $\lambda$ values conducting to stable values of the $g^p_{ij}$ and $h^p_{ijk}$, the main idea being that the imposed displacement needs to be not too small so as to correctly excite geometric nonlinearity, and not too large to stay in the range of moderate transformations (see~\cite{givois2019} for quantitative descriptions).

The STEP has been used in a number of different contexts, see {\em e.g.}~\cite{Mignolet08,KimCantilever,mignolet13,Lulf,LazarusThomas2012,givois2019,Vizza3d}. An important improvement in the computational complexity has been proposed in~\cite{Perez2014}, where the tangent stiffness matrix is used in order to divide by an order of magnitude the number of needed operations. It has also been proposed in~\cite{KimEstep} to compute the need\red{ed} quantities at the elementary levels, in order to speed up computations, however it seems intrusive. Finally, direct methods have also been proposed in order to compute intrusively the nonlinear modal coupling coefficients, see {\em e.g.}~\cite{Touze:compmech:2014} for a direct computation on MITC (mixed interpolation of tensorial components) shell elements, \cite{senechal2011} for 3D FE and~\cite{Dou2015} with beam elements and application to shape optimization.

\subsection{Implicit condensation and stress manifold}\label{subsec:IC}

The implicit condensation and expansion (ICE) method has been first introduced in a series of papers from two different groups. The first developments date back from the PhD thesis by Matthew Mc Ewan defended in 2001~\cite{McEwanThesis,McEwan01}, then continued and improved by Hollkamp and Gordon~\cite{Hollkamp2005,Hollkamp2008}, who introduced the acronym ICE for the method. Recently, it has been further investigated and used by Kuether, Allen {\em et al.}~\cite{KUETHER2014,kuether2015}, as well as by Frangi and Gobat~\cite{FRANGI2019}, who introduced the term {\em stress manifold} to describe the reduction subspace used to build the ROM. 
The method realizes an {\em implicit} condensation of the non-modeled degrees of freedom, and is shown to be fully equivalent to the classical {\em static condensation} if all equations are fully known~\cite{Hollkamp2008,mignolet13,YichangICE}, implying that it can never perform better than a static condensation. 

The ICE method is non-intrusive in nature and relies on two different steps. The first one can be realized by any FE code since using a standard procedure of static nonlinear computation. A series of body forces $\fvec_\text{e}$ that are proportional to the inertia of the linear modes, $\fvec_\text{e} = \beta_i \Mvec \phivec_i$ with $\beta_i\in\mathbb{R}$, are imposed to the structure, where $i=1...m$, $m$ being the number of selected master modes. The resulting structural deformation $\Xvec$ is computed in statics by the FE software and back-projected onto the eigenmodes in order to retrieve each modal displacement $x_j$, $j=1...m$. Since nonlinear couplings are present, the resulting $x_i$ are not directly proportional to the forcing. Instead, slave modes are excited through invariant-breaking terms, such that the resulting non-linear displacements of the master mode follows a stress manifold that implicitly realize\red{s} the condensation of the non-modeled coordinates. Appendix~\ref{app:static} gives a few more technical details on this computation and better highlights the link with explicit static condensation. A full mapping is constructed from this series of computation, with entries $\beta_i$ and outputs $x_j$, describing the stress manifold.  The second step is a fitting procedure that has to be undertaken, assuming the mapping $x_j(\beta_i)$ is invertible. Hence from the computed clouds of points, functional forms can describe the resulting nonlinear restoring force spanning the stress manifold.

In many applications, the fitting procedure is realized thanks to polynomial expansions. The versatility of the method can here lead to using higher-order polynomials or other test functions ({\em e.g. splines}), with the aim of getting more accurate results on a larger span of displacements. From the development of the method, one also easily understands that the fitted coefficients depend on the amplitudes of the scaling factors $\beta_i$ used to construct the stress manifold, as observed in~\cite{Hill00}. Indeed, this reduction subspace being curved, different polynomial fittings are obtained when varying the amplitude~\cite{YichangICE}. This is in contrast to the STEP, where prescribed displacements are used to compute the nonlinear modal coupling coefficients. Since the modal eigenspaces are straight planes, the STEP coefficients are constant on a large range of applied displacements~\cite{givois2019}.

\begin{figure}[ht!]
\pgfplotsset{compat=1.4}
\centering
\includegraphics[width=.33\linewidth]{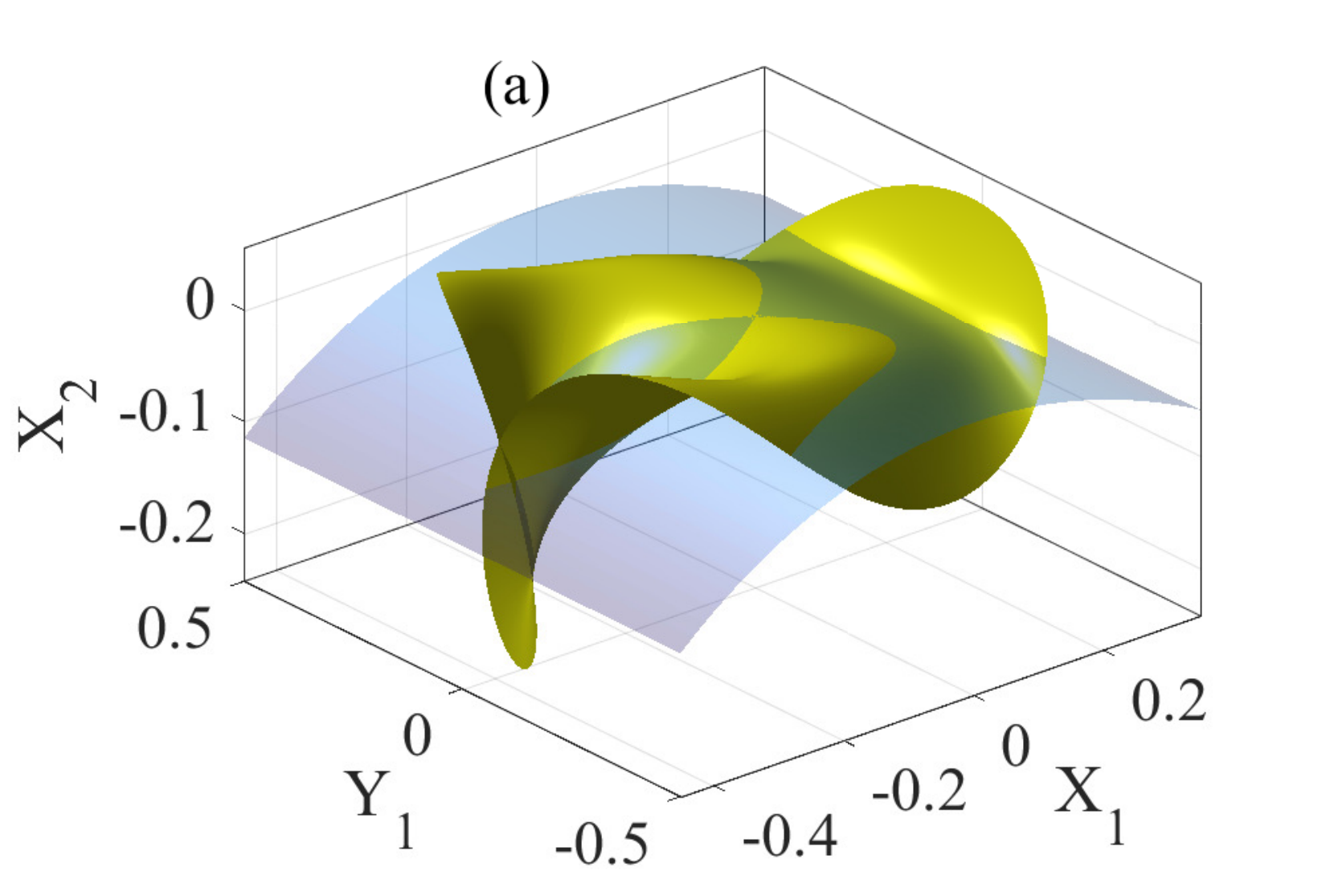}
$\quad$
\includegraphics[width=.33\linewidth]{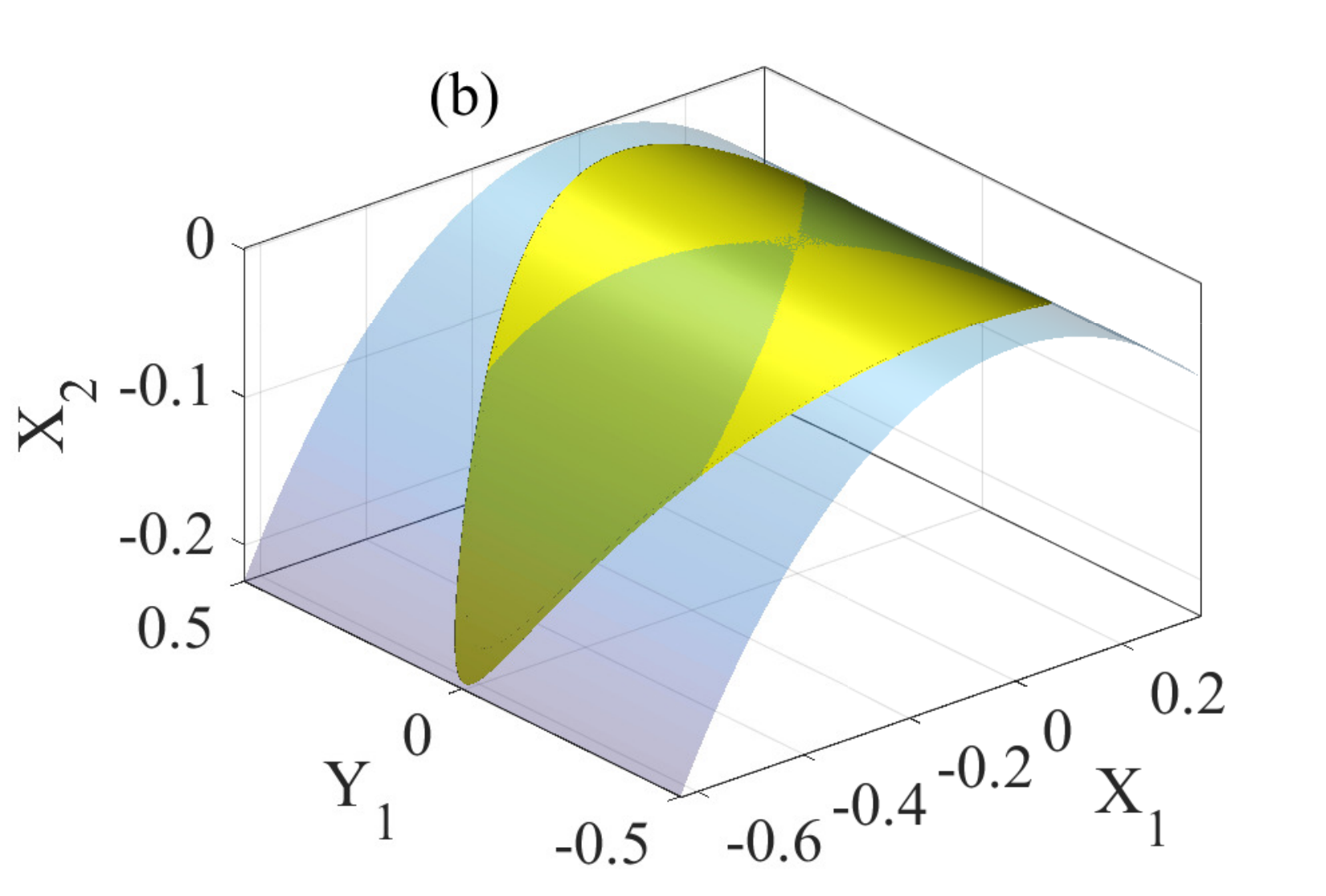}
$\quad$
\begin{tikzpicture}
\begin{axis}[
	unbounded coords=jump,
	width=.3\textwidth,
	xmin=0,xmax=6,ymax=3,ymin=-3,
	xlabel=$\rho$,
	ylabel=$\mathcal{R}$,
	grid=major,
	xtick={0,1,2,3,4,5,6},
	ytick={-3,-2,-1,0,1,2,3}
]
  \addplot[black,thick,samples=121,domain=0:6] (x,{(x*x-8/3)/(x*x-4)});
\end{axis}
\node[below right]
    at (current bounding box.north west) {\scriptsize{(c)}};
\end{tikzpicture}
\captionsetup{font={small}}
\caption{Comparison of stress manifold obtained with static condensation (light blue), and invariant manifold obtained from numerical continuation of periodic orbits (yellow) in phase space for a two-dof system of coupled nonlinear oscillators. In each figure, $\omega_{1}=1$ while $\omega_2$ is increased to meet the slow/fast assumption. (a) $\omega_{2}=2.5$, (b)  $\omega_{2}=10$. (c) Correction factor $\mathcal R$ as function of $\rho=\omega_2/\omega_1$. \red{Figure reworked from~\cite{YichangICE}.}}
\label{fig:manifolds00}
\end{figure}
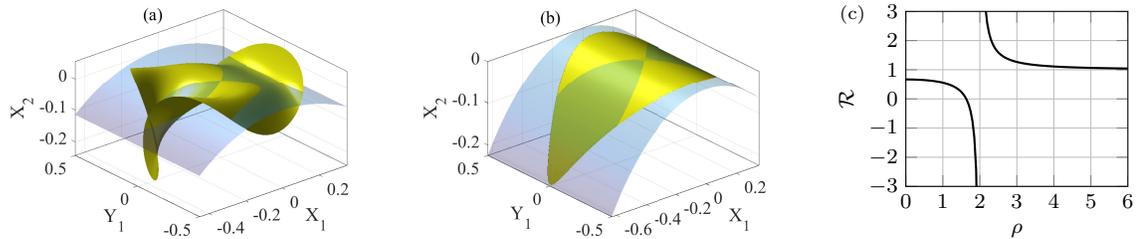

The ICE method has been compared to invariant manifolds in~\cite{HallerSF,YichangICE}. First, the stress manifold is not an invariant subspace. Second, the construction method is {\em static} in nature such that the resulting stress manifold does not depend on the velocities. This has important consequences on the accuracy and is for example illustrated in Fig.~\ref{fig:manifolds00}, which compares the stress manifold and the actual invariant manifold (IM) for a two-dofs nonlinear system, where the IM has been computed numerically with a continuation method. Two figures are shown, realized for two different values of the two selected eigenfrequencies $\omega_1$ and $\omega_2$ of the system. While $\omega_1=1$ for both figures, $\omega_2=$2.5 for Fig.~\ref{fig:manifolds00}(a) and $\omega_2=$10 for Fig.~\ref{fig:manifolds00}(b). In the first case, one can observe very important differences between the stress and the invariant manifolds, underlining the fact that neglecting the velocities leads to important mistakes, and that the ICE method produces a very simplistic manifold that is far from the correct reduction subspace. On the other hand when $\omega_2=$10, both manifolds tend to the same geometry.

This example illustrates the need of fulfilling a slow/fast assumption for the ICE method to produce correct results. The {\em slow/fast partition} refers to the fact that the eigenfrequencies of the slave modes are very large  as compared to those of the master modes. When this is verified, then
from the theorem demonstrated in~\cite{HallerSF} and the analytical and numerical calculations shown in~\cite{YichangICE}, the results given by ICE method are reliable and the stress manifold tends to the invariant manifold. On the other hand, incorrect predictions will be given by \red{the} ICE method if  \red{the} slow/fast assumption is not fulfilled.

Let us exemplify this result on the type of nonlinearity (hardening/softening behaviour). Assuming $m$ is the master mode, the reduced dynamics on the stress manifold, up to cubic order, can thus be written as~\cite{YichangICE}:
\begin{equation}\label{eq:ICEdyna}
\ddot{x}_m + \omega_m^2 X_m + g_{mm}^m X_m^2+\left( h_{mmm}^m   - \underset{s\neq m}{\sum_{s=1}^N}  2\;\frac{g_{ms}^m g_{mm}^s}{\omega_s^2} \right) x_m^3+{\mathcal{O}}(x_m^4)=0,
\end{equation}
where this equation has been obtained assuming all the coefficients of the model are known and applying explicit static condensation. Using the same notation as in Section~\ref{sec:NFmethod}, {\em i.e.} $\omega_{NL} = \omega_m (1 + \Gamma^{\mathrm{ICE}}_m a^2)$ for the amplitude-frequency relationship, the predicted type of nonlinearity given by condensation \red{can be written as}:
\begin{equation}\label{eq:GammaGeneSC}
\Gamma^{\mathrm{ICE}}_m= - \dfrac{5}{12\;\omega^2_{m}} \left(\dfrac{g^m_{mm}}{\omega_{m}}\right)^2  + \dfrac{3}{8\;\omega^2_{m}}\left( h_{mmm}^m - \sum_{\underset{s\neq m}{s=1}}^N  2\left(\dfrac{g^s_{mm}}{\omega_{s}}\right)^2 \right). 
\end{equation}
One can now compare the prediction given by  Eq.~\eqref{eq:GammaGeneSC} to that given by following the actual family of periodic orbits lying in the invariant manifold of the system, Eq.~\eqref{eq:gammaNF}, and observe that \red{the} first terms are exactly the same, the only difference being in the last summed terms. One can then construct the ratio $\mathcal R$ of these last summed term, denoted as {\em correction factors} in~\cite{VizzaMDNNM,Vizza3d,YichangICE}, thus exactly comparing the difference between the two approaches, which simply reads:
\begin{equation}\label{eq:corrfacR}
\mathcal{R} = \frac{\omega_s^2-\frac{8}{3}\omega_m^2}{\omega_s^2-4\omega_m^2} = \frac{\rho^2-\frac{8}{3}}{\rho^2-4},
\end{equation}
where it has been assumed for the sake of simplicity that a single slave mode $s$ exists. This correction factor has been expressed as a function of $\rho=\omega_s/\omega_m$, and interestingly it does not depend either on quadratic or on cubic coefficients. The behaviour of $\mathcal R$ as \red{a} function of $\rho$ is reported in Fig.~\ref{fig:manifolds00}(c)~\cite{YichangICE}. One can observe that the divergence when $\rho=2$ is not given by the ICE method. The divergence of the type of nonlinearity in \red{the} case of 2:1 internal resonance is a known effect already reported in Section~\ref{sec:IMSP}. In this area, a strong coupling exists and the reduction to a single master mode is meaningless, see {\em e.g.}~\cite{Arafat03,regacarbo00,touze-shelltypeNL,YichangICE,Lenci12,Haro} for more discussions on this subject. More importantly, one can clearly observe that $\mathcal R$ is tending to 1 when $\rho$ increases, showing that the prediction given by \red{the} ICE method is correct only when the slow/fast assumption is fulfilled. From Eq.~\eqref{eq:corrfacR}, a 1\% error on the type of nonlinearity is predicted when $\rho=$11.7 and a 10\% error when $\rho=$4.15. The incorrect prediction of the ICE method on the type of nonlinearity has been illustrated in~\cite{YichangVib} in the case of a linear beam resting on a nonlinear elastic foundation.

As a conclusion on the ICE method, one can note the important advantage of being fully non-intrusive and simple to implement on any FE code. \red{From a } theoretical viewpoint, the main  benefit of the method is to realize an implicit condensation of the non-modeled dofs. This is particularly meaningful for planar structures where the eigenspectrum between bending and in-plane modes is particularly well separated, but will become an obstacle when dealing with curved structures such as arches and shells. Also the method appears particularly appealing when reducing to a single mode, since allowing easily a higher-order polynomial fitting, the method could be able to follow backbone curves up to larger amplitudes than methods limited to third-order expansions. On the other hand, the main drawbacks of the method are that no velocity is taken into account when building the reduced dynamics, and a slow/fast assumption is mandatorily needed for providing accurate predictions, otherwise incorrect results are provided. The last drawback is connected to the fitting procedure and the construction methods in cases where the number of master modes becomes larger. As underlined in~\cite{YichangICE}, with more than one master coordinate, the method becomes very sensitive to the choice of the load scale factors  $\beta_i$, leading to a lack of robustness and a too strong dependence on small variations of inputs. 

A last known drawback of the method relies \red{on} its inability to properly take into account cases where inertia nonlinearity is important, as underlined for example \red{in} the case of the cantilever beam~\cite{KimCantilever,YichangVib,NicolaidouIceKE}, or  a micromirror in~\cite{AndreaROM}. A method to bring a correction for this specific case has been implemented recently in~\cite{NicolaidouIceKE}. In particular,  the technique proposed in~\cite{NicolaidouIceKE} is equivalent to using the quadratic
manifold approach with static modal derivatives orthogonal to the master modes if stopped \red{at} the second order. However, the derivation allows an easy computation of higher orders, beyond the second-order term contained in the quadratic manifold detailed in the next section. Nevertheless, the method still relies on the slow/fast separation of modes to deliver accurate predictions.

Finally, one can also note that investigations have tried to combine STEP and ICE method. For example, dual modes as introduced in~\cite{KIM2013,mignolet13,PEREZnotch,WangMigno18,WANG20181,WangLocal19} combine a set of bending modes (with nonlinear coefficients given by the STEP), to added {\em dual} modes obtained from static imposed body forces --as in the first step of the ICE method-- but then analysed with an SVD/POD method to determine the most important patterns. Results from dual modes and modal derivatives covered in the next section are reported in~\cite{KARAMOOZ,WANG2021}. Also a modified STEP (M-STEP) proposed in~\cite{mehner2000,Vizza3d}, selects only a subset of master nodes of the FE mesh to apply the prescribed displacements, letting the other free, so as to implicitely condense their nonlinear relationship to the masters. Application of M-STEP to symmetric structures such as beams and plates but also symmetric laminated panels with piezoelectric patches, have shown good results thanks to the fulfillement of the slow/fast assumption~\cite{Vizza3d,givois21-CS}. 
In the same trend, static condensation is also applied in~\cite{WANGPalacios2015} to reduce the information of a 3-D model of a slender structure to a 1-D equivalent beam model.

\subsection{Modal derivatives and quadratic manifold}\label{subsec:MDQM}

Modal derivatives (MD) have been first introduced by Idehlson and Cardona to solve structural vibrations problems with a nonlinear stiffness matrix~\cite{IDELSOHN1985,IDELSOHN1985b}, with the key idea of taking into account the amplitude dependence of mode shapes and eigenfrequencies. It has been generally used in a number of different context as added vectors that can be appended to the projection basis in order to enrich the representation and better take into account nonlinear effects~\cite{Slaats95,Tiso2011,Weeger2014,Weeger2016,Wu2016,SOMBROEK2018,Wu2019,KARAMOOZ}. Following~\cite{IDELSOHN1985,Weeger2016,VizzaMDNNM}, let us denote as $\tilde{\phivec}_i(\Xvec)$ this amplitude-dependent vector, such that for small amplitude one retrieve the usual eigenmode: $\left[ \tilde{\phivec}_i(\Xvec)\right]_{\Xvec=\zervec} = \phivec_i$. The $ij$-th MD $\vec{\Theta}_{ij}$ is the derivative of $\tilde{\phivec}_i(\Xvec)$ with respect to a displacement enforced along the direction of the $j$-th eigenvector $\phivec_j$:
\begin{equation}\label{eq:MDdef}
\vec{\Theta}_{ij} 
\doteq \frac{\partial \tilde{\phivec}_i (\Xvec)}{\partial x_j}
\bigg\rvert_{\Xvec=\vec{0}} \, .
\end{equation}

In order to arrive at a computable definition of the MD that can also be used in  a simulation-free context, one can first derive the eigenvalue problem with respect to amplitude. The first term retrieves the standard linear eigenvalue problem while the second one makes appear the MD~\cite{Weeger2016,VizzaMDNNM}. However the derivative of the eigenfrequency with respect to amplitude also comes into play, so that the problem becomes underdeterminate. To close the system, the usual procedure is to derive the same Taylor expansion on the mass normalization equation, leading to:
\begin{equation}
\begin{bmatrix}
\Kvec -\omega_i^2 \Mvec
&\;
-\Mvec
\phivec_i
\\[10pt]
-\phivec_i^T 
\Mvec
&
\; 0
\end{bmatrix}
\left\lbrace
\begin{matrix}
\vec{\Theta}_{ij}
\\[10pt]
\frac{\partial \omega_i^2}{\partial x_j}
\end{matrix}\right\rbrace
=
\left\lbrace
\begin{matrix}
-2\Gvec(\phivec_j, \phivec_i)
\\[15pt]
0
\end{matrix}\right\rbrace,
\label{eq:MD_sys}
\end{equation}
where it is possible to see that the introduction of the mass normalisation equation coincides with a constraint on the MD to be mass-orthogonal to the $i$-th mode, thus rendering the system solvable.
However in most of the studies, a simplification of this formulation is used by introducing the {\em static modal derivative} (SMD) $\vec{\Theta}^{(S)}_{ij}$ by neglecting the terms related to the mass matrix:
\begin{equation}
\Kvec \vec{\Theta}^{(S)}_{ij}=-2\Gvec (\phivec_j, \phivec_i).
\label{eq:SMD_sys}
\end{equation}
The main argument for this simplification resides in the fact that operation~\eqref{eq:SMD_sys} can be easily implemented in a non-intrusive manner in any FE software, which is not the case for \eqref{eq:MD_sys}. 

Since MDs and SMDs are directly linked to the quadratic terms of the nonlinear restoring force, the natural extension of the method is to define a nonlinear mapping where the linear part is conveyed by the standard eigenvectors while the quadratic part is expanded on the modal derivatives. This idea leads to the {\em quadratic manifold} (QM) approach developed in~\cite{Jain2017,Rutzmoser}, where a second-order nonlinear mapping is defined between master coordinates gathered in a $\xvec$ vector of small dimension, and the original physical coordinates, reading: 
\begin{equation}\label{eq:QMNLmap}
\Xvec = \vec{\Phi} \xvec +\frac{1}{2}\vec{\Theta}(\xvec,\xvec) = \sum^m_{i=1}\phivec_i x_i + \frac{1}{2}\sum^m_{i=1}\sum^m_{j=1} \bar{\vec{\Theta}}_{ij} x_i x_j,
\end{equation}
where $\bar{\vec{\Theta}}_{ij}=(\vec{\Theta}_{ij}+\vec{\Theta}_{ji})/2$ is the symmetrized MDs, which can be selected either as a full or a static MD. One can note in particular that the linear part is an expansion onto the master linear modes of interest, while the quadratic part takes into account the nonlinear dependence with amplitude through a spanning of the phase space expressed by the modal derivatives. The reduced-order model is obtained by deriving Eq.~\eqref{eq:QMNLmap} \red{twice} with respect to time and applying a standard Galerkin projection, see {\em e.g.}~\cite{Jain2017,Rutzmoser} for technical details and~\cite{VizzaMDNNM,YichangVib} for developed indicial expressions up to order three.

\begin{figure*}[h!]
\hspace{.4cm}~
\includegraphics[scale=.275]{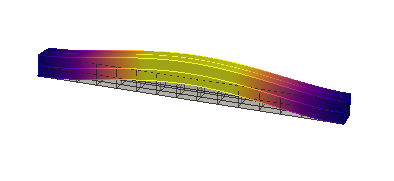}~
\includegraphics[scale=.275]{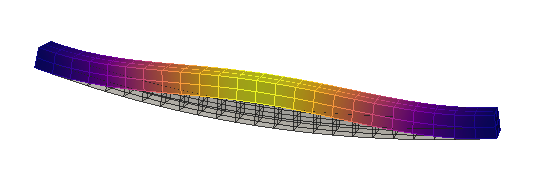}~
\includegraphics[scale=.275]{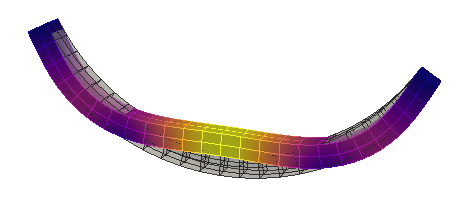}
\\
\includegraphics[width=.323\textwidth]{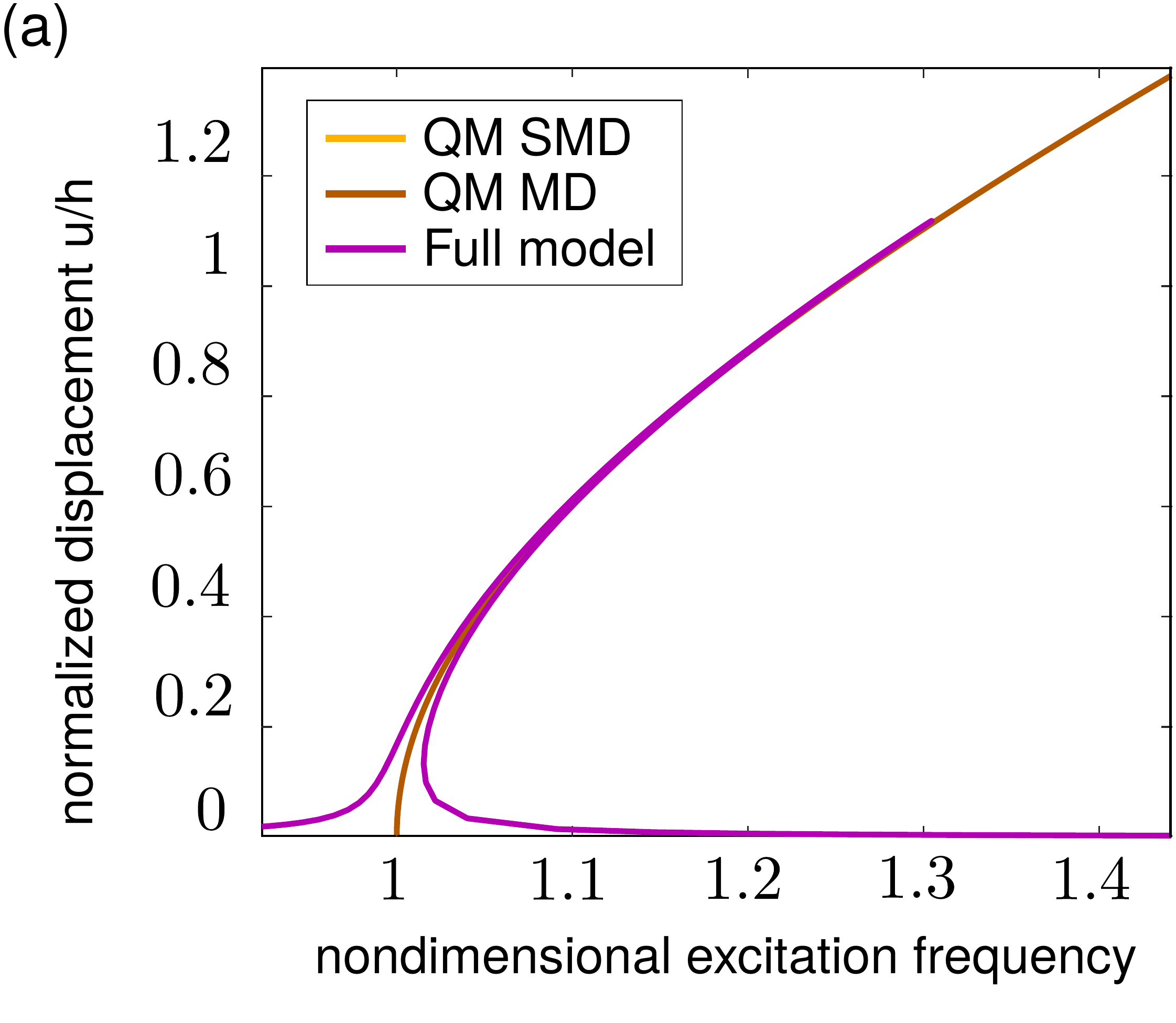}
\includegraphics[width=.323\textwidth]{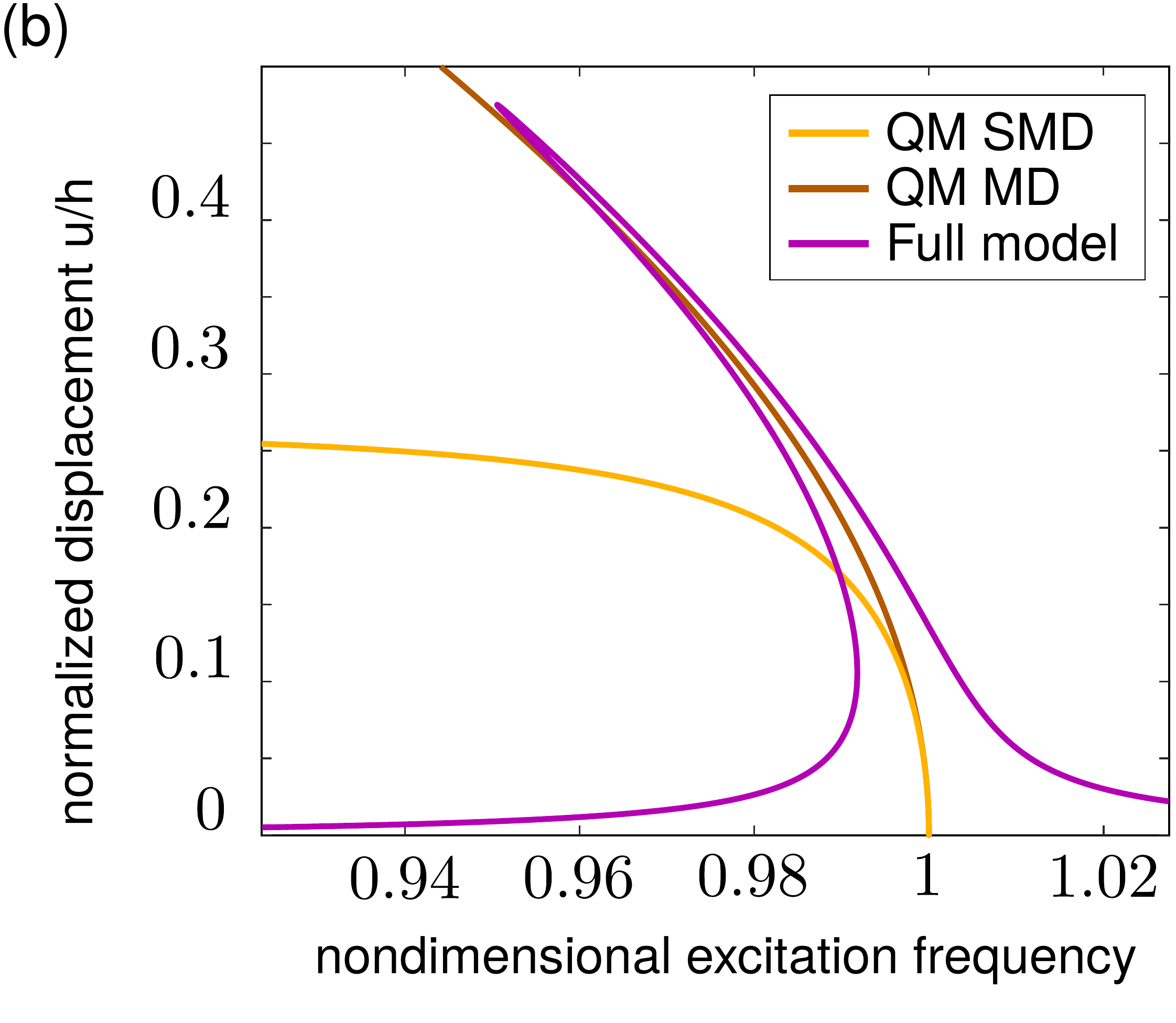}
\includegraphics[width=.323\textwidth]{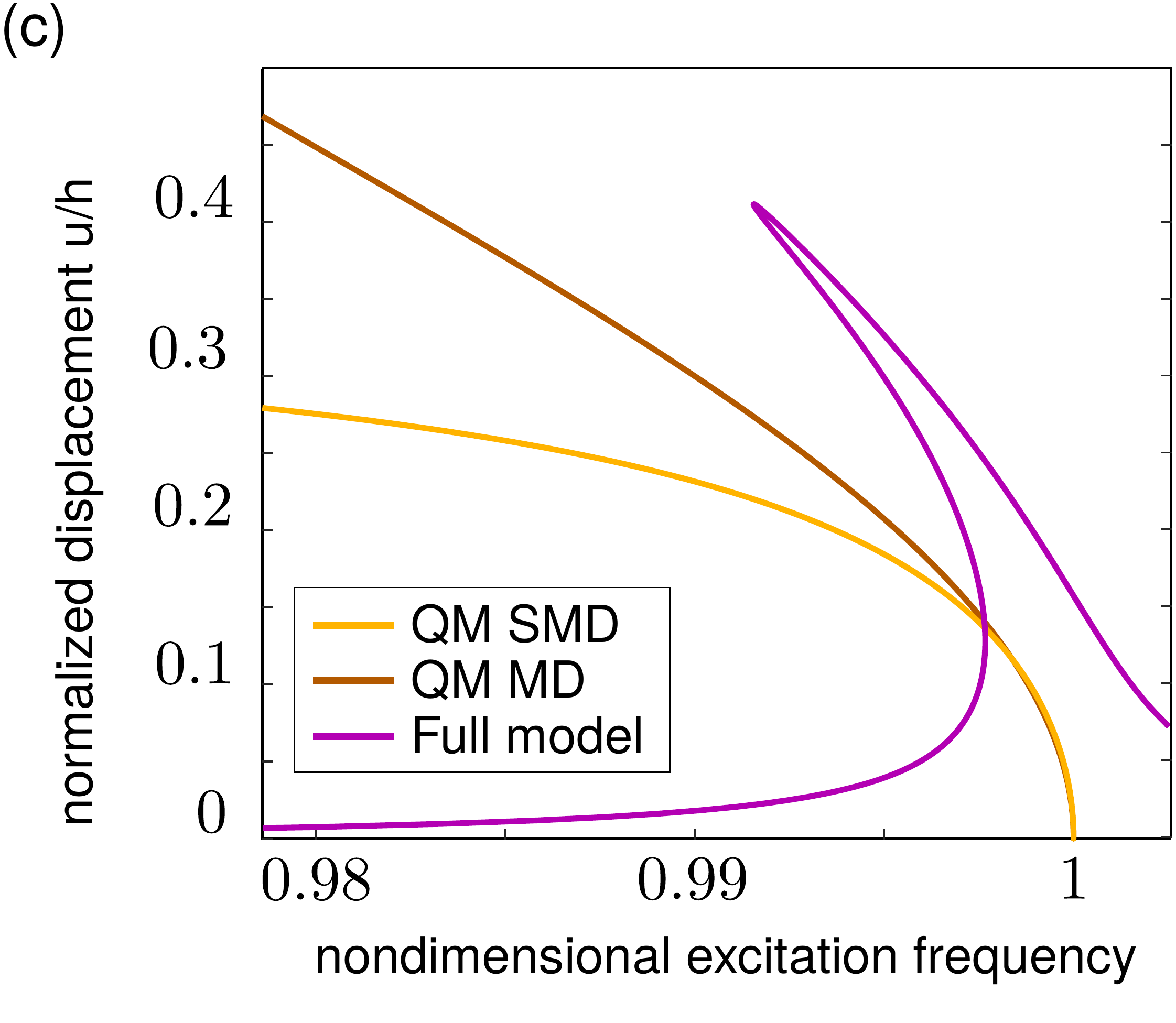}
\caption{Comparison of backbone curves obtained from QM with MDs (dark orange) and SMDs (yellow), for the three tested structures: (a) flat beam, (b) shallow arch, (c) non-shallow arch. Nondimensional amplitude of bending displacement (along $y$, nondimensionalised with respect to the thickness) of the central node of each beam as a function of $\omega/\omega_1$ where $\omega_1$ refers to the  eigenfrequency of the first bending mode. The backbone curves are contrasted to the FRF obtained on the full system (Full model, violet) with numerical continuation and a small amount of damping. Figure \red{reworked} from~\cite{VizzaMDNNM}.}
\label{fig:bbfrfbeams}
\end{figure*}

Figure~\ref{fig:bbfrfbeams} illustrates how the predictions given by QM-MD and QM-SMD may depart from the correct result given by the full model if the method is strictly used, with a single master coordinate, and truncating the dynamics up to third order. A clamped-clamped beam with increasing curvature is selected. For the flat beam case, both QM-MD and SMD allow predicting accurately the nonlinear response. On the other hand, for the second case with slight curvature, the SMD method departs from the correct prediction. Increasing again the curvature to arrive at a non-shallow arch, then both methods are not able anymore to reproduce the correct softening behaviour with a single master coordinate. \red{Note that incorrect predictions given by other SDOF reductions have also been reported before for buckled beams and truncation to single linear mode, see {\em e.g.}~\cite{nayfehcarbo97,LacarboBuck98,LACARBONARA1999}.}

The quadratic manifold approach with MD has been fully compared to the normal form approach in~\cite{VizzaMDNNM}, while rigorous theorems are provided in~\cite{HallerSF} in order to assess the merits and drawbacks of the method as compared to invariant manifold based techniques. To summarize the main findings, one can first note that the MDs and the associated nonlinear mapping (QM) do not take into account the velocity from the beginning of the development, either in Eq.~\eqref{eq:MDdef} or in \eqref{eq:QMNLmap}. This has important consequences and leads to the fact that the QM as derived in Eq.~\eqref{eq:QMNLmap} is independent from the velocity in phase space, and is thus not an invariant manifold. The second consequence is that the method needs again a slow/fast assumption between slave and master coordinates in order to predict correct results, as demonstrated in~\cite{HallerSF} and illustrated in~\cite{VizzaMDNNM} in several examples. In particular, phase space comparisons of the geometry of QM and IM are reported in~\cite{VizzaMDNNM}, illustrating the differences between reduction subspaces, when slow/fast assumption is not met. 

In order to give more quantitative understanding of the predictions given by the QM approach, the type of nonlinearity can be computed from the reduced dynamics.  Using \red{the} QM method with either MD or SMD leads to coefficients $\Gamma_\text{MD}$ and $\Gamma_\text{SMD}$ as~\cite{VizzaMDNNM}
\begin{subequations}\label{eq:zifullgamma}
\begin{align}
\Gamma_\text{MD}=&
-
\dfrac{5}{12\;\omega^2_{m}}
\left(\dfrac{g^m_{mm}}{\omega_{m}}\right)^2
+
\dfrac{3}{8\;\omega^2_{m}}\left(
h_{mmm}^m
-
\sum_{\underset{s\neq m}{s=1}}^n
2\left(\dfrac{g^s_{mm}}{\omega_{s}}\right)^2
\left(
1+
\dfrac{\omega^2_m(4\,\omega^2_s-3\,\omega^2_m)}{3(\omega^2_s-\omega^2_m)^2}
\right)\right), \label{eq:zifullgammaA}
\\
\Gamma_\text{SMD}=&
-
\dfrac{5}{12\;\omega^2_{m}}
\left(\dfrac{g^m_{mm}}{\omega_{m}}\right)^2
+
\dfrac{3}{8\;\omega^2_{m}}\left(
h_{mmm}^m
-
\sum_{\underset{s\neq m}{s=1}}^n
2\left(\dfrac{g^s_{mm}}{\omega_{s}}\right)^2
\left(
1+\dfrac{4\,\omega^2_m}{3\,\omega^2_s}
\right)\right). \label{eq:zifullgammaB}
\end{align}
\end{subequations}
Comparing these two formulas with the correct prediction given by following the family of periodic orbits foliating the associated invariant manifold, Eq.~\eqref{eq:gammaNF}, leads to the quantitative conclusion that QM can be used if the ratio between slave and master eigenfrequencies is larger than 4~\cite{VizzaMDNNM}. For flat structures where the separation between bending modes and in-plane modes is clear this is not problematic and MD can be used safely. On the other hand for curved beams or shell problems where more couplings between bending modes exist and where the slow/fast separation has no reason to be fulfilled, this could be problematic and lead to erroneous predictions. Illustrations on a simple analytical model (linear beam resting on nonlinear elastic foundation) are given in~\cite{YichangVib}, while the case of shallow spherical shells is investigated in~\cite{YichangNODYCON}.

Further developments on the amplitude of the harmonics of the solution also underlines that the QM with SMD leads to an incorrect treatment of the ``self-quadratic'' term $g^m_{mm} X_m^2$ with $m$ being the master mode (see~\cite{VizzaMDNNM} for details). Consequently when this term is present, which arises for curved structures such as arches or shells, the QM-SMD method cannot produce a correct prediction of the backbone curve, even if the slow/fast assumption is verified.

All these features are illustrated in Fig.~\ref{fig:bbfrfbeams}, taking as \red{an} illustrative example a clamped-clamped beam with increasing curvature. For the flat beam case, QM methods (either with MD and SMD) gives excellent results due to the fulfillment of the slow/fast assumption and the absence of quadratic couplings between bending modes. When a slight curvature is added, leading to a shallow arch with quadratic couplings, even though the slow/fast assumption is met, one can observe that QM-SMD method departs from the exact solution, due to the incorrect treatment of the self-quadratic term. Finally, increasing again the curvature, the slow/fast assumption is not met anymore, and both QM methods fail in predicting the correct type of nonlinearity. 

Modal derivatives have been applied to a number of problems, and we can mention in particular the latest contributions, apart from the ones already cited at the beginning of this section. Beams and panels are considered in~\cite{KARAMOOZ}, while the computation of backbone curves up to large amplitudes are addressed in~\cite{SOMBROEK2018}, underlining that numerous MDs need to be added as vector basis in order to catch internal resonance in particular. Transient analysis are used in~\cite{JainTiso2018} to enrich a quadratic manifold with linear modes and MDs. Finally shape imperfections and defects are taken into account in~\cite{MARCONI2020,MARCONI2021}, and comparisons with dual modes are commented in~\cite{KARAMOOZ,WANG2021}.

As a short conclusion on modal derivatives, one can first note the number of advantages it offers: it is a simulation-free approach, it is non-intrusive in nature and can be used easily from any FE code. In most of the applications, modal derivatives are used as added basis vectors allowing for a proficient method able to report involved nonlinear phenomena~\cite{SOMBROEK2018}. The framework of QM, while appearing as logical and attractive, shows however some inherent limitations. Most of the limitations (need of slow/fast assumption and incorrect treatment of self-quadratic term for SMD) fundamentally rel\red{y} on the fact that the method does not assume any dependence with the velocity, which is a strong drawback since in vibration theory velocities are independent variables and are needed in a phase space perspective in order to construct periodic orbits representing the oscillations. Interestingly, if one wants to insert the velocity dependence from the beginning, then the found nonlinear mapping is exactly the one derived from normal form theory exposed in section~\ref{sec:NFmethod}. This is further illustrated in the next section where a direct computation of the normal form is shown,  in order to apply the method from the FE discretization.

\subsection{Direct computation of normal form}\label{sec:DNF}

The direct normal form (DNF)\footnote{\red{One can note that the terminology DNF for direct normal form has been introduced before in~\cite{Elliott2018}, but with a different purpose as the one used here. In~\cite{Elliott2018}, "direct" is used to specify that the normal form is computed from the second-order oscillator equations, without using first-order, state-space formulation.}} approach has been first introduced in~\cite{artDNF2020}. The main idea is to propose a computational scheme that does not need the full modal basis calculation as a starting point, so that a direct nonlinear mapping from the physical space (dofs of the original problem) to the reduced subspace given by invariant manifolds, is retrieved. In a FE context where the number of dofs can be extremely large and attain millions, application of the formulas given in section~\ref{sec:NFmethod} \red{is} indeed inoperable, since the cost of diagonalizing the linear part is far too expensive. 
Consequently, adapting the method and giving a direct computation of the normal form is an important improvement, allowing to derive efficient reduced-order models. 

In the first development of the method reported in~\cite{artDNF2020}, real-valued mappings are used and the calculations followed closely the general guidelines given in~\cite{touze03-NNM,TOUZE:JSV:2006,TouzeCISM}, with the main difference that all the formulas are rewritten from the physical space. As reported in Section~\ref{sec:NFmodal}, a complete change of coordinate with summations up to all the dofs of the system is first derived, then followed by a truncation to retain the selected master modes. The advantage of proceeding like this is that the calculation is done once and for all, whatever the number of master coordinates to be retained.

Damping and external forcing are also considered in~\cite{artDNF2020}, following the general guidelines given in~\cite{touze03-NNM,TOUZE:JSV:2006,TouzeCISM}. Modal forcing is added to the reduced dynamical equations. For the damping, the general formulas given in~\cite{TOUZE:JSV:2006} have been adapted to the case of FE models in~\cite{artDNF2020} by (i) assuming small damping for the master modes so that only the first-order terms (in damping) of the general formulas given in~\cite{TOUZE:JSV:2006} are used, (ii) assuming Rayleigh damping as a specific input in order to comply with the general assumptions used to model losses in most of the FE codes. Of course, this specific damping is not restrictive and one can come back to the general formulas given for any modal damping in~\cite{TOUZE:JSV:2006}.

The DNF has then been completely rewritten in~\cite{AndreaROM}, using complex-valued formalism inherited from using first-order, state-space formulation from the beginning of the calculation. The main advantages of the complete rewriting is to use symmetric formulations throughout the derivation, while keeping clear the mechanical context by always referring, in all the developments, to mass and stiffness matrix $\Mvec$ and $\Kvec$, quadratic and cubic terms $\Gvec$ and $\Hvec$.  The homological equations have then been fully rewritten with these notations, opening the doors to numerous extensions of the method in order to go to higher orders or to add new physical phenomena described by new forces. The link with real-valued mappings derived in~\cite{artDNF2020} is completely detailed. Finally, the treatment of second-order internal resonances is highlighted, a development that had not been tackled before in~\cite{touze03-NNM,TOUZE:JSV:2006}. 

Let us describe with a few equations some of the main features of the method in the conservative framework, and how it compares to previous developments and in particular to modal derivatives described in Section~\ref{subsec:MDQM}. The starting point is a nonlinear change of coordinates  between the initial displacement-velocity vectors $\Xvec$, $\Yvec=\dot{\Xvec}$, and the {\em normal} variables  $(\R,\S)$ as:
\begin{subequations}
	\begin{align}
		\Xvec = \hat{\Psivec}(\R,\S), \\
		\Yvec = \hat{\Upsvec}(\R,\S),
	\end{align}
\end{subequations}
with $\hat{\Psivec}$, $\hat{\Upsvec}$ polynomial mappings in $\R$, $\S$. Following the real-valued expressions, and by considering only second-order terms for the sake of brevity, the nonlinear mapping reads~\cite{artDNF2020,AndreaROM}:
\begin{subequations}\label{eq:rv_map}
	\begin{align}
		\Xvec =&\, \sum_{k=1}^{N}\phivec_{k} R_k + 
		\sum_{i=1}^{N}\sum_{j=1}^{N}
		\left( \hat{\avec}_{ij} R_i R_j + \hat{\bvec}_{ij} S_i S_j + \hat{\cvec}_{ij} R_i S_j\right), \\
		\Yvec =&\, \sum_{k=1}^{N}\phivec_{k} S_k + 
		\sum_{k=1}^{N}\sum_{l=1}^{N}
		\left( \hat{\alphavec}_{ij} R_i R_j + \hat{\betavec}_{ij} S_i S_j + \hat{\gammavec}_{ij} R_i S_j\right).
	\end{align}
\end{subequations}
It is worth mentioning that the nonlinear mapping $\hat{\Psivec}$, $\hat{\Upsvec}$ are not independent from one another in the context of vibratory systems so fewer calculations are required to compute them for each polynomial order. Also, summations up to $N$ (the number of dofs) are given in Eq.~\eqref{eq:rv_map}; nevertheless in practice only a small subset of $m$ master {\em normal} coordinates $(R_k,S_k)_{k=1,...,m}$, with $m \ll N$; are selected. Cancelling all slave normal coordinates in Eq.~\eqref{eq:rv_map} leads to rewrit\red{ing} sumations up to $m$.

The unknown vectors of coefficients $\hat{\avec}_{ij}, \hat{\bvec}_{ij}, \hat{\cvec}_{ij}, \hat{\alphavec}_{ij}, \hat{\betavec}_{ij}$ and $\hat{\gammavec}_{ij}$ are derived from the second-order homological equations, see~\cite{touzeLMA,touze03-NNM,AndreaROM} for more details. 
These vectors of coefficients are the equivalent in physical coordinates of the vectors in Eq.~\eqref{eq:nonlinear_change_modal_compact}, and for this reason they are denoted here by the $\hat{}$ superscript. Thanks to the above-mentioned dependence of the velocity mapping $\hat{\Psivec}$ from the displacement one $\hat{\Upsvec}$,  only three of those six vectors are independent in the general case, reducing to only two in the case of a conservative system.
The six unknown vectors can thus be fully computed from the two following vectors of coefficients   $\hat{\Psivec}^{(\mathrm{P})}_{ij}$ and $\hat{\Psivec}^{(\mathrm{N})}_{ij}$, defined as:
\begin{subequations}
	\begin{align}
		& \hat{\Psivec}^{(\mathrm{P})}_{ij} = \left[ (+\omega_{i}+\omega_{j})^{2}\Mvec-\Kvec \right]^{-1}\Gvec\left(\phivec_i,\phivec_j\right), \\
		& \hat{\Psivec}^{(\mathrm{N})}_{ij} =  \left[ (+\omega_{i}-\omega_{j})^{2}\Mvec-\Kvec \right]^{-1}\Gvec\left(\phivec_i,\phivec_j\right).
	\end{align}
	\label{eq:realmappingscalculation}
\end{subequations}
The following solutions are found:
\begin{subequations}
	\begin{align}
		\hat{\avec}_{ij} =&\, \frac{1}{2}
		\left( \hat{\Psivec}^{(\mathrm{P})}_{ij} +\hat{\Psivec}^{(\mathrm{N})}_{ij}    \right), \\
		\hat{\bvec}_{ij} =&\, -\frac{1}{2\omega_{i}\omega_{j}}
		\left( \hat{\Psivec}^{(\mathrm{P})}_{ij} - \hat{\Psivec}^{(\mathrm{N})}_{ij}  \right),\\
		\hat{\cvec}_{ij}  =&\, \zervec,\\
		\hat{\alphavec}_{ij} = &\,\zervec, \\
		\hat{\betavec}_{ij} = &\,\zervec, \\
		\hat{\gammavec}_{ij} = &\,  
		\frac{\omega_{j}+\omega_{i}}{\omega_j}\hat{\Psivec}^{(\mathrm{P})}_{ij} + 
		\frac{\omega_{j}-\omega_{i}}{\omega_{j}}\hat{\Psivec}^{(\mathrm{N})}_{ij} .
	\end{align}\label{eq:realmappings}
\end{subequations}

One can see in particular that the quadratic manifold built from MD as explicited in~Eq.~\eqref{eq:QMNLmap} is a simplification of the more general formula given in Eqs.~\eqref{eq:rv_map} where the velocity-dependence has been properly taken into account. Putting back the velocities as independent coordinates leads to bypass the slow/fast assumption that was limiting the QM approach, thus proposing a uniformly valid and simulation-free method. Interestingly, a conference paper on MDs~\cite{RixenMDNF19} has proposed a first step in this direction, by generalizing the quadratic manifold by taking the velocities into account, finally arriving at formulas equivalent to Eqs.~\eqref{eq:realmappings}. The developments reported in~\cite{artDNF2020,AndreaROM} settles down the full formulation, unfold\red{s} the link with the parametrisation method of invariant manifolds, and open\red{s} the door to higher-order developments by giving the recursive general formulations. \red{From} the computational point of view, the real-valued formalism can also be fully written in a non-intrusive manner, such that the method as derived in~\cite{artDNF2020} can be used from any FE code, provided the code allows a user to script for performing matrix operations online, such that all outputs can be computed without the need to export the full mass and stiffness matrices.

\begin{figure}[h!]
\centering
\includegraphics[width=\linewidth]{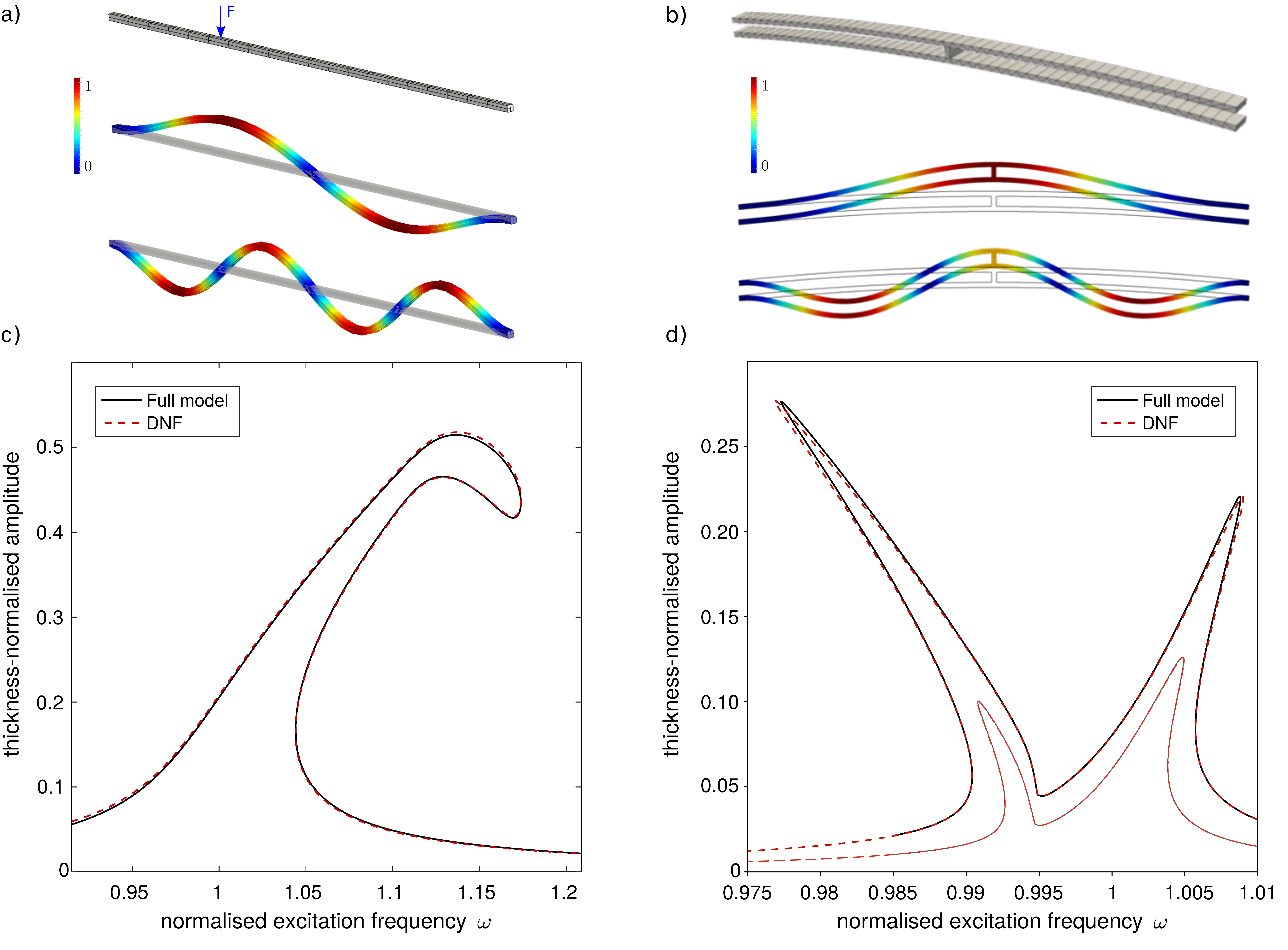}
\caption{Comparisons of full-order solutions and ROMS using DNF for two different structures including forcing and damping. (a) Mesh, second and fourth mode of a clamped clamped beam. (c) FRF in the vicinity of the second eigenfrequency, with the characteristic loop due to the excitation of the 1:3 resonance between modes 2 and 4, already mentioned in Fig.~\ref{fig:resonanceexample}(b). \red{Figure reworked from~\cite{artDNF2020}}. (b) Mesh, first and second mode of a MEMS-like arch. (d) FRF featuring a 1:2 internal resonance between the first two modes, for two different levels of excitation amplitude. Reprinted from~\cite{AndreaROM}. Frequency axis normalized by the eigenfrequency of the driven mode ($\omega_2$ in (b) and $\omega_1$ in (d)), amplitude axis normalized by the thickness (in (d) the thickness of the sub-beam has been selected).}
\label{fig:ziDNFexample}
\end{figure}

%

The DNF method has been already applied to numerous different examples, with or without internal resonance. A fan blade and a clamped-clamped beam \red{are} considered in~\cite{artDNF2020}. The fan blade shows an industrial example with complex geometry while the beam is used to show how the method can handle internal resonances occurring with nonlinear frequencies at large vibration amplitude. Numerous other examples of beams have been considered in~\cite{YichangVib}, including a 1:1 internal resonance between the two polarizations of a beam, an arch with increasing curvature, and a cantilever beam. In~\cite{YichangVib}, the results provided by DNF are also compared to \red{the} QM method with MDs and \red{the} ICE technique. Finally, complex MEMS structures are considered in~\cite{AndreaROM}, including a micromirror with large rotations and complex assemblies of beams and arches featuring 1:2 and 1:3 internal resonance. Fig.~\ref{fig:ziDNFexample} shows two examples of the results obtained. The left column, Figs.~\ref{fig:ziDNFexample}(a,c) reports the case of a clamped-clamped beam with a pointwise excitation with harmonic content in the vicinity of the second eigenfrequency. The frequency response function (FRF) displayed in Fig.~\ref{fig:ziDNFexample}(c) shows the particular loop corresponding to the activation of a strong 1:3 internal resonance with mode 4. Figs.~\ref{fig:ziDNFexample}(b,d) reports the case of an arch MEMS resonator featuring a 1:2 internal resonance between the first two bending modes. The arch is then excited with a modal forcing (having the shape of the first eigenmode), in the vicinity of $\omega_1$, and the typical FRF of systems with 1:2 internal resonance is retrieved by the ROM with only two modes.

In most of the examples reported in~\cite{artDNF2020,YichangVib,AndreaROM}, the second-order DNF has been used, where the nonlinear mapping is truncated at order two while the reduced dynamics is truncated at the third-order, showing that such a simple development already allows \red{one} to obtain excellent results \red{for} numerous test cases up to comfortable vibration amplitudes. Another important remark formulated in~\cite{artDNF2020} relates to the fact that using second-order DNF, the analyst does not need to care about internal resonances higher than second-order. Indeed, since no further development is sought, all possible internal resonances from the third-order are not treated and can thus be excited. Hence the loss in accuracy on the geometry of the manifold (limited to order two) is compensated with a  more easy treatment of internal resonances in the reduced dynamics. The reported gains in computational times are in general impressive, with a burden of the order of 1-2 days for the computation of the full-order solution to obtain a complete frequency response curve, as compared to 1-2 minutes for the ROMs. Scalability of the methods up to millions of dofs is also illustrated in~\cite{AndreaROM}, illustrating that up to 3 millions of dofs the time needed to build the ROM is less than one hour. As a general note, one can observe that memory requirements and computational burden of the method are not important, generalizing the QM method with MD and proposing ROMS at the same light computational cost. Note that the second-order DNF will be released as a command in the version of the open FE code \texttt{code\_aster}~\cite{ASTER} from version 15.4 (June 2021).

\subsection{Direct computation of SSM}\label{sec:directSSM}

Direct computation of ROMS using parametrisation method of invariant manifolds has been proposed in~\cite{VERASZTO,JAIN2021How}, following the previous developments led on the definition of SSMs, reported in Section~\ref{subsec:SSM}. In the first contribution~\cite{VERASZTO}, the direct formulation is given by selecting a single master mode, and up to the third-order. The general formula are written for an initial system with mixed coordinates, since the master mode is assumed to be expressed in the modal basis, while all the slave modes are given in the physical space. When damping is not considered, reduction to \red{a} LSM is given. Interestingly, the authors use the graph style for parametrising the LSM, probably due to the initial choice for expressing the dynamics. The reduction formula they give are thus very close to the general formula given in~\cite{PesheckBoivin}, but fully rewritten from the physical space for the slave coordinates only, thus making appear the mass and stiffness matrices (truncated to slave coordinates), as well as the quadratic nonlinearity. The reduced dynamics is equivalent to Eq.~\eqref{eq:romIM1mdev}. When forcing and damping is added, the formulation is extended to reduction onto SSM as derived in their previous papers~\cite{Haller2016,PONSIOEN2018}. Although generalizing important previous developments, the starting point selected by the authors in this contribution is not generically used and is a priori not standard for FE models, and a first linear transform needs to be sorted out for an easy use of the given formulas. In~\cite{VERASZTO}, the method is applied to a chain of oscillator\red{s} and to a Timoshenko beam, but not directly to large FE models.

The direct computation has been fully tackled in~\cite{JAIN2021How}, starting from the equations of motion in physical space, such that direct applications to FE models are proposed. General formulations are given without specific a priori restrictions neither on the eigensolutions (self-adjointness of linear operators not mandatorily required), nor on the form of the nonlinearities. Constant, periodic and quasi-periodic forcing terms are included in the framework. The equations of motion are written as first-order system in order to fit the guidelines of dynamical system theory, so that the starting point is a system written as
\begin{subequations}
\begin{align}
\Bvec \dot{\zvec} & = \Avec \zvec + \Fvec (\zvec) + \varepsilon  \Fvec^{\mathrm{ext}}  (\zvec,\phivec), \\
\dot{\phivec} &= \Omegavec,
\end{align}
\end{subequations}
where $\zvec$ is the $2N$-dimensional vector including displacements and velocities, $\Avec$ and $\Bvec$ are the first-order matrices composed of the usual stiffness, mass and damping matrices of the $N$-dof mechanical system, $\Fvec$ gathers the nonlinear terms, $\varepsilon$ is a small parameter and $\Fvec^{\mathrm{ext}}$ is the external forcing composed of $K$ forcing radian frequencies grouped in the vector $\Omegavec = (\omega_1, ..., \omega_K)$. 

The invariance equation \eqref{eq:invariance} is then rewritten to fit this starting point with $\Avec$ and $\Bvec$ matrices, thus reading~\cite{JAIN2021How}:
\begin{equation}
\Bvec (D \Wvec) \Rvec = \Avec \Wvec + \Fvec \circ \Wvec,
\end{equation}
where $\Wvec$ stands for the nonlinear mapping as in Eq.~\eqref{eq:NLmap}, while $\Rvec$ represents the reduced dynamics on the invariant manifold. The computational procedure then proceeds following canonical rules with asymptotic polynomial expansions. Identifications of co-homological equations are a bit more involved since their writing in physical space is less obvious, nevertheless with correct projections onto kernels, closed expressions are attainable. Both graph style and normal form style are highlighted, and recursive expressions to deal with higher-order expansions are given, thus offering an automated version of the whole procedure. All the developments are coded in an open software, SSMtools 2.0, which is interfaced with a FE solver. In~\cite{JAIN2021How}, A FE model of a von K{\`a}rm{\`a}n beam having 30 dofs is used as first example, then a shallow parabolic panel with 1320 dofs is selected and the softening frequency response curve is computed. Finally, a FE model of an aircraft wing with 133 920 dofs is investigated, and the hardening behaviour is reported. Even though the framework is given in its most general formulation, examples reported in~\cite{JAIN2021How} are restricted to single master mode dynamics without internal resonance, and to geometric nonlinearity, such that extensions to handle different kind of nonlinear terms still remains to be shown.

\section{Open problems and future directions}\label{sec:future}

This section aims at underlining open problems related to the use of invariant manifold theory for reduced-order modeling of the nonlinear dynamics of structures, in order to point to possible further developments and future research directions.\\

A first open question is related to the folding of invariant manifolds in phase space and the consequences on the dynamics, the parametrisation and the performance of the ROM. The theoretical developments on the parametrisation method underline a main difference between graph style and normal form style. 
While normal form style  is theoretically able to overcome the potential foldings of the IM,  such a distinction has never been clearly emphasised and illustrated in nonlinear vibration theory. Indeed, most of the foldings reported in the literature appear through loops in the backbone curve (see Figs~\ref{fig:resonanceexample}(b), \ref{fig:manifolds}) and are related to the appearance of an internal resonance between the nonlinear frequencies of the system. However the folding seems to be apparent and only due to the projection of (at least) 4-D dimensional manifolds into 3-d representation. Finding out a clear folding without internal resonance would underline the difference in the parametrisation styles.

A second open and important problem is related to the {\em a priori} estimate of the quality of the ROM developed using IM approach. For linear reduction methods such as POD for example, a priori estimates are easy to find since the magnitude of the singular values gives the amount of energy captured by the ROM, which can be directly linked to the accuracy of the reduction. For nonlinear methods, the picture is less clear since using curved manifolds to represent the dynamics, one cannot rely on linear ideas for error estimate. The problem is addressed in~\cite{Haro}, and a priori and a posteriori estimates are proposed based on errors on the invariance equation or errors on the orbits. An upper bound for validity limits of normal transform is also proposed in~\cite{LamarqueUP}. However all these preliminary ideas did not \red{yet} translate to a simple tool, and are also related to the ease of developing high-order approximations. In this realm, the automated computations of higher orders proposed in~\cite{PONSIOEN2018,JAIN2021How} \red{gives} a different point of view first developed by these authors. Instead of using  added vectors to reach convergence of the ROM, the idea is to state that the dimension of the manifold is given by the dynamics at hand, such that convergence is simply reached by adding higher-orders up to a converged backbone (or frequency response curve). Whereas this can work in many cases, the problem of internal resonance between the nonlinear frequencies of the system (thus appearing at higher vibration amplitudes), still remains an obstacle for blind application of this idea. Indeed, the parametrisation method can only check resonance relationship\red{s} between the linear eigenfrequencies of the problem, and can warn that additional master modes are needed based on this inspection. 
Moving to higher amplitudes, resonance can be fulfilled between nonlinear frequencies. 
Developing new methodologies to give accurate and uniform a priori estimates is thus still a question that deserves further investigations.

Future research directions along these lines should enlarge the scope of the methodologies given in sections~\ref{sec:DNF}-\ref{sec:directSSM}. Indeed, strong results from dynamical system theory ensure that long-term behaviour of the dynamical solutions are contained within attracting center sets~\cite{Carr,Temam1989}. Hence the reduction methods based on these theorems are the most accurate way of deriving efficient ROMs, the main question being more their accurate computation. In this way, recent progress reported in \ref{sec:DNF}-\ref{sec:directSSM} shows that effective methods can now be used for a direct computation from the physical space. Enlarging these recent results appears thus as a logical path to the improvement of ROM computation.

In this realm, the following directions should be investigated soon as direct applications of the general method. First of all, applications to different physical problems, including different types of nonlinear forces, should be investigated, as for example nonlinear damping laws~\cite{AmabNLdamp1,AmabNLdamp2,colin2020}, coupling with other physical forces such as piezoelectric couplings~\cite{LazarusThomas2012,GivoisPiezo,shami2021}, piezoelectric material nonlinearities \cite{vonwagnerDSTA,leadenham2015,frangi2020}, non-local models for nanostructures \cite{ribeiro17-JCNLD,RIBEIRO2019}, often used in energy-harvesting problems, electrostatic forces in MEMS dynamics~\cite{Younisbook}, centrifugal and Coriolis effects in rotating systems~\cite{DELHEZ2021,tacher2021}  with applications to blades~\cite{Turbomachine,PESHECK2001Blade,PesheckASME2002,thomas16-ND}, large strain elastic nonlinear constitutive laws \cite{MORIN201866}, fluid-structure interaction~\cite{MATTHIES2002,Noels2006} and coupling with nonlinear aeroelastic forces~\cite{Dimitriadis}; or thermal effects~\cite{PerezThermal,JAIN2020T}, to cite a few of the most obvious directions where the general reduction strategy could be easily extended. Extensions to structures with symmetries, in order to get more quantitative informations and highlight the link with mode localization could be also used with such tools~\cite{Wei1,Wei2,Vakakis93,georgiades09}.

Another interesting research direction would be to enlarge the scope of invariant-based ROMs to tackle more complex dynamics involving a larger number of master modes. While most of the reported applications uses 1-3 master modes, jumping to 10-20 master modes and investigate the transition to more complex dynamics and the limits of the methods with regard to chaotic vibrations and wave turbulence raises a number of open questions. Also in this realm, a link with more involved mathematical analysis related to the existence of inertial manifold and global attractors in structural dynamics~\cite{Chuechov2002,CHUESHOV2004,Chuechov2008,BalaWebster}, should be of interest.

A closer connection and better understanding of other theoretical efforts could also help in unifying the concepts and calculation methods. In this direction, the link with Koopman operator and Koopman modes is an interesting topic that has been first investigated in~\cite{CIRILLO2016}, based on the general results derived for example in~\cite{Mezic05,LAN2013,Mezic_annual,Mauroy16}. In short and following~\cite{CIRILLO2016}, the Koopman operator replaces a finite dimensional nonlinear dynamical system by an infinite dimensional linear system. The eigensolutions of the Koopman operator are infinite and contains all the eigensolutions of the original dynamical system, say eigenvectors $\phivec_1,...,\phivec_N$ with associated eigenvalues $\lambda_1, ..., \lambda_N$, as well as all possible combinations of the form $\phivec_1^{k_1} \times \phivec_2^{k_2} \times ... \times  \phivec_N^{k_N}$ for $(k_1,...k_N)$ integers, with eigenvalues $k_1\lambda_1 + ... + k_N \lambda_N$. By making appear the nonlinear resonance relationships as Koopman modes allows redefining the problem with a large (infinite) dimension instead of seeing them appearing order by order. Investigating further the computational properties of this equivalence might be helpful for deriving other numerical methods.

For nonlinear systems, identification methods is also a very active field of research in order to extract important model characteristics from experimental data~\cite{NOEL2017}. In this field, invariant manifold theory is already used through NNMs and identification of backbone curves with  generalization of {\em e.g.} phase separation or force appropriation techniques, see~\cite{PAI2011,PEETERS2011,NOEL2016} and references therein. The use of normal form to select {\em ex nihilo} a suitable nonlinear model before experimental identification, especially in the case of internal resonances, is also a powerful tool \cite{thomas03-JSV,thomas07-ND,monteil15-AA,jossic18}, which has been recently coupled to phase locked loop experimental continuation method \cite{denis18-MSSP,Givois11}. The use of SSM has also been reported for model identification in~\cite{Szalai2017}, and a recent contribution propose\red{s} to enlarge the scope by using {\em spectral foliations} in model identification in order to better take into account transients and orthogonal directions to the SSMs~\cite{Szalai2020}.

A long term research direction is the application of invariant manifold theory to nonsmooth problems occurring in mechanics, mainly through contact and friction problems. The main limitation appears on the smoothness of the dynamical system, since all the theorems used in~\cite{Haro} assumes a sufficiently smooth map, and the smoothness order is directly linked to the order of the SSM, see {\em e.g.}~\cite{Haller2016}. Research in these directions considers Filippov systems, see {\em e.g.}~\cite{NonSmoothBrogliato1999,NumericalMethodsBrogliato2009,LEINE2006,Leine2013}, while penalisation and regularisation methods are often used and could be adapted. Recent attempts to settle down a nonsmooth modal analysis and extends the invariant manifolds to impacting systems are investigated in~\cite{thorin2017,legrand2017,yoong2018,thorin2018}, showing how nonsmooth modes of vibration can be defined. 

%
%
%
%
\section{Conclusion}

In this contribution, a review \red{of} the nonlinear methods for deriving accurate and efficient ROMs for geometrically nonlinear structures, is given. Nonlinear methods differ from linear methods by defining a nonlinear mapping between the initial and reduced coordinates. In a phase space perspective, this leads to projection onto a curved manifold instead of using orthogonal vectors to decompose the dynamics. Though proposing a more involved calculation at first sight, reduction to nonlinear manifolds \red{is} then expected to produce more accurate results with fewer master coordinates since embedding the geometric complexity into the nonlinear mapping.

A special emphasis has been put on methods based on invariant manifold theory, and the  strong results provided by  dynamical system theory in order to derive efficient and accurate predictive ROMs. Indeed, invariant manifold methods differ from other\red{s} by the fact that \red{the related} theorems ensure that the long-term dynamics of the mechanical systems are enclosed in the vicinity of these subsets. Consequently full-order solution\red{s} displaying low-order dynamics exactly relies on these manifolds, such that computing their characteristics is the key to derive the most accurate ROMs. This point of view is different from ad-hoc methods that can be compared on their predictive accuracy. Here the problem is not to find the correct representative set which is known theoretically, but to compute it efficiently. As underlined, the point of view is geometrical in nature. The curvatures of the invariant manifolds in phase space have a strong meaning and relate to the non-resonant couplings. Capturing them accurately is  computationally  more involved  as compared to linear reduction method, but offers better performance and stronger reduction.

The presentation followed historical developments, and a focus has been set on the derivation of the parametrisation method of invariant manifolds that offers a unified and comprehensive point of view, used for model-order reduction using ad-hoc terminology (LSM for conservative systems and SSM for damped systems). Then for the importance of applications, the special case of FE structures has been specially developed.

In the context of FE structures, new questions arise due to the fact that existing powerful FE codes might be used non-intrusively for deriving ROM, an appealing feature offering great versatility. Also, specific developments have been led in the field of computational mechanics, with the development of the STEP, implicit condensation and modal derivatives. All these methods are reviewed and systematically compared to invariant manifolds, showing that they suffer from a lack of generality and need extra assumptions such as a slow/fast separation to be used blindly. However, they all have important benefits in the ease-of-use, rapidity and efficiency of the computation, non-intrusiveness, and give excellent result when used with the correct assumptions fulfilled.

The paper concludes with the latest developments in the field showing how one can use invariant manifold based ROMs, directly from a FE mesh, and possibly in a non-intrusive manner. The proposed methods are in general simulation-free, and can be used with a computational cost that is of the same order of {\em e.g.} modal derivative based techniques. As a conclusion, we advocate for a more general use of nonlinear techniques for efficient ROM computation for geometrically nonlinear structures. This point of view makes a direct link between the large dimension of initial problems meshed with FE and the generally small dimensional subsets where the important dynamics is contained and allows one to compute efficient ROMS that can be used for a lot of different purposes including analysis and design. Open problems and future directions are briefly listed at the end of the paper, underlining that many interesting developments can be conducted to generalize the methods to a large number of cases.

\section*{Acknowledgments}
The authors are especially thankful to Alex Haro for detailed discussions we had in January 2021. \red{Walter Lacarbonara is also warmly thanked for inviting the first author to write a review paper on reduction methods.} Andrea Opreni and Attilio Frangi are thanked for their major contributions to the developments of efficient writing and coding of the DNF for FE systems and the interesting applications to MEMS as well as first-draft readings. Lo{\"i}c Salles is thanked for all the discussions and collaborations on the subject and for launching the subject again in January 2019. Yichang Shen is thanked for his involvement in the project. \red{Steve Shaw and David Wagg read first versions of the paper and brought out numerous interesting comments to the authors that helped improving the presentation, they are sincerely thanked for this precious help. Claude Lamarque, G\'erard Iooss and Paul Manneville are thanked for the non-countable discussions we had on normal forms in the last years.}

\section*{Funding}
The work received no additional funding.

\section*{Conflict of interest}
The authors declare that they have no conflict of interest.
\section*{Data availability statement}
The data that support the findings of this study are available from the corresponding author, upon reasonable request.
\bibliographystyle{spmpsci}
\bibliography{biblioROM}

\appendix

\section{Symmetry of the quadratic and cubic tensors}
\label{sec:ghsym}

This appendix is devoted to the demonstration of the symmetry properties of the nonlinear tensors of coefficients $\G$, $\H$, $\g$ and $\h$, that appear in the equations of motion written in terms of FE coordinates, Eq.~\ref{eq:eom_phys}), or with  modal coordinates, Eq.~(\ref{eq:eom_modal}). 

The FE coordinates are first considered. Let us denote by $\kvec(\Xvec)=\Kvec\Xvec+\fvec_\text{nl}(\Xvec)$ the internal force vector. Using indicial notations and Einstein summation convention, it can be written explicitely, for $i,j,l,s=1,\ldots N$:
\begin{equation}
\label{eq:ks}
k_s=K_{si}X_i+G^s_{ij}X_iX_j+H^s_{ijl}X_iX_jX_l,
\end{equation}
where $k_s$ is the $s$-th component of the internal force vector $\kvec$, $K_{si}$ are the components of the stiffness matrix $\Kvec$; while $G^s_{ij}$, $H^s_{ijl}$ are the quadratic and cubic coefficients defined in Eq.~(\ref{eq:tensor_product_phys}). In a 3D finite element context, the physical displacement vector $\uvec(\yvec)$ is interpolated on a family of shape functions $\Nvec_i$, such that $u_\alpha(\yvec)=N_{\alpha i}(\yvec)X_i(t)$, $\alpha=1,2,3$. Using Voigt notations, the Green-Lagrange strain tensor , Eq.~\eqref{eq:3DMMC}, and its variation can be written, for $\alpha=1,\ldots 6$ (see \cite{givois21-CS}):
\begin{equation}
E_\alpha=B^{(1)}_{\alpha i}X_i+\frac{1}{2}B^{(2)}_{\alpha ij}X_iX_j,\quad \delta E_\alpha=B^{(1)}_{\alpha i}(\yvec)\delta X_i+B^{(2)}_{\alpha ij}X_i\delta X_j,
\end{equation}
with $B^{\alpha(1)}_{i}$ of size $6\times N$ and $B^{\alpha(2)}_{ij}$ of size $6\times N\times N$ are two discretized gradients operators, defined by:
\begin{equation}
\label{eq:Bsij}
B^{(1)}_{\alpha i}=\begin{bmatrix}
N_{1i,1} \\
N_{2i,2} \\
N_{3i,3} \\
N_{2i,3}+N_{3i,2} \\
N_{1i,3}+N_{3i,1} \\
N_{1i,2}+N_{2i,1}
\end{bmatrix}
\qquad
B^{(2)}_{\alpha ij}=\begin{bmatrix}
N_{li,1 }N_{lj,1}\\
N_{li,2 }N_{lj,2}\\
N_{li,3 }N_{lj,3}\\
N_{li,2 }N_{lj,3}+N_{li,3 }N_{lj,2}\\
N_{li,3 }N_{lj,1}+N_{li,1 }N_{lj,3}\\
N_{li,2 }N_{lj,1}+N_{li,1 }N_{lj,2}
\end{bmatrix}
\end{equation}
where each row of the matrices corresponds to the corresponding value of index $\alpha=1,\ldots 6$ and $N_{\alpha i,\beta}=\partial N_{\alpha i}/\partial y_\beta$ is a space derivative of the shape functions.


Then, the virtual work of the internal forces can be written:
\begin{align}
\delta\mathcal{W}_\text{int} & =\int_\Omega C_{\alpha\beta}E_\alpha\delta E_\beta\,\dd\Omega = k_s\delta X_s \\
 & =\Bigg[\underbrace{\int_\Omega C_{\alpha\beta}B_{\alpha i}^{(1)}B_{\beta s}^{(1)}\,\dd\Omega}_{K_{si}} \,X_i + \underbrace{\int_\Omega C_{\alpha\beta}\left(B_{\alpha i}^{(1)}B_{\beta js}^{(2)}+\frac{1}{2}B_{\alpha s}^{(1)}B_{\beta ij}^{(2)}\right)\,\dd\Omega}_{G^s_{ij}} \,X_i X_j + \underbrace{\int_\Omega \frac{1}{2}C_{\alpha\beta}B_{\alpha ij}^{(2)}B_{\beta ls}^{(2)}\,\dd\Omega}_{H^s_{ijl}} \,X_i X_jX_l \Bigg]\delta X_s 
\end{align}
where $C_{\alpha\beta}$ is the elasticity tensor in Voigt notations. The above equation defines the tensor components $K_{si}$, $G^s_{ij}$ and $H^s_{ijl}$ of the internal forces.

The elasticity tensor is symmetric ($C_{\alpha\beta}=C_{\beta\alpha}$) and Eq.~(\ref{eq:Bsij}) shows that $B^{(2)}_{\alpha ij}=B^{(2)}_{\alpha ji}$: on can invert the two latin subscripts. This leads to a symmetric stiffness matrix and allows {\em any change of the order of the indices for the cubic coefficients} ($4!=24$ possibilities if $i\neq j\neq l\neq s$):
\begin{equation}
\label{eq:Hsijlsym}
H^s_{ijl} = H^l_{sij} = H^j_{lsi} = H^i_{jls} = H^s_{jil}=\ldots
\end{equation}

For the quadratic coefficients, no symmetry appears. However, the coefficients  $G^s_{ij}$ and $G^s_{ji}$ refer to the same monomial $X_iX_j$ in Eq.~(\ref{eq:ks}). Due to the commutativity property of the usual product, one then understands that only the summation of these two quadratic coefficient $G^s_{ij}$ and $G^s_{ji}$ matters.  This leads several authors (see e.g. \cite{muravyov,touze03-NNM,givois2019}) to adopt a so-called upper triangular form for those tensors for which Eq.~(\ref{eq:ks}) is rewritten as
\begin{equation}
\label{eq:ksutf}
k_s = \sum_{i=1}^N K_{si}X_i + \sum_{i=1}^N \sum_{j=i}^N \hat{G}^s_{ij} X_i X_j + \sum_{i=1}^N \sum_{j=i}^N \sum_{l=j}^N \hat{H}^s_{ijl} X_i X_j X_l.
\end{equation}
With this selection, an unequivocal representation of the monomials is given, and the coefficients are attributed such that only those with increasing indices ($l\geq j \geq i$) are non-zero, while the other ones ($l\leq j \leq i$) are set to zero. For the quadratic coefficients, this leads to, for all $s$:
\begin{gather}
\hat{G}^s_{ii}=G^s_{ii}, \quad \forall i, \\
\hat{G}^s_{ij}=G^s_{ij}+G^s_{ji}=\int_\Omega C_{\alpha\beta}\left(B_{\alpha i}^{(1)}B_{\beta js}^{(2)}+B_{\alpha j}^{(1)}B_{\beta is}^{(2)}+B_{\alpha s}^{(1)}B_{\beta ij}^{(2)}\right)\,\dd\Omega,\quad  \forall j>i, \label{eq:Ghatdef}\\
\hat{G}^s_{ji}=0,\quad \forall j>i.
\end{gather}
In the above Eq.~(\ref{eq:Ghatdef}), any change of the order of the indices $s,i,j$ is allowed, which leads to the following properties:
\begin{equation}
\label{eq:symG}
\hat{G}^i_{ij}=2\hat{G}^j_{ii},\quad \hat{G}^j_{ij}=2\hat{G}^i_{jj},\quad \hat{G}^l_{ij}=\hat{G}^j_{il},\quad \hat{G}^l_{ij}=\hat{G}^i_{lj},\quad \forall i<j<l
\end{equation}
Analog expressions are obtained for the cubic coefficients $\hat{H}^s_{ijl}$ (see \cite{muravyov,givois2019}). 

In the present article, we found convenient to use the standard form~(\ref{eq:ks}) instead of the upper triangular form~(\ref{eq:ksutf}) for all our demonstrations and we enforced the symmetry on the quadratic coefficients by redefining them as the symmetric part of their upper triangular counterparts: $G^s_{ij}=G^s_{ji}=\hat{G}^s_{ij}/2$, $\forall j\neq i$. This leads to allows {\em any change of the order of the indices also for the quadratic coefficients}:
\begin{equation}
\label{eq:Gsijsym}
G^s_{ij}=G^s_{ji}=G^i_{js}=G^i_{sj}=G^j_{is}=G^j_{si}.
\end{equation}

All the above reasoning {\it equally applies} to the modal coefficients $g^s_{ij}$  and $h^s_{ijl}$ of Eq.~(\ref{eq:eom_modal}), since, according to Eq.~(\ref{eq:g_ij}),(\ref{eq:h_ijk}), rewritten in full indicial form with Einstein notation, one has:
\begin{equation}
g^p_{mn}=G^s_{ij}\phi_{sp}\phi_{im}\phi_{jn},\quad h^p_{mnq}=H^s_{ijl}\phi_{sp}\phi_{im}\phi_{jn}\phi_{lq},\quad\forall p,m,n,q=1,\ldots N,
\end{equation}
where $\phi_{jn}$ refers to the $j$-th component of the $n$-th eigenvector. Consequently,  the same symmetry as in Eqs.~(\ref{eq:Gsijsym}) and (\ref{eq:Hsijlsym}) applies for the modal coefficients:
\begin{gather}
g^s_{ij}=g^s_{ji}=g^i_{js}=g^i_{sj}=g^j_{is}=g^j_{si}, \label{eq:zesymg}\\
h^s_{ijl} = h^l_{sij} = h^j_{lsi} = h^i_{jls} = h^s_{jil}=\ldots \label{eq:zesymh}
\end{gather}
Again, {\it any permutation of the indices $s,i,j,l$ is possible.}

All the above symmetries are also a consequence of the existence of a potential energy \cite{muravyov,givois2019}:
\begin{equation}
\mathcal{V}=\frac{1}{2}\int_\Omega C_{\alpha\beta}E_\alpha(\Xvec) E_\beta(\Xvec)\,\dd\Omega,
\end{equation}
which is a quartic polynomial in $\Xvec$. Then, the $i$-th component of the internal force vector can be directly derived from the elastic energy as:
\begin{equation}
k_i(\Xvec)= \dfrac{\partial \mathcal{V}}{\partial X_i}.
\end{equation}
Using Schwarz's theorem, one can write:
\begin{equation}
\dfrac{\partial^2 \mathcal{V}}{\partial x_i\partial x_j} = \dfrac{\partial^2 \mathcal{V}}{\partial x_j\partial x_i}\quad\Rightarrow\quad \dfrac{\partial k_i}{\partial x_j} = \dfrac{\partial k_j}{\partial x_i}. 
\end{equation}
Then, identifying the coefficients of identical monomials in the above last equations with the upper triangular form for $k_i$ leads to prove all the symmetry properties of Eq.~(\ref{eq:symG}) and their analogs for the cubic coefficients.

\section{Classification of nonlinear terms}\label{app:class}

This section is devoted to give more details of the terminology used throughout the text to classify the different nonlinear coupling terms and monomials appearing in the dynamics. It is based on previous developments reported {\em e.g.} in~\cite{touze03-NNM,TouzeCISM}, and all the textbooks dealing with normal form theory, where the reader can find more details. Systems with geometric nonlinearity are essentially driven by a large assembly of nonlinearly coupled oscillators, thus generating a very large number of coupling terms (of the order of $N^4$ terms). However, all the terms does not play the same role and it is important to identify the contributions of each monomials. For this discussion on the terminology, we use the equations of motion written in modal space, Eq.~\eqref{eq:eom_modal}, and more specifically, 
$\forall \; p \, \in \, [1,...,N]$:
\begin{equation}
\ddot{x}_p + \omega_p^2 x_p +  \sum_{i=1}^{N} \sum_{j \geq i }^{N} g_{ij}^p x_i x_j + \sum_{i=1}^{N} \sum_{j\geq i}^{N} \sum_{k \geq j}^{N}  h_{ijk}^p x_ix_jx_k = 0.
\label{eq:EDOmodalprojapp}
\end{equation}
For the discussion herein, assume that $m$ is the master mode such that most of the energy is contained within $x_m$. {\em Invariant-breaking terms} have already been commented in the main text, they are the monomials on $p$-th oscillator equation with the form $g^p_{mm} x_m^2$ and $h^p_{mmm} x_m^3$. As soon as $x_m \neq 0$, then all $x_p$ where this terms exist will also be excited and will thus have a non-zero amplitude. These terms break the invariance of the linear eigensubspaces and can be directly tracked in the equations defining the geometry of the invariant manifolds, as underlined in sections~\ref{sec:IMSP} and \ref{sec:NFmodal}.

In order to go ahead in this classification, the link with internal resonance must be properly understood. Let us start by underlining that any nonlinear term can be interpreted as a forcing term for the corresponding oscillator equation. Continuing with the same example, the term $g^p_{mm} x_m^2$ is a forcing term on oscillator $p$. Interestingly, at the lowest order of approximation, $x_m \sim \expo^{\pm i \omega_m t}$, such that $x_m^2$ will create a forcing with frequency components $2\omega_m$ and $0$. From this we can conclude that for oscillator $p$, if $\omega_p \simeq 2\omega_m$ then the forcing term $g^p_{mm} x_m^2$ will be a resonant forcing term, exciting component $p$ in the vicinity of its eigenfrequency thus creating large amplitude response. One then call the monomial $g^p_{mm} x_m^2$ a {\em resonant monomial}. On the other hand, as long as $\omega_p \neq 2\omega_m$, then the forcing term is non-resonant, and the corresponding monomial is non-resonant. 

This simple example generalizes thanks to normal form theory. The resonance relationship  are then linked with internal resonance between the eigenfrequencies of the system, all of them being connected to a specific order of nonlinearity, conducting to so-called second-order internal resonance and third-order internal resonance. Due to the fact that for linear conservative system, the eigenspectrum is purely imaginary $\{ \pm i \omega_r \}$, some third-order relationships are always fulfilled, all those of the form:
\begin{equation}
\forall {r,p}=[1,...,N]: \quad +i\omega_r = +i\omega_p -i \omega_p + i\omega_r.
\label{eq:trivialres}
\end{equation}
These resonances are called {\em trivial resonance} and their associated monomials are called {\em trivially resonant monomials}. The main consequence in terms of normal form is that all these monomials cannot be cancelled from the normal form of the system. The normal form is not linear but stay nonlinear with only these trivially resonant monomials in case of no other internal resonances between the eigenfrequencies, following the general results from Poincar{\'e} and Poincar{\'e}-Dulac's theorems.

From these developments we can derive the following classification:
\begin{itemize}
\item {\em trivially resonant monomials}: for the $m$-th oscillator, all the terms $x_m^3$, $x_m x_p^2$, $\forall p=1,...N$, are trivially resonant monomials, corresponding to the trivial resonance relationships \eqref{eq:trivialres}. They cannot be cancelled from the normal form and stay in the resulting equations as ascertained in Eq.~\eqref{eq:ROMNF}. Note in particular than none of these are invariant-breaking, recovering the fact that the dynamics is well expressed in an invariant-based span of the phase space for Eq.~\eqref{eq:ROMNF}, which is not the case for the equations of motion in modal space, Eq.~\eqref{eq:EDOmodalprojapp}.
\item {\em resonant monomial}: a resonant monomial is the nonlinear term connected to an internal resonance between the eigenfrequencies of the system. For any internal resonance, a few monomials exist which can be tracked from the simple interpretation of nonlinear term as forcing.
\item {\em non-resonant monomial}: a nonlinear term that is not connected to an internal resonance.
\end{itemize}
When an internal resonance exist, the corresponding resonant monomials create what is generally called a {\em strong}, resonant coupling, and they convey the energy exchange between the coupled oscillators, leading to more complex form of the dynamics with bifurcation in a larger phase space. Otherwise, the coupling is termed {\em weak} or non-resonant.

\section{Parametrisation of invariant manifold}\label{app:IMparam}

In this appendix,  more details on the parametrisation method for the computation of invariant manifold of vector fields in the vicinity of a fixed point, are given. The presentation follows strictly the one given by Haro {\em et al.}, using their notations and developments, so that all the credit of the presentation reported here is given to~\cite{Haro}. Here a simple summary (with slight simplifications in the presentation) is provided and the interested reader is referred to~\cite{Haro} for more details. 

The unknowns $\Wvec$ and $\fvec$ are searched as polynomial expansions of order $k$. They are computed order by order, under the form
\begin{subequations}\label{eq:expansion}
\begin{align}
\Wvec (\svec) &= \zvec_{\star} + \Lvec \svec + \sum_{k \geq 2} \Wvec_k (\svec), \label{eq:expansionW}\\
\fvec (\svec) &= \Lambdavec_L \svec + \sum_{k \geq 2} \fvec_k (\svec),\label{eq:expansionf}
\end{align}
\end{subequations}
where $\Lvec$ is the restriction of the matrix of eigenvectors to the master modes of interest, and $\Lambdavec_L$  the linear diagonalized part of the dynamics restricted to the master modes contained in $\Lvec$. $\Wvec_k (\svec)$ represents an $n$-dimensional vector of $d$-variate homogeneous polynomials of degree $k$, starting at second order (quadratic terms), while $\fvec_k$ is $n$-dimensional.

By replacing \eqref{eq:expansion} into the invariance equation \eqref{eq:invariance} and identifying terms with the same power $k$, on obtains the so-called order-$k$ homological equation as
\begin{equation}
\DD\Fvec(z_{\star}) \Wvec_k (\svec) - \Lvec \fvec_k (\svec) - \DD\Wvec_k (\svec) \Lambdavec_L \svec = -\Evec_k (\svec),
\label{eq:cohomologik}
\end{equation}
where the order-$k$ error term $\Evec_k (\svec)$ has been introduced as:
\begin{equation}
\Evec_k (\svec) = \left[  \Fvec (\Wvec_{<k} (\svec))   \right]_k - \left[ \DD\Wvec_{<k} (\svec) \fvec_{<k} (\svec)  \right]_k,
\end{equation}
and where the shortcut notation $\left[\; .\; \right]_k$ refers to the selection of $k$-th order terms only, while $\Wvec_{<k}$ refers to all orders strictly smaller than $k$.

More insight can be given by projecting onto the modal coordinates, then allowing separating contributions due to master and slave coordinates. Denoting as $\Pvec$ the vector of eigenfunctions, one can introduce:
\begin{equation}
\xivec_k (\svec) = \Pvec^{-1} \Wvec_k (\svec),
\end{equation}
the coefficients of the nonlinear mapping expressed in the modal basis, as well as 
\begin{equation}
\etavec_k (\svec) = \Pvec^{-1} \Evec_k (\svec),
\end{equation}
the expression of the error-$k$ vector in the modal basis. 

The vector $\xivec_k$ can be split as follows:
\begin{equation}
\xivec_k (\svec) = \left[\begin{array}{c}
\xivec_k^L (\svec) \\
\xivec_k^N (\svec)
\end{array}
\right],
\end{equation}
where the first $d$ lines $\xivec_k^L$ is the tangent part, related to the original linear matrix $\Lvec$ containing the master mode coordinates, and the last $n-d$ lines, $\xivec_k^N$, refers to the normal part (slave coordinates). The normal part of Eq.~\eqref{eq:cohomologik} is now called the {\em normal co-homological equation}, and reads
\begin{equation}
\Lambdavec_N \xivec_k^N (\svec) - D \xivec_k^N (\svec) \Lambdavec_L \svec = \etavec_k^N (\svec),
\label{eq:cohomnormal}
\end{equation}
where the second member vector $\etavec$ has been split following the same notation. One must first solve this equation as only depending on one unknown $\xivec_k^N$. The remaining part is called the {\em tangent co-homological equation}, and reads
\begin{equation}
\Lambdavec_L \xivec_k^L (\svec) - D \xivec_k^L (\svec) \Lambdavec_L \svec  - \fvec_k (\svec) = \etavec_k^L (\svec),
\label{eq:cohomtangent}
\end{equation}

In order to fully analyze the co-homological equations and express their solutions, let us first denote as $\xivec_k  (\svec) = [\xi_k^1  (\svec), ..., \xi_k^n  (\svec)]^t$ the $n$ components of the vector $\xivec_k  (\svec)$, which all are homogeneous polynomial of degree~$k$. The same notation is used for the second vector of unknowns, $\fvec_k  (\svec) = [f_k^1 (\svec), ..., f_k^d (\svec)]^t$, which is $d$-dimensional and also composed of homogeneous polynomial of degree $k$. The known vector $\etavec_k  (\svec) = [\eta_k^1  (\svec), ..., \eta_k^n  (\svec)]^t$ is developed as well following the same indicial notation. The normal part of the cohomological equation \eqref{eq:cohomnormal} can now simply be rewritten component by component, for $i=d+1, ..., n$:
\begin{equation}
\lambda_i \xi_k^i  (\svec) - D\xi_k^i  (\svec) \Lambdavec_L \svec = \eta_k^i (\svec).
\end{equation}
Let us denote as $\xi_{\mvec}^i$ the coefficient of the monomial term associated to the $i$-th line, $i=d+1, ..., n$, and to the vector of integers $\mvec$ such that $\mvec = [m_1, ..., m_d]^t$, where all $m_j$ are integers and $|\mvec| = m_1 + m_2 + ... + m_d = k$ is the order $k$ of the polynomials considered. Saying things differently, $\xi_{\mvec}^i$ is the coefficient of the monomial term $s_1^{m_1}s_2^{m_2}...s_d^{m_d}$, of order $k$, and each $\xivec_k^i  (\svec)$ is composed of the summations of all possible combinations of these order-$k$ monomial terms. Since Eq.~\eqref{eq:cohomnormal} is diagonal with respect to $\xi_{\mvec}^i$,   one can write, for $i=d+1, ..., n$ and for $|\mvec|=k$:
\begin{equation}
\left( \lambda_i - \mvec \lambdavec_L  \right) \xi_{\mvec}^i = \eta_{\mvec}^i,
\end{equation}
with the shorcut notation $\mvec \lambdavec_L = m_1 \lambda_1 + ... + m_d \lambda_d$.

A {\em cross resonance} occurs if there exist pairs $(\mvec,i)$ such that $\lambda_i = \mvec \lambdavec_L$. If the system has no cross resonance, then an explicit solution for the unknown coefficient $\xi_{\mvec}^i$ is found as
\begin{equation}
\xi_{\mvec}^i = \frac{\eta_{\mvec}^i}{\lambda_i - \mvec \lambdavec_L}.
\end{equation}

If a cross resonance exist, then a strong coupling exist between one slave coordinate and the set of master coordinates, and there is an obstruction in solving the normal cohomological equation. This means that a strong nonlinear coupling exist between one master and one slave coordinate such that the initial choice is not good, and the remedy consists in enlarging the number of master coordinates by considering the slave resonant modes as masters.

Following the same notations, for the tangent co-homological equation, one  arrives at, for all $i=1, ..., d$:
\begin{equation}
\lambda_i \xi_k^i  (\svec) - D\xi_k^i  (\svec) \Lambdavec_L \svec  - f_k^i (\svec )= \tilde{\eta}_k^i (\svec),
\end{equation}
with the shortcut notation $\tilde{\etavec}_k^L (\svec) = \etavec_k^L (\svec)$. Again, this last equation can be explicited in terms of the unknowns which are each of the coefficients of the monomial terms. Using the same notation one arrives at the following, for all $i=1, ..., d$
\begin{equation}
\left( \lambda_i - \mvec \lambdavec_L  \right) \xi_{\mvec}^i  - f_{\mvec}^i= \tilde{\eta}_{\mvec}^i.
\label{eq:cohomtangentindic}
\end{equation}
The pairs $(\mvec,i) \in \N^d \times {1,...,d}$, with $|\mvec| \geq 2$  such that $\lambda_i = \mvec \lambdavec_L $ creates an {\em internal resonance}, referring to a nonlinear resonance relationship between the eigenvalues of the master coordinates. When such an internal resonance exist, then there is an obstruction to the linearisation of the reduced dynamics.

Since in \eqref{eq:cohomtangentindic} two unknowns are present, namely the coefficient of the monomials of both the nonlinear change of coordinate $\xi_{\mvec}^i$ and the reduced dynamics $f_{\mvec}^i$, the solution to this equation is not unique, and can also be given even if there exist some resonances. This explains why there exist many different ways of solving the problem, thus leading to the different {\em styles} of parametrisation. The two main style of solutions are given as the {\em graph style} and the {\em normal form style}.

The graph style is simply obtained by stating that from order $k=2$, all the corrections contained in the master coordinates $\xivec_k^L (\svec)$, are vanishing:  $\xivec_k^L (\svec)=0$. This means that the master coordinates are only linearly related to the original ones. With this assumption one can then replace in Eq.\eqref{eq:cohomtangentindic} to arrive at the terms allowing one to write the reduced dynamics as, for $i=1,...d$, $|m|=k$:
\begin{equation}
f_{\mvec}^i= -\tilde{\eta}_{\mvec}^i, \quad \xi_{\mvec}^i = 0.
\end{equation}
With this choice, one recovers the classical technique promoted from center manifold theorem giving as initial guess a functional relationship (graph) between slave and master coordinates.

In the normal form style, the idea is to simplify as much as possible the reduced-order dynamics, by keeping only the resonant monomials, and discarding all other non-essential terms for the dynamical analysis. This leads to a more complex calculation, and a full nonlinear mapping between original coordinates and reduced ones, and at the end one arrives at a normal form for the reduced vector fields $\fvec$. The drawback is that calculations are a bit more involved (which is particularly true when there are numerous internal resonances to handle). The advantage is that the parametrisation is able to go over the foldings of the manifold.

To this end one solves Eq.~\eqref{eq:cohomtangentindic} following the rules (depending on the presence of internal resonance or not):
\begin{subequations}
\begin{align}
\mbox{If} \; \lambda_i \neq \mvec \lambdavec_L, \quad f_{\mvec}^i &=0, \quad \xi_{\mvec}^i =\frac{\tilde{\eta}_{\mvec}^i}{\lambda_i - \mvec \lambdavec_L},\\
\mbox{If} \; \lambda_i = \mvec \lambdavec_L, \quad f_{\mvec}^i &= - \tilde{\eta}_{\mvec}^i, \quad \xi_{\mvec}^i = 0.
\end{align}
\end{subequations}
The formulas given for this invariant manifold computation with normal form style can be extended to the case $d=n$, and one then strictly recovers the usual full normal form of the original system~\cite{Haro}.

\section{Comparison of reduced dynamics}\label{app:IMtoNF}

The aim of this appendix is to demonstrate the equivalence between the reduced dynamics given by the graph style as derived using the Shaw and Pierre approach, Eq.~\eqref{eq:romIM1mdev}, to the one obtained thanks to real normal form approach, Eq.~\eqref{eq:NFsingleNNM}. In both case reduction to a single master mode is used. Theoretically speaking, the two methods compute the same manifold and should thus provide the same dynamics. However, their formulations are differing since they are not expressed with the same variables. Let us start from the dynamics obtained using the graph style, rewritten here for the ease of reading:
\begin{equation}\label{eq:romIM1mdeva1}
\ddot{x}_m + \omega_m^2 x_m + g^m_{mm} x_m^2 + x_m \left(\underset{s\neq m}{\sum_{s=1}^N} 2\,g^m_{ms} g^s_{mm}\left[\frac{2\omega_m^2 - \omega_s^2}{\omega_s^2 (\omega_s^2 - 4 \omega_m^2)}x_m^2 +  \frac{2}{\omega_s^2 (\omega_s^2 - 4 \omega_m^2)}y_m^2 \right] \right) + h^m_{mmm} x_m^3 = 0.
\end{equation}
Using standard symmetry relationships on the quadratic coefficients, namely  $g^m_{sm} = g^m_{ms} = g^s_{mm}$, the equation can be rewritten as
\begin{equation}\label{eq:romIM1mdeva2}
\ddot{x}_m + \omega_m^2 x_m + g^m_{mm} x_m^2 + h^m_{mmm} x_m^3 +  \left(\underset{s\neq m}{\sum_{s=1}^N} (g^s_{mm})^2 \frac{2}{\omega_s^2} \left[\frac{2\omega_m^2 - \omega_s^2}{ (\omega_s^2 - 4 \omega_m^2)}x_m^3 +  \frac{2}{ (\omega_s^2 - 4 \omega_m^2)}x_m y_m^2 \right] \right)  = 0.
\end{equation}
On the other hand the reduced dynamics given by normal form writes:
\begin{equation}\label{eq:singleNF}
\ddot{R}_m + \omega_m^2 R_m + h_{mmm}^m R_m^3 + \sum_{s=1}^n 
({g}^s_{mm})^2 \frac{2}{\omega^2_{s}} \left( \frac{2 \omega^2_m - \omega^2_s}{\omega^2_s-4 \omega^2_m}\;R_m^3 + \frac{2}{\omega^2_s-4 \omega^2_m}\;\dot{R}_m^2 R_m \right)\;=\;0.
\end{equation}
Comparing term by term the two equations, one can observe two main differences: the presence of the quadratic term $g^m_{mm} x_m^2$, and the summation which excludes the term $s=m$ in the first equation. The nonlinear relationship between the modal and the normal variables, in this case of a single master coordinate, reads:
\begin{equation}\label{eq:cvmm}
x_m = R_m - \frac{1}{3\omega_m^2}g^m_{mm} R_m^2 - \frac{2}{3\omega_m^2}g^m_{mm} S_m^2 + {\mathcal O}(R_m^4,S_m^4),
\end{equation}
where the shorcut notation $S_m = \dot{R}_m$ is used. This equation is simply obtained from Eqs.~\eqref{eq:nonlinear_change_modal_compact}, assuming a single-mode motion, and replacing $a^m_{mm}$, $b^m_{mm}$ and $\gamma^m_{mm}$ with their exact analytical values given in \cite{touze03-NNM}. Note that this expansion is valid up to fourth-order, since the cubic coefficients are vanishing: $r^m_{mmm}=u^m_{mmm}=0$. Consequently no cubic terms are present in \eqref{eq:cvmm}. Replacing \eqref{eq:cvmm} in \eqref{eq:romIM1mdeva2}, and  denoting as $T=\ddot{x}_m + \omega_m^2 x_m + g^m_{mm} x_m^2$ the term that will produce extra quadratic and cubic terms, one arrives at:
\begin{align}\label{eq:devecq}
T  = &\ddot{R}_m + \omega_m^2 R_m  + \frac{2g^m_{mm}}{3} R_m^2 - \frac{2g^m_{mm}}{3\omega_m^2} S_m^2  - \frac{2g^m_{mm}}{3\omega_m^2} \left(\dot{R}_m^2 + R_m \ddot{R}_m \right) - \frac{4g^m_{mm}}{3\omega_m^4} \left(\dot{S}_m^2 + S_m \ddot{S}_m \right)\nonumber \\
  & - \frac{2(g^m_{mm})^2}{3\omega_m^2} R_m^3 -\frac{4(g^m_{mm})^2}{3\omega_m^4} R_m \dot{R}_m^2,
\end{align}
where the first line gathers linear and quadratic terms and the second the cubic terms. One can first observe that the linear terms will produce the same as those in \eqref{eq:singleNF}, a general property of identity-tangent nonlinear mapping. In order to obtain the equivalence between the two formulations, the goal is to show that the quadratic terms are vanishing. This is easily achieved by assuming that, at lower order, $\dot{R}_m = i\omega_m R_m$,  $\ddot{R}_m = -\omega_m^2 R_m$ and so on. Replacing all the combinations, one obtains that the quadratic terms exactly cancels. Also the cubic terms appearing in \eqref{eq:devecq} are exactly the one obtains from the summation in Eq.~\eqref{eq:singleNF} for $s=m$. Consequently the two equations are strictly equivalent up to the third order.

\section{Implicit static condensation}\label{app:static}

This appendix aims at giving a few more details on the static condensation method, underlining the link that the ICE method shares with explicit condensation and the role of invariant-breaking terms to produce the curvatures of the stress manifold. This section use explanations reported in~\cite{YichangICE}, and the interested reader is referred to this paper for more details.
 
Let us first underline the role of invariant-breaking terms in the construction of the stress manifold. For the sake of simplicity, let us assume modal expansion for the equations of motion, Eq.~\eqref{eq:EDOmodalproj}, and that only $x_m$ is selected as master mode. The method consists in applying a static body force of the form $\fvec_\text{e} = \beta_m \Mvec \phivec_m$ for several values of $\beta_m\in\mathbb{R}$, to compute with the FE code the corresponding displacement $\X(\beta_m)$ and the modal coordinates $x_i(\beta_m)$ by modal expansion. It results in solving the following system:
\begin{subequations}\label{eq:iceexpl0}
\begin{align}
& \omega_m^2 x_m +
\sum_{i=1}^{N} \sum_{j=1}^{N} g_{ij}^m x_i x_j +
\sum_{i=1}^{N} \sum_{j=1}^{N} \sum_{k=1}^{N}  h_{ijk}^m x_ix_jx_k
= \beta_m, \label{eq:iceexpla} \\
\forall s \neq m,\quad &\omega_s^2 x_s + \sum_{i=1}^{N} \sum_{j=1}^{N} g_{ij}^s x_i x_j +
\sum_{i=1}^{N} \sum_{j=1}^{N} \sum_{k=1}^{N}  h_{ijk}^s x_ix_jx_k = 0, \label{eq:iceexplb}
\end{align}
\end{subequations}
in which the forcing is aligned with the $m$-th eigenvector, resulting in a zero forcing of the other oscillators, for $s\neq m$. Because of this last property and the implicit function theorem, Eqs.~(\ref{eq:iceexplb}) leads to the existence of a static nonlinear relationship between the slave coordinates $x_s$ and the master one $x_m$, expressed formally as:
\begin{equation}
\label{eq:csxm}
x_s = c_s (x_m),\quad \forall s\neq m.
\end{equation}
Replacing in Eq.~\eqref{eq:iceexpla}, the reduced-order dynamics simply reads:
\begin{equation}
\ddot{x}_m + \underbrace{\omega_m^2 x_m + \sum_{i=1}^{N} \sum_{j=1}^{N} g_{ij}^m c_i (x_m) c_j(x_m) + \sum_{i=1}^{N} \sum_{j=1}^{N} \sum_{k=1}^{N}  h_{ijk}^m c_i (x_m) c_j(x_m) c_k(x_m) }_{\beta_m(x_m)}= 0,
\label{eq:ice_reddyn_gen}
\end{equation}
in which $\beta_m(x_m)$ can be identified as a polynomial in $x_m$ to obtain the reduced order dynamics at any order. From the above developments, it is clear that the dynamics of Eq.~(\ref{eq:ice_reddyn_gen}) is equivalent to the one of the full model with all the slave modal coordinates $x_s$ statically condensed  into the master dynamics.

The general explicit solution to Eq.~\eqref{eq:csxm} with closed formulation is generally out of reach such that  the $c_s$ functions are known implicitely. Nevertheless, they can be searched for as polynomial expansions and the first terms can be found. In particular, the quadratic term of the development is easy to find and is sufficient to derive the third-order dynamics, it reads~\cite{YichangICE}:
\begin{equation}
c_s (x_m) = -\frac{g^s_{mm}}{\omega_s^2} x_m^2 + \mathcal{O}(x_m^3),
\label{eq:ice_masterslave}
\end{equation}
which, substituted into Eq.~\eqref{eq:ice_reddyn_gen} with a subsequent truncation up to third order leads to Eq.~\eqref{eq:ICEdyna}.

Moreover, Eq.~\eqref{eq:iceexplb} shows that in the nonlinear terms, the important invariant-breaking terms $g^s_{mm} x_m^2$ and $h^s_{mmm} x_m^3$ have a large magnitude, and these terms create a non-zero static response for $x_s$. These couplings make the manifold depart from the linear eigensubspace to create the stress manifold.

\end{document}